\documentclass[3p]{elsarticle}
\usepackage{amsmath, amsfonts, amssymb}
\usepackage{latexsym}
\usepackage{graphicx}
\usepackage{color}
\usepackage{picins}
\usepackage{cases}

\newcommand{\trace}{\mathop{\rm Tr}\nolimits}

\newcommand{\diag}{\mathop{\rm Diag}\nolimits}

\newcommand{\rank}{\mathop{\rm Rank}\nolimits}

\newcommand{\C}{{\mathbb{C}}}

\newcommand{\R}{{\mathbb{R}}}
\newcommand{\N}{{\mathbb{N}}}

\DeclareRobustCommand\openone{\leavevmode\hbox{\small1\normalsize\kern-.33em1}}

\newcommand{\be}{\begin{equation}}
\newcommand{\ee}{\end{equation}}
\newcommand{\bea}{\begin{eqnarray}}
\newcommand{\eea}{\end{eqnarray}}
\newcommand{\beas}{\begin{eqnarray*}}
\newcommand{\eeas}{\end{eqnarray*}}

\newtheorem{theorem}{Theorem}
\newtheorem{lemma}[theorem]{Lemma}
\newtheorem{proposition}[theorem]{Proposition}

\newcommand{\ig}{{\rm Img}\ }
\newcommand{\rg}{\rank}
\newcommand{\tr}{\trace}
\newcommand{\cO}{{\mathcal O}}

\newlength{\pwdt}
\newlength{\twdt}
\newcommand{\pscl}{0.5}

\newlength{\pwdtt}
\newlength{\twdtt}
\newcommand{\psclt}{0.5}

\begin{document}
\begin{frontmatter}
\title{Impressions of convexity -- An illustration for commutator bounds}
\author{David Wenzel}
\address{Fakult\"at f\"ur Mathematik, TU Chemnitz\\
09107 Chemnitz, Germany}
\ead{david.wenzel@s2000.tu-chemnitz.de}
\author{Koenraad M.R.\ Audenaert}
\address{Mathematics Department,\\
Royal Holloway, University of London,\\
Egham TW20 0EX, United Kingdom}
\ead{koenraad.audenaert@rhul.ac.uk}
\begin{keyword}
Convexity \sep Commutator \sep Norm inequality \sep Complex Interpolation
\MSC 15A45
\end{keyword}
\begin{abstract}
We determine
the sharpest constant $C_{p,q,r}$
such that for all complex matrices $X$ and $Y$, and for
Schatten $p$-, $q$- and $r$-norms the inequality
$$
\|XY-YX\|_p\leq C_{p,q,r}\|X\|_q\|Y\|_r
$$
is valid. The main theoretical tool in our investigations is complex interpolation theory.
\end{abstract}
\end{frontmatter}

\section{Introduction}

In this paper we determine
the sharpest constant $C_{p,q,r}$
such that for all complex matrices $X$ and $Y$ the inequality
\begin{equation}\label{eqNIpqr}
\|XY-YX\|_p\leq C_{p,q,r}\|X\|_q\|Y\|_r
\end{equation}
is valid. Here, all norms are Schatten norms, i.e.
\[\|X\|_p=(\sigma_1^p+\cdots+\sigma_d^p)^{1/p}\]
with $\sigma_i$ the decreasingly ordered singular values $\sigma_1\geq\ldots\geq\sigma_d\geq 0$ of $X$.

This question is a straightforward continuation of a line of investigation about analogous inequalities
considered previously with special choices for the norm indices $p,q$ and $r$.
For instance, in \cite{Aud} one of us raised the conjecture
that in the case $q=p$ one has
\begin{equation}\label{eqCppr}
C_{p,p,r}=2^{\max\{1/p,1-1/p,1-1/r\}}.
\end{equation}
We want to show the validity of this conjecture and carry over the developed ideas
to the general situation. We will also take a closer look at the cases of equality in (\ref{eqNIpqr}),
studied previously for $p=q=r=2$ in \cite{CVW}.

\medskip

The main technique used in this paper is complex interpolation \textit{a la} Riesz-Thorin,
applied in a rather intricate way to the problem at hand.
To achieve optimal clarity, the exposition will partially leave the usual format,
with two effects. While certain steps in the proofs
later turn out to be redundant, we have chosen to keep them in because of their use in the development
of the complete proof and their importance in obtaining a better understanding of
what is going on behind the scenes. Secondly, some
parts are not following the usual structure and should be
understood as a written presentation that will guide the
reader through our thoughts.

\subsection{Notations}\label{sec1.1}

We will use some abbreviations in formulas:
the Lie bracket $[X,Y]=XY-YX$ for the commutator,
$\sigma(X)$ for the vector of singular values of $X$,
$\tr X$ for its trace, $X^T$ for its transpose and $X^*$ for its adjoint.
Moreover, $\cO$ will denote a zero matrix of appropriate size,
$I_n$ a $n\times n$ identity matrix and $A\oplus B=\diag(A,B)$ will be written
for the construction of block diagonal matrices.
For any norm index $p\in[1,\infty]$, $p'$ denotes the conjugate index of $p$, i.e.\ the number
$p'\in[1,\infty]$ satisfying $\frac{1}{p}+\frac{1}{p'}=1$. As is well-known,
the Schatten $p'$ norm is the dual norm of the Schatten $p$ norm.
Note that we took the formal equality $\|X\|_p=\|\sigma(X)\|_p$ as sufficient reason for
denoting the usual $\ell_p$ norm of a vector also by $\|\cdot\|_p$.

\subsection{Illustrations}\label{sec1.2}

Throughout the paper our proofs will be of a very pictorial nature,
because there are so many special cases to be considered, and it so happens
that these cases can be presented graphically in a very clear way.
We hope that this will allow the reader to gain a better understanding of the several
steps and at the same time quickly obtain an overall view of the whole proof.

As stated before, the topic of this paper is finding the best constant $C_{p,q,r}$ in (\ref{eqNIpqr}),
where $p$, $q$ and $r$ are norm indices, $1\le p,q,r$. The triplet of values $(p,q,r)$ can be depicted as
a point in $\R^3$, or more precisely in $[1,\infty]^3$.
The proofs of our theorems require subdividing this infinite cube in several regions, and rather than
just define these regions in the usual way (with equalities and inequalities), we will
augment every definition with a graphical illustration, of points and regions in $\R^3$
or $\R^2$ (when we restrict to the case $p=q$),
where every real axis corresponds to one of these norm indices.
In addition we'll use these pictures to display many other quantities that are important in the proofs, but
that will become clear later on.

\medskip

Of course, we need some device
to portray the whole real line or even only the semi-bounded
interval $[1,\infty]$ in a finite space. So we need to cheat a little bit and
we will distort reality by mapping norm indices $p\in[1,\infty]$
to positions in the image given by the reciprocal of the conjugate index $\frac{1}{p'}$.

\begin{figure}[b]
\centering\includegraphics[scale=0.5]{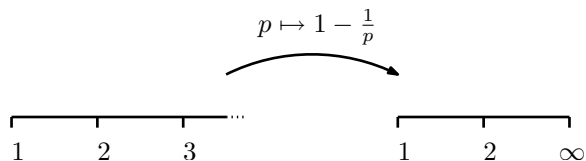}

\caption{The scaling of norm indices for 1D imaging purposes.}\label{figImg}
\end{figure}

Applying this mapping
\[\ig\!: [1,\infty]\rightarrow\R,\quad p\mapsto 1-\frac{1}{p}\]
in illustrations has several advantages (see Figure \ref{figImg}).
Firstly, we obtain finite pictures as $[1,\infty]$ is mapped onto $[0,1]$.
Moreover, the unreachably far away index $p=\infty$ becomes the handy
point $\ig \infty=1$.
The mapping preserves the order of the norm indices, i.e.
\[p<q\quad\Rightarrow\quad \ig p<\ig q.\]
So, we are given just an appropriate scaling and the smallest possible
index $p=1$ is of course the left-most point in the images.
Last but not least, the index $p=2$ is mapped
exactly to the middle of the line segment, befitting its special role
as the only self-conjugate index.

\begin{figure}
\centering\includegraphics[scale=0.5]{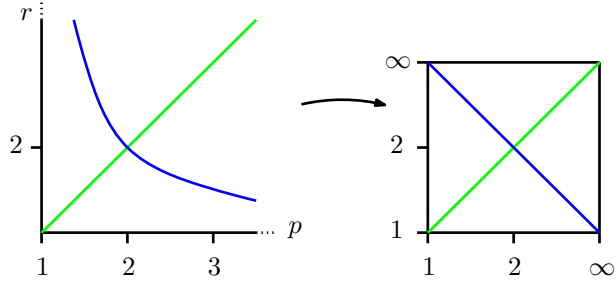}

\caption{The curves $r=p$ (green) and $r=p'$ (blue) in the original and the 2D scaled setting.}\label{figImg2}
\end{figure}

\medskip
As the first object of interest (\ref{eqCppr}) involves two norm indices $p$ and $r$
we are going to use two-dimensional images by applying the scaling
function twice independently:
\[\ig\!^2 : [1,\infty]\times[1,\infty]\rightarrow\R^2,\quad (p,r)\mapsto (\ig p, \ig r).\]
The result is a finite square whose center corresponds to the well known
special case $p=r=2$ that was proved in \cite{BW2}.
There are some other nice side effects. The points satisfying $r=p$ still form a straight
line in the graphics. Moreover, the curve $r=p'=(1-1/p)^{-1}$ is mapped to the square's
other diagonal (see Figure \ref{figImg2}).

\medskip
Later on, when we study (\ref{eqNIpqr}) in full generality, we will use this same image scaling to three
dimensional pictures:
\[\ig\!^3 : [1,\infty]^3\rightarrow\R^3,\quad (p,q,r)\mapsto (\ig p, \ig q, \ig r).\]
There are again several curves that have lines as images. Furthermore, we will encounter some surfaces
that are conveniently mapped to planes.

\subsection{Basics on norm interpolation}\label{sec1.3}

We want to briefly introduce a concept that is a key to our proofs and will
be used extensively in the remainder of the paper. More detailed explanations
and additional applications can be found in \cite{G}.

\medskip
In 1926 M.\ Riesz established a theorem that allows to interpolate between two
inequalities involving the usual $\ell_p$ vector norms. Stated in our notations:
\begin{theorem}[Riesz-Thorin]\label{thmRT}
Let $1\leq p_0\leq p_1\leq\infty$ and $1\leq q_0\leq q_1\leq\infty$ be given such that
\begin{equation}\label{eqRTrealpq}
q_1\leq p_1\quad {\sl and}\quad q_2\leq p_2.
\end{equation}
If for a linear operator
\begin{equation}\label{eqRTop}
T:\R^k\rightarrow\R^n
\end{equation}
there are $M_0,M_1>0$
such that
\begin{equation}\label{eqRTbaseineq}
\|Tx\|_{p_0}\leq M_0\|x\|_{q_0}\quad {\sl and}\quad \|Tx\|_{p_1}\leq M_1\|x\|_{q_1}
\end{equation}
for all arguments $x$, then for any $\theta\in[0,1]$ and every vector $x$ the inequality
\begin{equation}\label{eqRTintineq}
\|Tx\|_{p}\leq M_0^{1-\theta}M_1^\theta\|x\|_{q}
\end{equation}
holds with $p\in[p_0,p_1], q\in[q_0,q_1]$ defined by
\begin{equation}\label{eqRTconvcomb}
\frac{1}{p}=\frac{1-\theta}{p_0}+\frac{\theta}{p_1}\quad {\sl and}\quad \frac{1}{q}=\frac{1-\theta}{q_0}+\frac{\theta}{q_1}.
\end{equation}
\end{theorem}

The theorem was enshrined in the fundamental methods of analysis,
when Riesz' student G.O.\ Thorin extended the theorem to complex arguments and operators,
obtaining an analogon of Theorem \ref{thmRT} with
(\ref{eqRTop}) replaced by
\[T:\C^k\rightarrow\C^n.\]
His proof, based on an ingenious use of Hadamard's three line lemma from the theory of analytic functions,
reveals the surprising fact that the condition (\ref{eqRTrealpq}) is no longer necessary in the complex case
(essentially because condition (\ref{eqRTbaseineq}) must now hold for all \textit{complex} vectors); an
assertion that is completely wrong in the real case!

Afterwards, the result was extended to operators $T$ defined on subspaces and, by help
of density arguments, to operators acting on infinite-dimensional spaces, in particular the $L^p$-spaces.
Moreover, it was shown that, if $x$ and $Tx$ are matrices, the underlying norms may be replaced by
their Schatten type analogues. This holds due to a general equivalence between
sequence spaces and the corresponding Schatten classes as far as interpolation is concerned \cite{arazy}.

Recently, in \cite{W}, one of us restated the theorem in terms of a special structure, the tensor product of argument
vectors, with the purpose of investigating (\ref{eqNIpqr}) with $p=q=r$.
Although formulated in a more specific way in that paper,
its wider validity was noted. Indeed, one can replace (\ref{eqRTop}) by
\[T:\C^{(i+k)\times(j+l)}\rightarrow\C^{m\times n}\]
and substitute
\[x=X\otimes Y\]
with matrices $X\in\C^{i\times j}, Y\in\C^{k\times l}$.
That is, we are given a linear operator on the whole set of matrices, but
only apply it to arguments that are tensor products
(also called Kronecker product for matrices).

The proof is an adaption of Thorin's proof as presented in \cite{G}, combined with
the fact that the generated simple functions (actually vectors in the finite-dimensional
case) respect the tensor structure of the arguments.
As for the original theorem, $X$ and $Y$ may be taken from subspaces
of $\C^{i\times j}$ or $\C^{k\times l}$, respectively.

\medskip
In the aftermath of the WATIE 2009 conference we learned about the multilinear
version of the Riesz-Thorin theorem. In the multilinear case
(\ref{eqRTop}) is replaced by the multilinear operator
\[T:\C^{k_1}\times\cdots\times\C^{k_m}\rightarrow\C^n,\]
(\ref{eqRTbaseineq}) and (\ref{eqRTintineq}) by inequalities like
\[\|T(x^{(1)}, ..., x^{(m)})\|_{p_\theta} \leq M_\theta \|x^{(1)}\|_{q_\theta^{(1)}}\cdots\|x^{(m)}\|_{q_\theta^{(m)}}\]
and (\ref{eqRTconvcomb}) then consists of $m$ inequalities for fixing $q^{(j)}$ \cite{BS}.

Closer inspection revealed that the statement is actually equivalent to the
usual interpolation but applied to tensor products,
owing to the property $\|X\otimes Y\|_p=\|X\|_p\|Y\|_p$ of Schatten norms.
Later on, we will see that the multilinear interpretation is too comprehensive
for our needs, whereas the original interpolation theorem and its
diagonal extension via tensor products serve their purpose very well.
Be sure to read the acknowledgement for some more insights.

\medskip
Our scaling function (Section \ref{sec1.2}) is especially convenient for picturing certain salient aspects
related to norm interpolation. The Riesz-Thorin theorem, in particular (\ref{eqRTconvcomb}), tells us that
in terms of reciprocals, the interpolated
index $\frac{1}{p}$ is a convex combination of the base indices
$\frac{1}{p_0}$ and $\frac{1}{p_1}$. This is the reason why the Riesz-Thorin theorem is sometimes
called a convexity theorem.

If the norms of argument and target space are different (i.e.\ $p_i\neq q_i$),
we need to consider a joint convex combination of the index reciprocals.
Due to the way $\ig$ is defined in the images of this paper $\ig\!^2 (p,q)$
conveniently lies on a straight line between $\ig\!^2 (p_0,q_0)$ and $\ig\!^2 (p_1,q_1)$.

We call the points $(p_i,q_i)$ \textit{interpolation base points} and all points $(p,q)$
subject to (\ref{eqRTconvcomb}) \textit{interpolants}. We carry over this nomenclature to
the associated inequalities and to the images of the points in our pictures.

\begin{figure}[b]
\vspace*{20pt}
  \centering\includegraphics[scale=0.75]{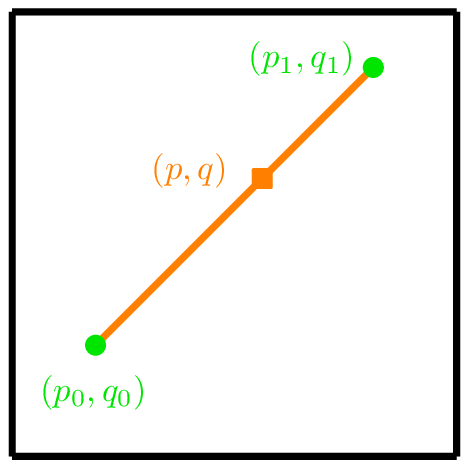} \quad\quad\quad \includegraphics[scale=0.75]{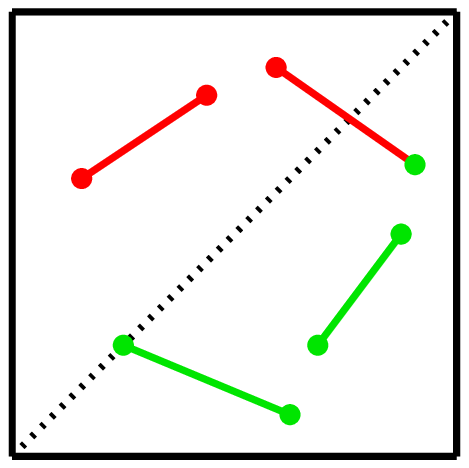}

  \caption{Left: Illustrating interpolation base points $(p_j,q_j)$ and interpolants $(p,q)$ inbetween;
  Right: Possible choices for base points (green) and obtained interpolants as well as
  points yielding no statement (red) for the real case.}\label{figIntCR}
\end{figure}

Note that for real interpolation both base points are necessarily located in the lower
triangle determined by the main diagonal $q=p$ (Figure \ref{figIntCR}).

\subsection{Overview}
For the sake of clarity, in Section \ref{sec2},
we start with treating the original and simpler conjectured inequality (\ref{eqCppr}), about the
constant $C_{p,p,r}$. Since only two parameters enter the treatment, the pictures are 2-dimensional.
In Section \ref{sec3}, the approach used in Section \ref{sec2} is generalised to treat as much
of the general 3-parameter problem as possible. In the course of this process, we will encounter
a number of parameter regions that could not, as yet, be treated using the interpolation methods
applied in Section \ref{sec2}. To overcome this hurdle, two things are needed.
Firstly, the value of $C_{p,q,r}$ in certain extremal points of parameter space
must be established. This is done in Section \ref{sec4} using a combination of
basic linear algebra methods and esoteric knowledge about certain magical symbols.
Secondly, the remaining areas of parameter space have to be covered, and this is done in Section
\ref{sec4.3} using more advanced versions of Riesz-Thorin interpolation. Thus, the proof of our main
theorem is finished at that point.
We hasten to add that for certain
areas in parameter space the $C_{p,q,r}$ constant depends on the dimension $d$ of the matrices.
Furthermore, in some instances the interpolation method did not yield the sharpest possible bound.
In Section \ref{sec5}, the cases of equality are considered, and we wrap up with a conclusion
(Section \ref{sec6}) and a list of recommended readings.

\section{The original conjecture and its proof}\label{sec2}

This section is dedicated to the derivation of (\ref{eqCppr}), which is the following theorem.

\begin{theorem}\label{thmCppr}
With the notations of equation (\ref{eqNIpqr}),
\[C_{p,p,r}=\max\left\{2^{1/p},2^{1-1/p},2^{1-1/r}\right\}.\]
\end{theorem}

This is the original conjecture stated in \cite{Aud}.
The proof we give here is somewhat longer than what could have been,
but in this way it clearly demonstrates the power and applicability of interpolation.
Near the end of this section, the reader will notice that the proof may be shortened a bit.

\subsection{The claim and some special situations}\label{sec2.1}

To begin with, we ensure that the value claimed for $C_{p,p,r}$ can be attained.
For this, take a look at the examples
\begin{equation}\label{eqspec2}
X=\left(\begin{array}{cc}
1 & 0 \\
0 & -1
\end{array}\right),\quad
Y=\left(\begin{array}{cc}
0 & 1 \\
0 & 0
\end{array}\right),\quad
XY-YX=\left(\begin{array}{cc}
0 & 2 \\
0 & 0
\end{array}\right),
\end{equation}
yielding the quotients $\frac{\|XY-YX\|_p}{\|X\|_p\|Y\|_r}=2^{1-1/p}$
as well as $\frac{\|XY-YX\|_p}{\|Y\|_p\|X\|_r}=2^{1-1/r}$ and
\begin{equation}\label{eqspec1}
X=\left(\begin{array}{cc}
0 & 1 \\
0 & 0
\end{array}\right),\quad
Y=\left(\begin{array}{cc}
0 & 0 \\
1 & 0
\end{array}\right),\quad
XY-YX=\left(\begin{array}{cc}
1 & 0 \\
0 & -1
\end{array}\right),
\end{equation}
giving the value $2^{1/p}$.
Hence, the constant $C_{p,p,r}$ cannot be smaller than asserted.

\begin{figure}[p]
\centering\includegraphics[scale=0.5]{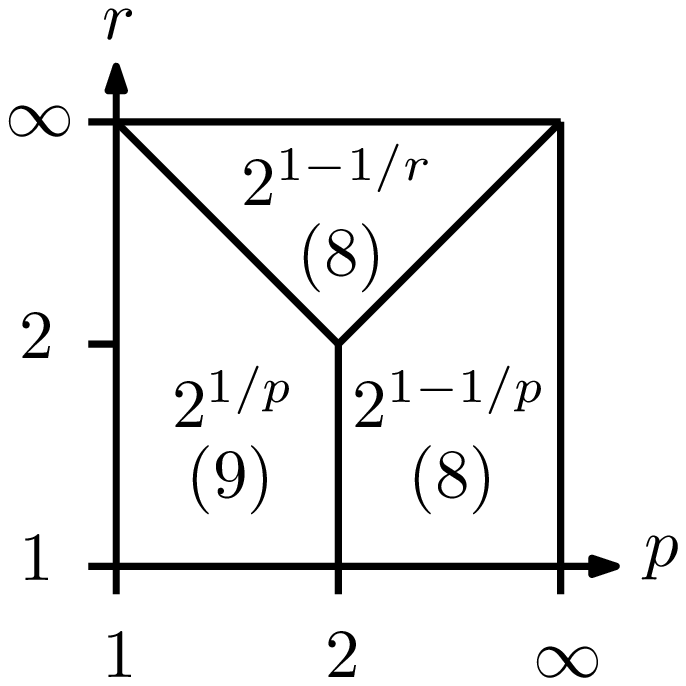}\includegraphics[scale=0.5]{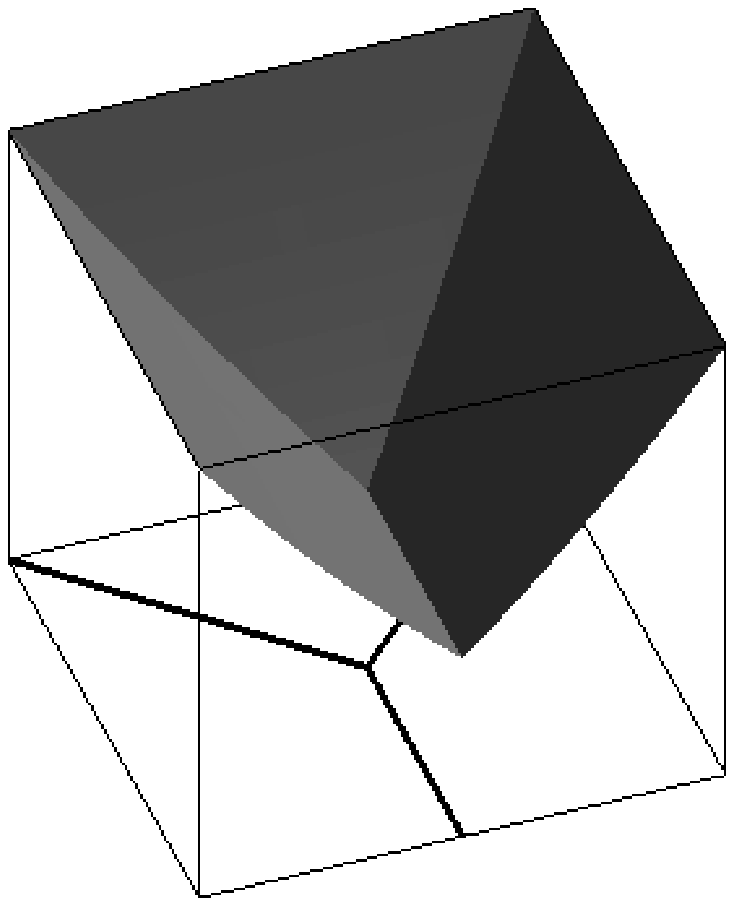}

\caption{The three segments where the constant $C_{p,p,r}$ takes on different values
according to Theorem \ref{thmCppr} and the reference to
examples achieving equality (left), and the graph of $C_{p,p,r}$ as a function of $(p,r)$
(right).}\label{figConjppr}
\end{figure}

\medskip
Because of the appearance of a maximum, over the three terms as stated in Theorem \ref{thmCppr},
the set of all pairs $(p,r)$ of norm indices is naturally subdivided into three segments:
\begin{align*}
& 1\leq p\leq 2\quad \wedge\quad r\leq p'  & \Rightarrow\quad & C_{p,p,r} = 2^{1/p}  \\
& 2\leq p\leq \infty\quad \wedge\quad r\leq p  & \Rightarrow\quad & C_{p,p,r} = 2^{1-1/p}  \\
& r\geq p'\quad \wedge\quad r\geq p  & \Rightarrow\quad & C_{p,p,r} = 2^{1-1/r}.
\end{align*}
That this is an equivalent statement is easily verified analytically
and is illustrated in Figure \ref{figConjppr}.

\medskip

The conjecture is already known to be true in some special cases, namely
\begin{itemize}
\item $p=r=2$; this case is the origin of the investigations and was shown in full generality in \cite{BW2};

\item $p=r\in[1,\infty]$, proven in \cite{W};

\item $p=2,r\in[1,\infty]$, proven in \cite{Aud}.
\end{itemize}

\begin{figure}[p]
\centering\includegraphics[scale=0.5]{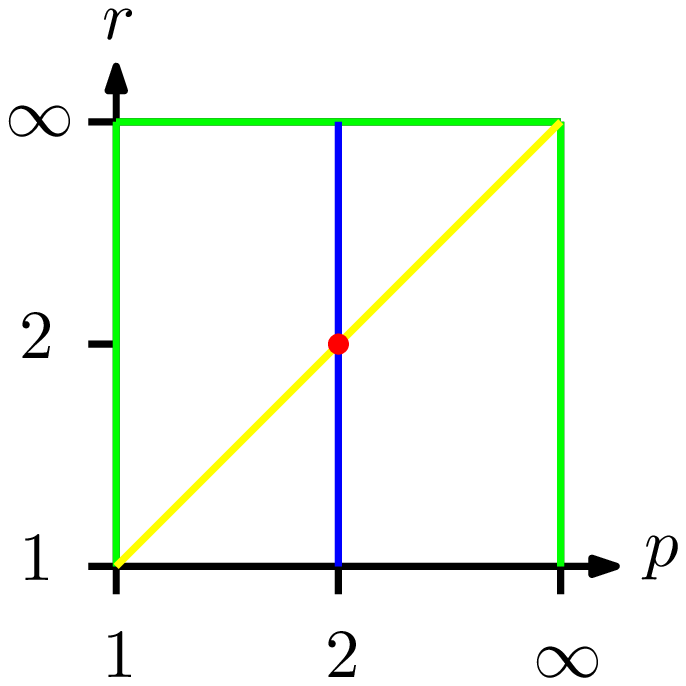}\includegraphics[scale=0.5]{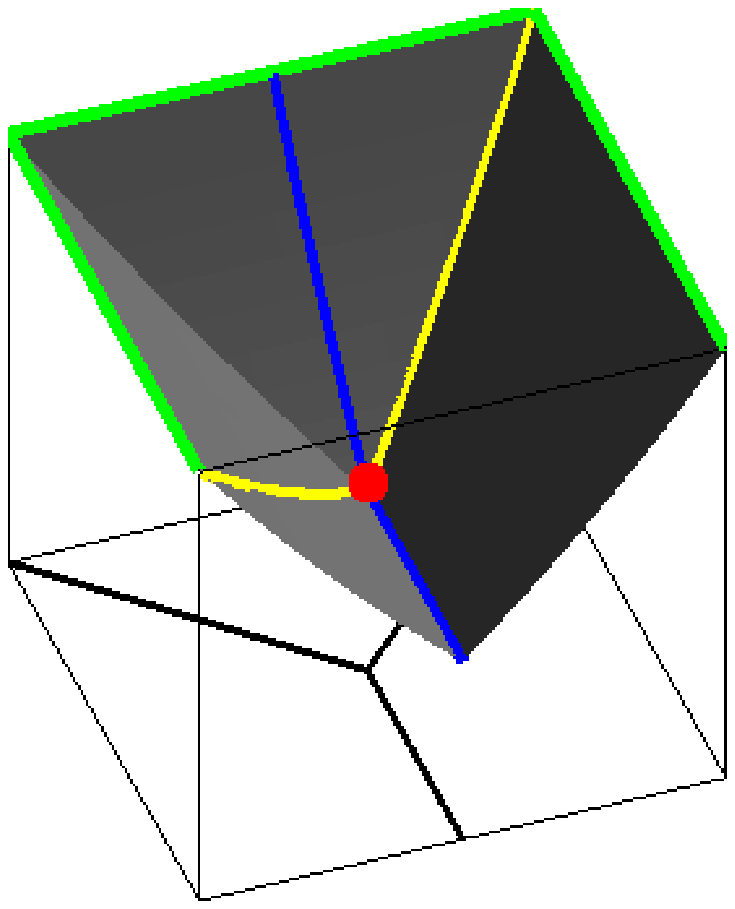}

\caption{Known special cases for $C_{p,p,r}$: red: $p=r=2$,
yellow: $p=r\in[1,\infty]$,
blue: $p=2,r\in[1,\infty]$ and the trivial cases in
green.}\label{figConjknow}
\end{figure}

The conjecture holds trivially for
\begin{itemize}
\item $p=1$, \\
as $\|XY\|_1\leq \|X\|_1\|Y\|_\infty \leq\|X\|_1\|Y\|_r$ and
the triangle inequality $\|XY-YX\|_1\leq\|XY\|_1+\|YX\|_1$
together with (\ref{eqspec1}) give $C_{1,1,r}=2$;

\item $r=\infty$, \\
since also $\|XY\|_p\leq \|X\|_p\|Y\|_\infty$ holds and (\ref{eqspec2})
realizes equality;

\item $p=\infty$, \\
because of equally simple conclusions.
\end{itemize}
These pairs $(p,r)$ and their corresponding constants are depicted in Figure \ref{figConjknow}.

For all 2D images depicting $(p,r)$ of Theorem \ref{thmCppr} (as in the left of Figures \ref{figConjppr} and
\ref{figConjknow}) we will subsequently omit axis labels to avoid unnecessary information overflow.


\subsection{A re-interpretation of known cases}

First we reconsider the case $p=2,r\in[1,\infty]$, but from a different point of view. The validity was
obtained by one of us as a consequence of an even stronger inequality \cite{Aud}. We want
to deduce the value of $C_{2,2,r}$ in a different way,
show-casing the two major techniques (complex interpolation and norm index monotonicity, see below)
we will repeatedly use in the rest of the paper.

We will also demonstrate the strong link between the promised pictures and the associated
argumentation and formulas.

\bigskip

\parpic[l]{\includegraphics[scale=\pscl]{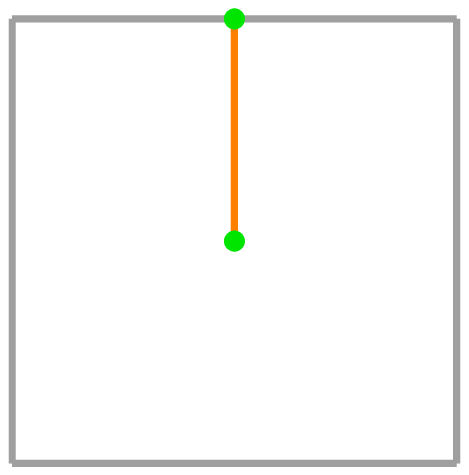}}\picskip{7}
We know the values of $C_{2,2,2}$ and $C_{2,2,\infty}$ from the inequalities
\[\|[X,Y]\|_2\leq \sqrt{2} \|X\|_2\|Y\|_2 \,\,\,{\rm and}\,\,\, \|[X,Y]\|_2\leq 2 \|X\|_2\|Y\|_\infty\]
for all $d\times d$ matrices $X$ and $Y$.
The respective pairs of parameter values $(2,2)$ and $(2,\infty)$ are
represented by the green points in the picture at the left.

Now fix an arbitrary $X$ with $\|X\|_2=1$ and consider the commutator as a linear operator
\[K_X:\C^{d\times d}\rightarrow\C^{d\times d}, Y\mapsto XY-YX.\]
Clearly, we have
\[\|K_X(Y)\|_2\leq \sqrt{2} \|Y\|_2 \quad{\rm and}\quad \|K_X(Y)\|_2\leq 2 \|Y\|_\infty\]
for any $Y$.
As these correspond to the premises (\ref{eqRTbaseineq}) of the Riesz-Thorin theorem, in its usual form,
(see Theorem \ref{thmRT} and the comments on generalization thereafter)
we endeavour to apply this theorem
for $p=2,r\in(2,\infty)$
(the points on the orange line in the picture).
For this we require the validity of (\ref{eqRTconvcomb}), which is, in our case:
for any $\theta\in(0,1)$
\[\frac{1}{2}=\frac{1-\theta}{2}+\frac{\theta}{2}\quad{\rm and}\quad
\frac{1}{r}=\frac{1-\theta}{2}+\frac{\theta}{\infty}.\]
As the first equality is trivially true we have that the parameter
\[\theta = 1-\frac{2}{r}\]
is in one-to-one correspondence to all possible interpolants $(2,r)$.
Consequently, we obtain inequality (\ref{eqRTintineq}), that is
\[\|K_X(Y)\|_2\leq \sqrt{2}^{1-\theta} 2^\theta \|Y\|_r\]
or equivalently
\[\|XY-YX\|_2\leq 2^{1-1/r} \|X\|_2\|Y\|_r\]
and hence
$C_{2,2,r}\leq 2^{1-1/r}$, as required. Note that assuming $X$ to be normalised
incurs no loss of generality.

\medskip
\parpic[l]{\includegraphics[scale=\pscl]{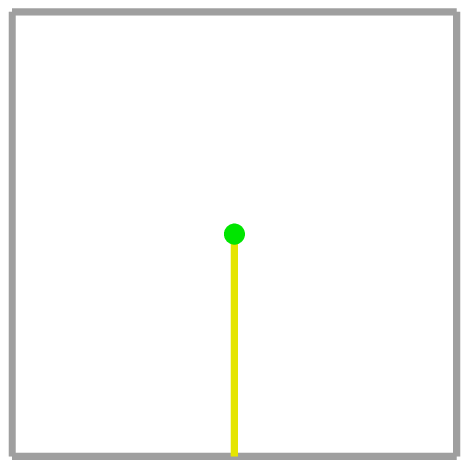}}\picskip{8}
For the remaining case $r\in[1,2)$ we can use a simpler concept,
which we would like to call \textit{norm index monotonicity}, or just monotonicity for short. By this
we mean the well-known relation
\[\|A\|_p\leq \|A\|_q\quad {\rm for\ any\ } p\geq q\]
and arbitrary matrices $A$. This procedure could be regarded as an interpolation
with only one base.

In this manner we obtain directly from the knowledge of $C_{2,2,2}=\sqrt{2}$ that
\[\|XY-YX\|_2\leq \sqrt{2} \|X\|_2\|Y\|_2 \leq \sqrt{2} \|X\|_2\|Y\|_r,\]
which gives $C_{2,2,r}\leq\sqrt{2}$ for $p<2$ (points on the yellow line).

For both cases, $r>2$ and $1\le r<2$, the
proof is now easily completed by providing an example of two matrices that achieve equality,
as we have already done in Section \ref{sec2.1}.

\medskip

We will keep on the arrangement for picturing known base points green and indicating an
interpolation process by an orange line and the use of the monotonicity argument by a yellow line.

\subsection{Towards a full proof}\label{sec2.3}

In this section we give an intuitive overview of the proof
of Theorem \ref{thmCppr}, but on the other hand also provide the necessary details for
more demanding readers. To accommodate both audiences we have adopted an unusual style that may be called
a scientific graphic novel. In an attempt to avoid boring the reader too much,
the level of detail will be reduced in due course
when coming across cases that are similar to already covered ones.

Roughly speaking, the proof can be subdivided in four parts, each part corresponding to one of the four
quadrants of the parameter space: the lower left quadrant, corresponding to $p,r\le 2$,
the lower right, $p\ge 2$, $r\le 2$,
the upper right $p,r\ge2$, and the upper left quadrant $p\le 2$, $r\ge 2$.
We begin with the lower left quadrant.

\parpic[l]{\includegraphics[scale=\pscl]{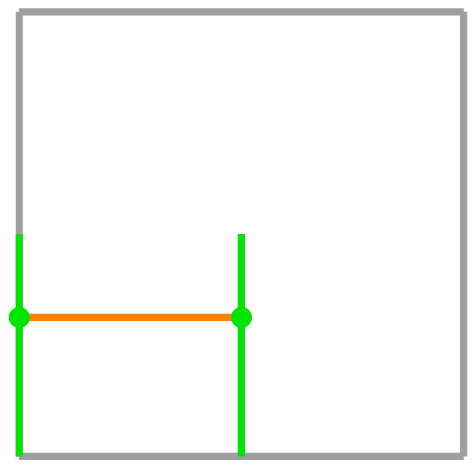}}\picskip{8}
The conjecture can easily be shown to be true for $p,r\leq 2$ by ordinary Riesz-Thorin interpolation.
For this fix $r\in[1,2]$ arbitrarily.\\
We obtain two points on the green lines in the picture, for which we have
\[
\|[X,Y]\|_1\leq 2 \|X\|_1\|Y\|_r\quad{\rm  and}\quad \|[X,Y]\|_2\leq \sqrt{2} \|X\|_2\|Y\|_r.
\]
Regard the commutator as a map $K_Y(X)=[X,Y]$ with some fixed $Y$ with $\|Y\|_r=1$. So,
\[
\|K_Y(X)\|_1\leq 2 \|X\|_1\quad{\rm  and}\quad \|K_Y(X)\|_2\leq \sqrt{2} \|X\|_2.
\]
As the norm indices of original and target space coincide for both inequalities,
we need to satisfy
\[\frac{1}{p}=\frac{1-\theta}{1}+\frac{\theta}{2}\]
twice. Hence,
$\theta=2-2/p$ and from the Riesz-Thorin theorem we immediately get
$C_{p,p,r}\leq 2^{1-\theta} \sqrt{2}^\theta = 2^{1/p}$.

\medskip
\parpic[l]{\includegraphics[scale=\pscl]{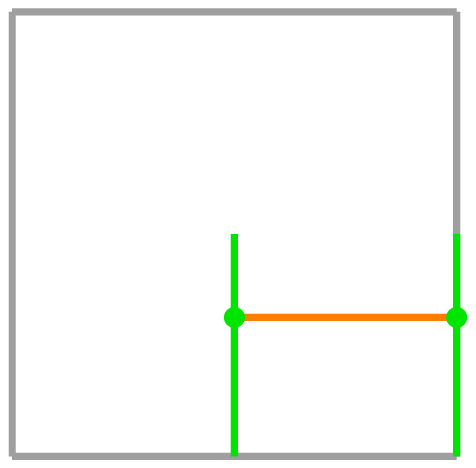}}\picskip{6}
Interpolation also works in the case $p\geq2,r\leq 2$. Again,
fix $r\in[1,2]$. Here, we have
\[
\|[X,Y]\|_2\leq \sqrt{2} \|X\|_2\|Y\|_r {\rm \ \ and\ \ } \|[X,Y]\|_\infty\leq 2 \|X\|_\infty\|Y\|_r.
\]
Interpolation of $K_Y$ now requires
\[\frac{1}{p}=\frac{1-\theta}{2}+\frac{\theta}{\infty}\]
which amounts to $\theta=1-2/p$
and yields
$C_{p,p,r}\leq \sqrt{2}^{1-\theta} 2^\theta = 2^{1-1/p}$.

Note that applicability of Theorem \ref{thmRT} comes from fixing the variable $Y$ (for given norm index $r$) or $X$
(for given $p$, as in the last subsection),
which is expressed by a vertical or horizontal line in our graphics
of norm index pairs $(p,r)$. We remark that jointly interpolating $p$ and $r$
for bivariate inequalities is not supported by the original theorem.
Hence, slanted (non-horizontal, non-vertical) lines for interpolation are forbidden here.

\begin{figure}
\centering\includegraphics[scale=0.4]{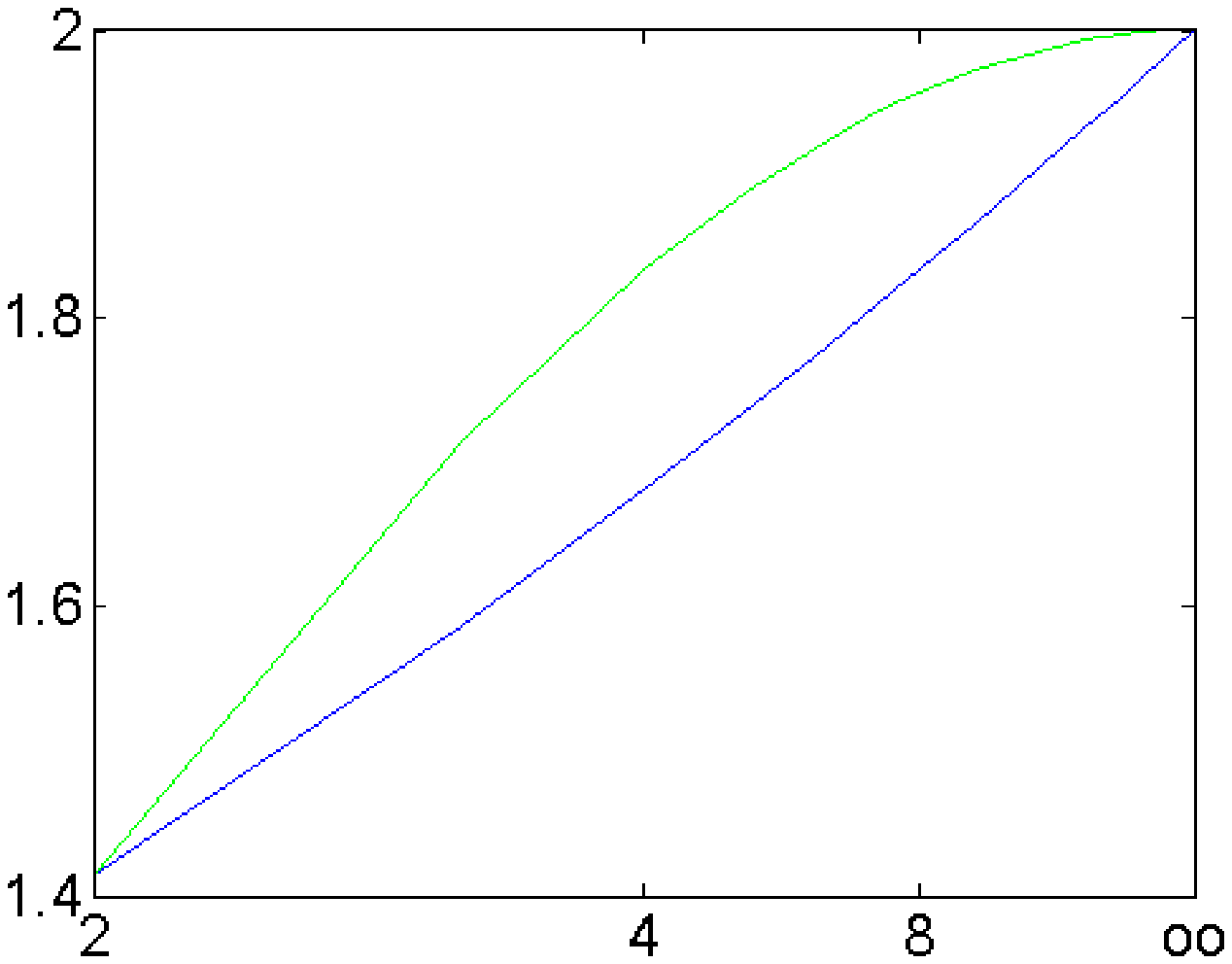} \includegraphics[scale=0.4]{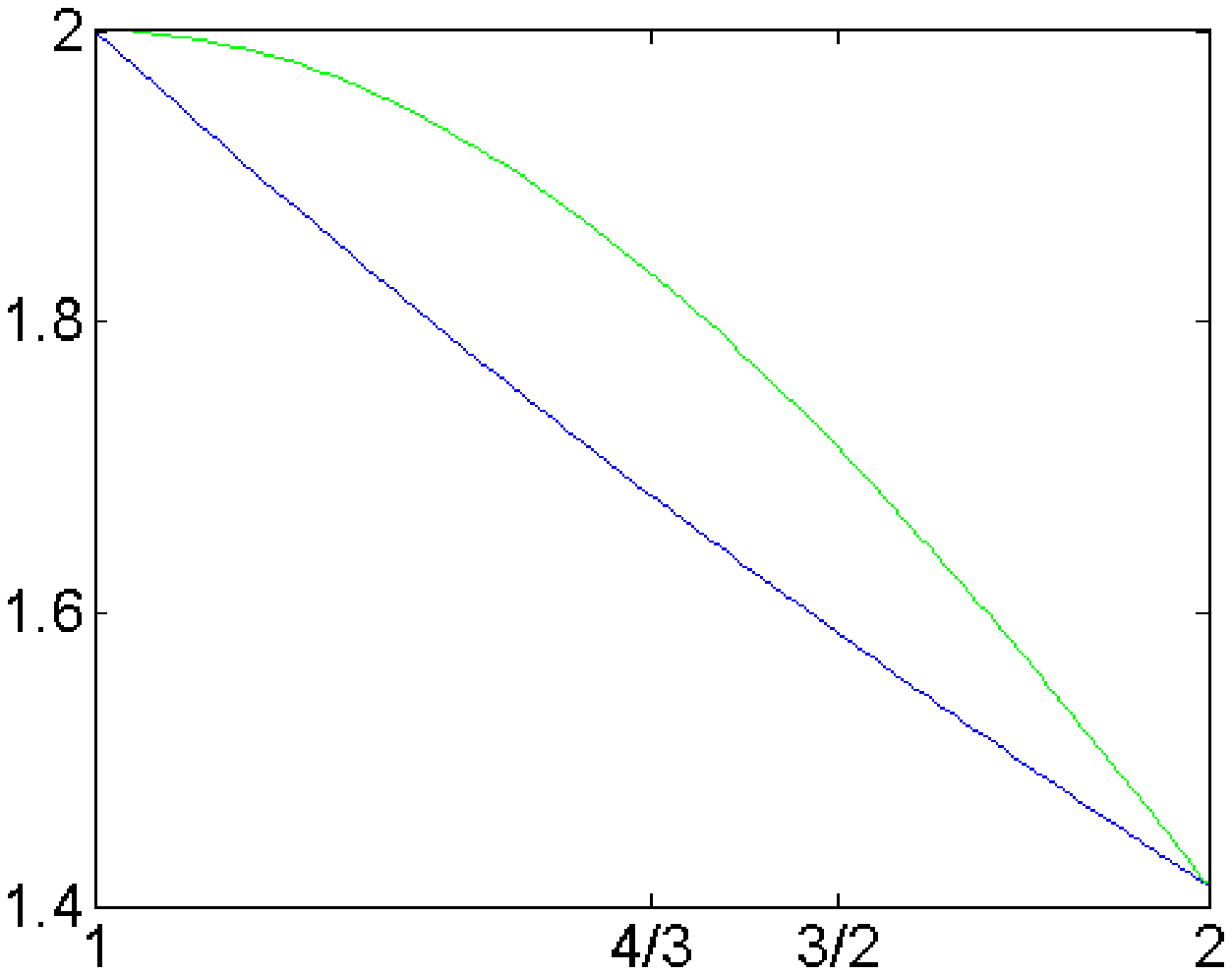}

\caption{Claimed values (blue) for $C_{p,p,r}$ and interpolation bounds (green) for the lines
$r=p$ in the upper right quadrant of parameter space (left)
and $r=p'$ in the upper left quadrant (right)\vspace*{20pt}}\label{figDiffQ34}
\end{figure}

\medskip
\parpic[l]{\includegraphics[scale=\pscl]{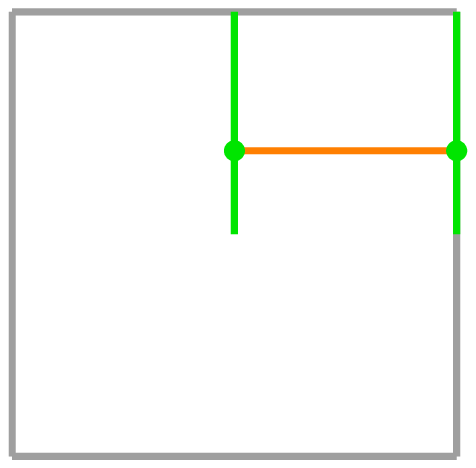} \includegraphics[scale=\pscl]{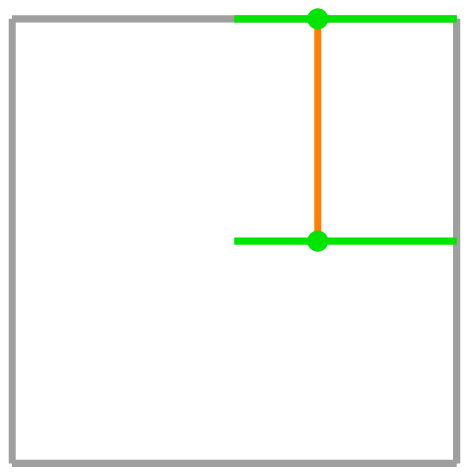}}\picskip{6}
Now, having covered half of the proof, interpolation will not work for
$p\geq 2, r\in(2,\infty)$ as it did in the previous cases.
Regardless whether we interpolate $K_Y$ (left) or $K_X$ (right), i.e.\ fixing $r$ or $p$, the
obtained bound always gives a larger value than the one claimed in Theorem \ref{thmCppr}:
$2^{1-2/rp}\geq\max\{2^{1-1/p},2^{1-1/r}\}$.
See Figure \ref{figDiffQ34} for an illustration of the difference.

\medskip
\parpic[l]{\includegraphics[scale=\pscl]{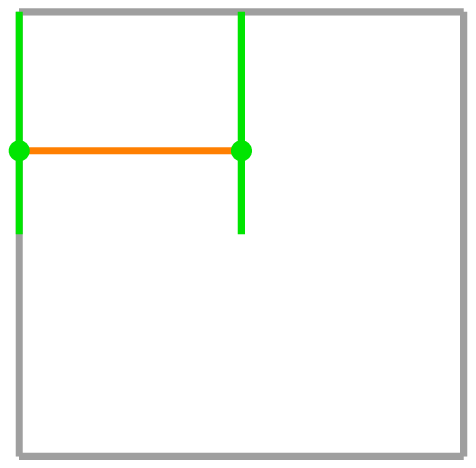} \includegraphics[scale=\pscl]{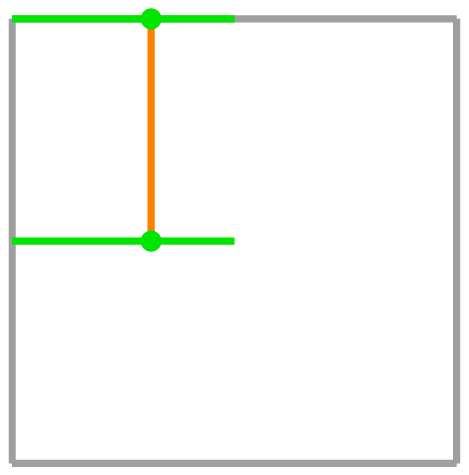}}\picskip{6}
Analogously, for $p\leq 2, r\geq 2$, interpolation does not yield the desired bound either, but gives
the larger value of $2^{1-2/r+2/rp}$.\\
This actually was to be expected, since interpolation produces smooth bounds, while the claimed constant
is not smooth as function of $p$ and $r$.
To wit, the graph of the constant exhibits cusps at the lines $r=p$ and $r=p'$, lines that are intersected
by the current directions of interpolation; for this reason we call these lines \textit{cusp lines}.

Nevertheless, we can obtain sharp values for $C_{p,p,r}$ in a more complicated, two-step interpolation process,
combining one of the more advanced
versions of the Riesz-Thorin method with its ordinary version,
and carefully choosing the right step at the right time.

\medskip
\parpic[l]{\includegraphics[scale=\pscl]{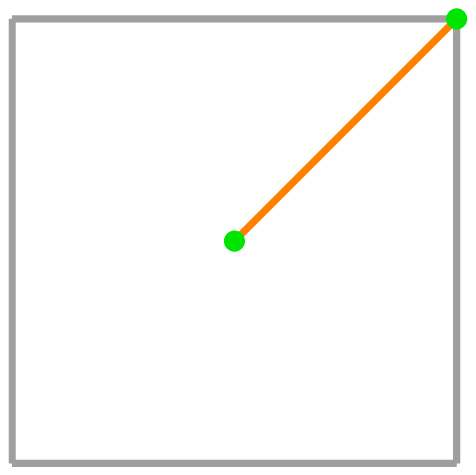}}\picskip{6}
The value of $C_{p,p,r}$ along one of the cusp lines, namely the main diagonal $r=p$,
can be obtained with help of the tensor structure interpolation
mentioned in Section \ref{sec1.3}, between $C_{1,1,1}, C_{2,2,2}$ and $C_{\infty,\infty,\infty}$.
By applying the usual interpolation
statements on these special arguments, interpolation along diagonal lines becomes possible.
This has already been done in \cite{W}, with the result
$C_{p,p,p} = 2^{1-1/p}$ for $p>2$ (and $C_{p,p,p} = 2^{1/p}$ for $p<2$).\\
Note that for $p<2$
this result has also been obtained  in the above, but in a single step, using ordinary
Riesz-Thorin interpolation.
This shows that the more complicated approach does not always lead to sharper bounds.

\medskip
\parpic[l]{\includegraphics[scale=\pscl]{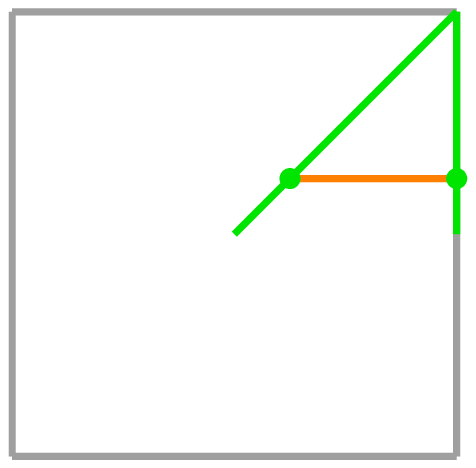}}\picskip{8}
Instead of interpolating over the whole upper right quadrant we can now
do this in a triangle only, and get sharp values, exhibiting the cusp at the diagonal.
Fix $r\in(2,\infty)$.
We have
\[
\|[X,Y]\|_r\leq 2^{1-1/r} \|X\|_r\|Y\|_r {\rm \ and\ } \|[X,Y]\|_\infty\leq 2 \|X\|_\infty\|Y\|_r
\]
and consider $K_Y$.
By
\[\frac{1}{p}=\frac{1-\theta}{r}+\frac{\theta}{\infty}\]
and consequently $\theta=1-r/p$
we get
$C_{p,p,r}\leq \left(2^{1-1/r}\right)^{1-\theta} 2^\theta = 2^{1-1/p}$, as claimed.

\medskip
\parpic[l]{\includegraphics[scale=\pscl]{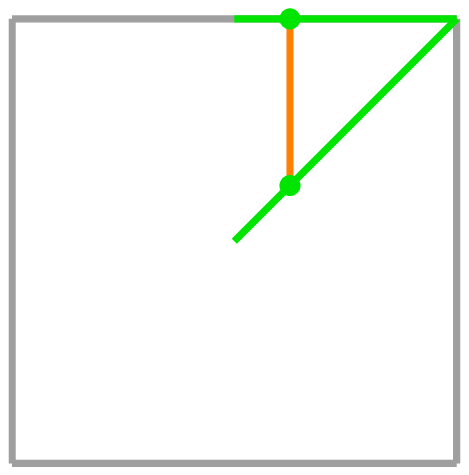}}\picskip{8}
Similarly, for the second triangle, fix $p\in(2,\infty)$. We know the bounds for
the two base points:
\[
\|[X,Y]\|_p\leq 2^{1-1/p} \|X\|_p\|Y\|_p {\rm \ and\ } \|[X,Y]\|_p\leq 2 \|X\|_p\|Y\|_\infty.
\]
Now, by interpreting the commutator as the linear map
$K_X$ we obtain
\[\frac{1}{r}=\frac{1-\theta}{p}+\frac{\theta}{\infty}\]
or $\theta=1-p/r$
and
$C_{p,p,r}\leq \left(2^{1-1/p}\right)^{1-\theta} 2^\theta = 2^{1-1/r}$.

\medskip
\parpic[l]{\includegraphics[scale=\pscl]{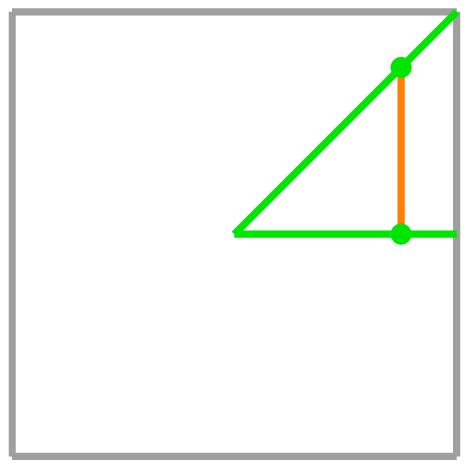} \includegraphics[scale=\pscl]{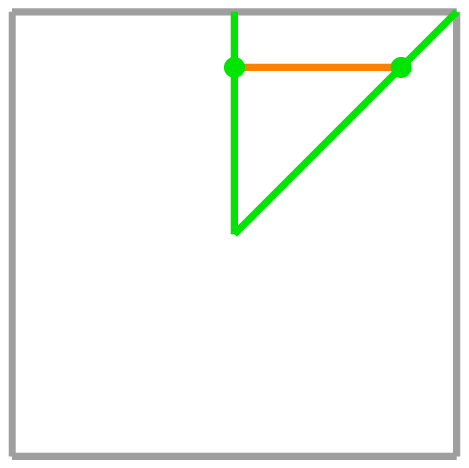}}\picskip{6}
We ought to remark that it is easier to interpolate the two triangles along the other respective direction
(i.e.\ fixing the other variable as we have done above), as the
constant in the corresponding inequalities is the same for both interpolation base points
and, hence, automatically yields exactly this value for all interpolant inequalities.
We don't even need to determine the relation linking the parameter $\theta$ with $p$ or $r$.
Note that this only works because some of the base points were determined before by other interpolation
steps, whereas the variant given first relies only on the trivial estimates on the boundary
of the square.

\medskip
\parpic[l]{\includegraphics[scale=\pscl]{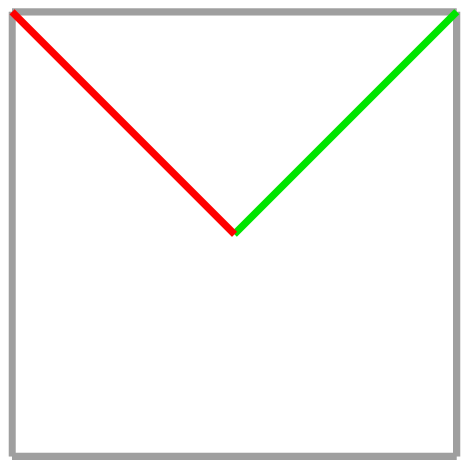}}\picskip{6}
In the last, upper left quadrant, another cusp line appears. In order to proceed
in a similar way as we did with the upper right quadrant, the anti-diagonal $r=p'$ is needed.
Fortunately we can obtain these values by a simple duality argument, as follows.
For $p\geq 2$ we have

\[\|[X,Y]\|_p\leq 2^{1-1/p}\|X\|_p\|Y\|_p\]

obtained from
tensor product interpolation (green line). Then
for any $Y$ with $\|Y\|_p=1$ one has

\[2^{1-1/p}=\sup_X \frac{\|K_Y(X)\|_p}{\|X\|_p}.\]

We conclude

\[2^{1/p'}=\sup_X \frac{\|K_Y^*(X)\|_{p'}}{\|X\|_{p'}} = \sup_X \frac{\|K_Y(X)\|_{p'}}{\|X\|_{p'}},\]

giving
$\|[X,Y]\|_{p'}\leq 2^{1/p'}\|X\|_{p'}\|Y\|_p$ for $p'\leq 2$,
which is the assertion for the anti-diagonal (red line);
see the proof of Proposition \ref{propCpqrsym} for details on the above equality.

\medskip
\parpic[l]{\includegraphics[scale=\pscl]{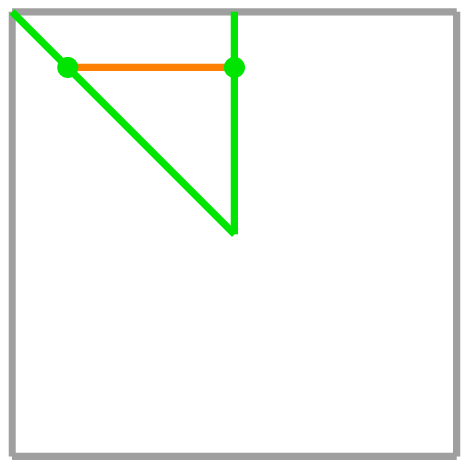} \includegraphics[scale=\pscl]{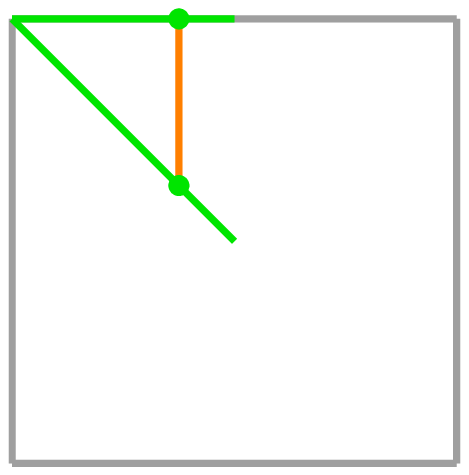}}
It should be clear now that in complete analogy to the upper right quadrant we have to
interpolate the two triangles separately. Again it doesn't matter along which direction
we fix one of the norm indices.
\vspace*{25pt}

\medskip
\parpic[l]{\includegraphics[scale=\pscl]{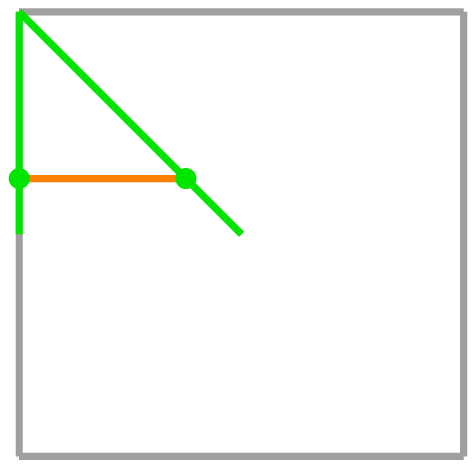} \includegraphics[scale=\pscl]{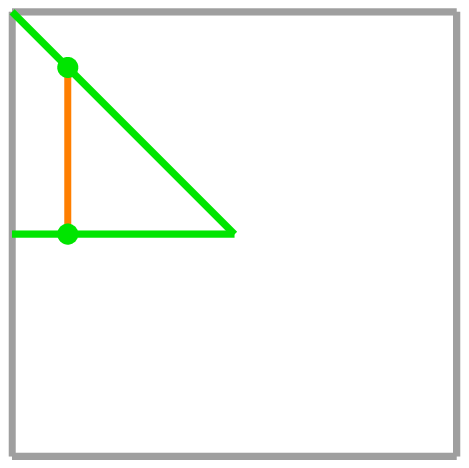}}\picskip{6}
Of course, one of the directions is easier than the other.
We leave it to the reader
to find out which one is preferable.\\
We skip the details as the result does also follow from a plain duality argument
(as we have done for the anti-diagonal). By this argument we directly get
$C_{p,p,r}=C_{p',p',r}$, which may be interpreted as a reflection symmetry of the constant about the line $p=2$.

\subsection{A little short-cut}

\medskip
\parpic[l]{\includegraphics[scale=\pscl]{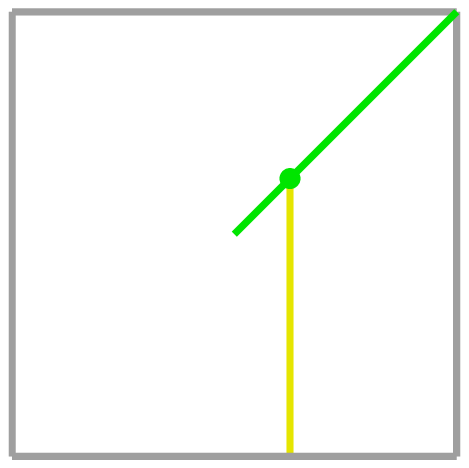} \includegraphics[scale=\pscl]{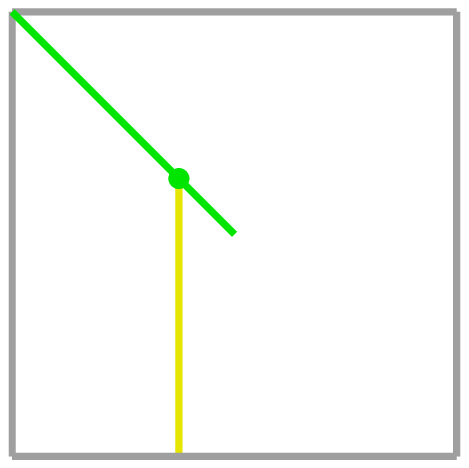}}\picskip{6}
The case $p\geq2, r\leq p$ and similarly $p\leq2, r\leq p'$ can be done in an even easier way using norm index monotonicity.
By this, the value $2^{1-1/p}$ (known for points on the green line in the left picture)
extends to all $1\leq r<p$.
So, after the diagonal interpolation (as the really first step), this attempt may replace the previous investigations
of the lower right quadrant and one triangle. Similarly, once we obtain the anti-diagonal by duality, the observations for the
lower left quadrant follow automatically.

\medskip
\parpic[l]{ \includegraphics[scale=\pscl]{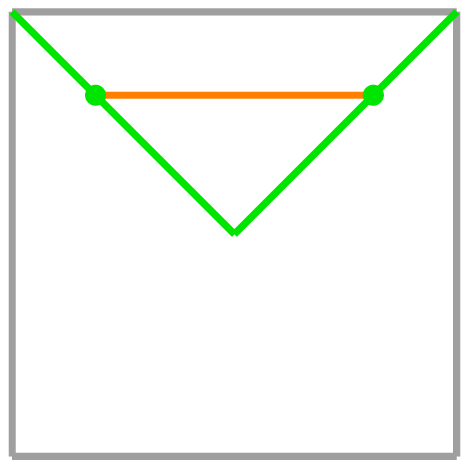}}
Also note that the remaining two single triangles can be merged into a single step. This cannot be done
with help of the norm index monotonicity, and interpolation (with fixed $r$) becomes necessary.
However, since both base inequalities admit the same constant, this turns out to be pretty easy.\vspace*{10pt}

\medskip
\parpic[l]{ \includegraphics[scale=\pscl]{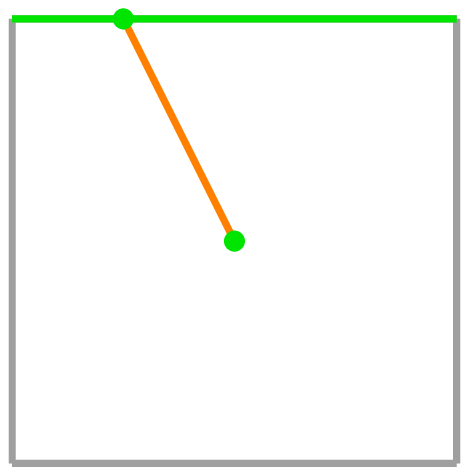}}\picskip{6}
Finally, we remark that the diagonal tensor interpolation, the dual anti-diagonal values
and the triangle interpolation between both can also be merged into a single step
by applying the multilinear extension of the Riesz-Thorin theorem.
However, we will not give more details about this since the treatment
of (\ref{eqRTconvcomb}) requires the synchronization of three equalities
and the calculation of the value of the interpolated bound is no longer that easy.

\section{Generalisation}\label{sec3}
In the previous section we have proven Theorem \ref{thmCppr}, which is really a special case of inequality
(\ref{eqNIpqr}). In the present section
we want to try our two main tools, as well as some slightly more delicate things,
to see how much extra mileage they allow us on the road towards a full proof of inequality (\ref{eqNIpqr}).
In that sense, this section is really a continuation of Section \ref{sec2}.
The main result of these investigations can be found at the end of this section.

\medskip
\parpic[l]{\includegraphics[scale=\psclt]{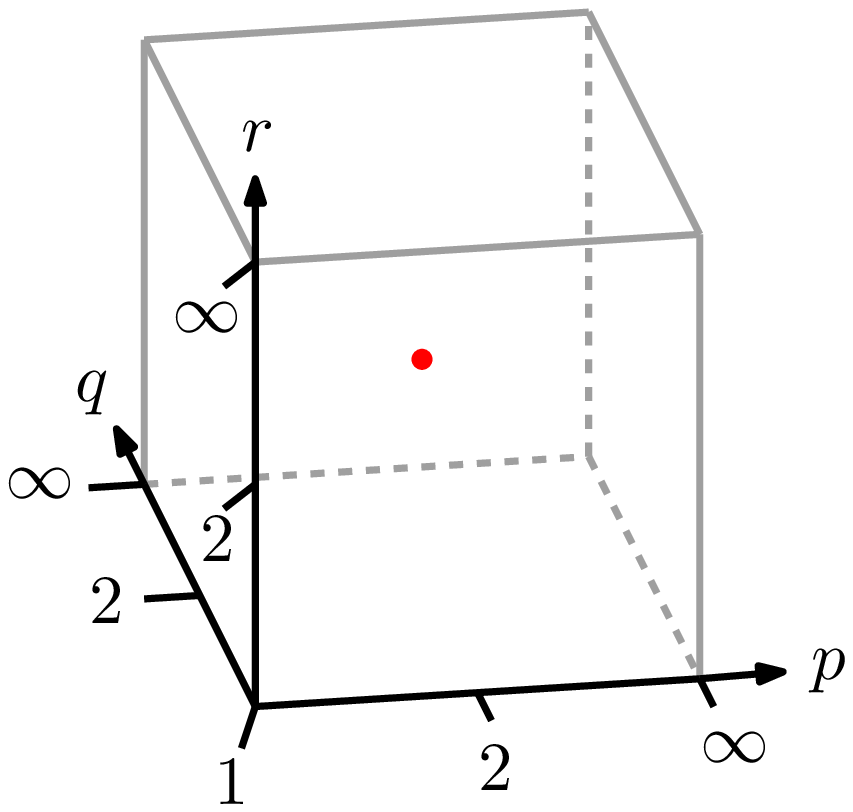}}\picskip{10}
For the general situation (\ref{eqNIpqr}) three norm indices $p,q$ and $r$
have to be depicted, requiring three-dimensional images.
In what follows, we transform the cube $[1,\infty]^3$ by the mapping $\ig\!^3$
and represent its image using a perspective projection from a fixed viewing direction.
Under these circumstances we can again drop axis labels, just as in Section \ref{sec2}.
Note that $(2,2,2)$ is again represented by the cube's center (red point).\\

\vspace*{30pt}
As of now, regions of parameter space will be colored differently depending on the rule
that determines $C_{p,q,r}$.

\medskip
\parpic[l]{\includegraphics[scale=\psclt]{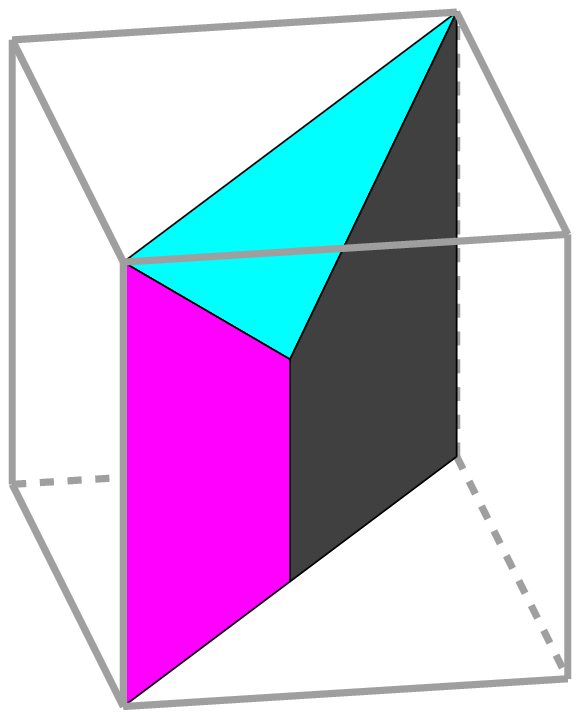} \includegraphics[scale=\psclt]{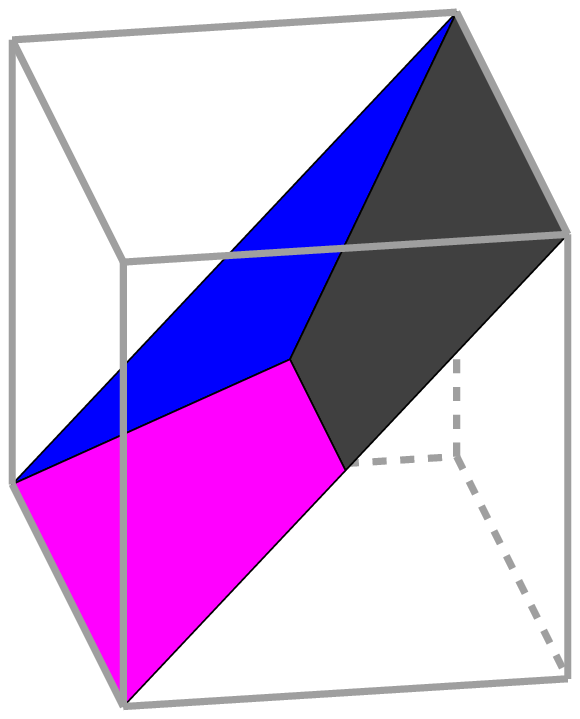}}\picskip{11}
First, we picture the (now proven) originally conjectured special case in the general context. We know the values of
\[C_{p,p,r}=\max\left\{{\color[rgb]{1,0,1}2^{1/p}},{\color[rgb]{0.25,0.25,0.25}2^{1-1/p}},{\color[rgb]{0,1,1}2^{1-1/r}}\right\}\]
and by swapping the roles of $X$ and $Y$ also of
\[C_{p,q,p}=\max\left\{{\color[rgb]{1,0,1}2^{1/p}},{\color[rgb]{0.25,0.25,0.25}2^{1-1/p}},{\color[rgb]{0,0,1}2^{1-1/q}}\right\}.\]

\parpic[l]{\includegraphics[scale=\psclt]{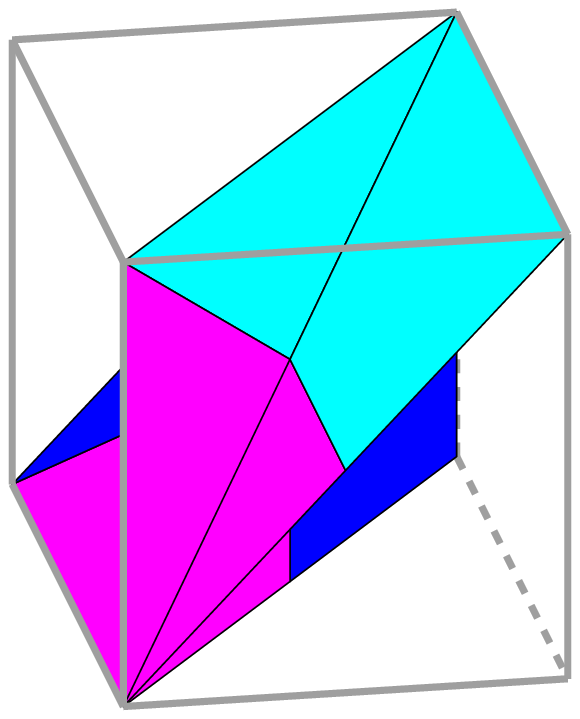}}\picskip{11}
These constants are represented by triplets on the planes $q=p$ and $r=p$.
Due to the properties of our scaling function $\ig\!^3$ the latter are indeed planes
(recall similar statements for $\ig\!^2$ given in Section \ref{sec1.2}). \\
We combine the two results and moreover modify them in a way
that turns out to be more suitable for what follows.
Naturally, one has ${\color[rgb]{0.25,0.25,0.25}2^{1-1/p}}={\color[rgb]{0,0,1}2^{1-1/q}}$ and
${\color[rgb]{0.25,0.25,0.25}2^{1-1/p}}={\color[rgb]{0,1,1}2^{1-1/r}}$
in the two planes, respectively.
\vspace*{30pt}

\medskip

\subsection{Monotonicity conquers (almost) all}\label{sec3.1}

\medskip
The validity of the conjecture naturally extends to some of the cases with $q\neq p$, by applying
the norm index monotonicity argument.

\medskip
\parpic[l]{\includegraphics[scale=\psclt]{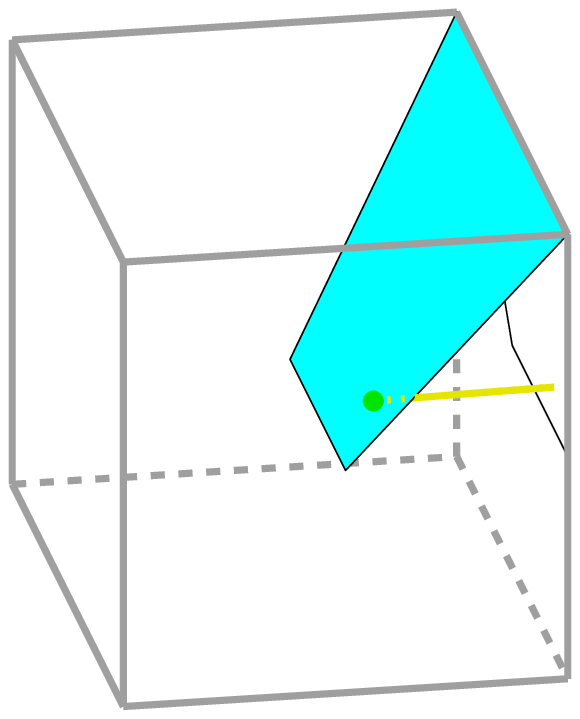} \includegraphics[scale=\psclt]{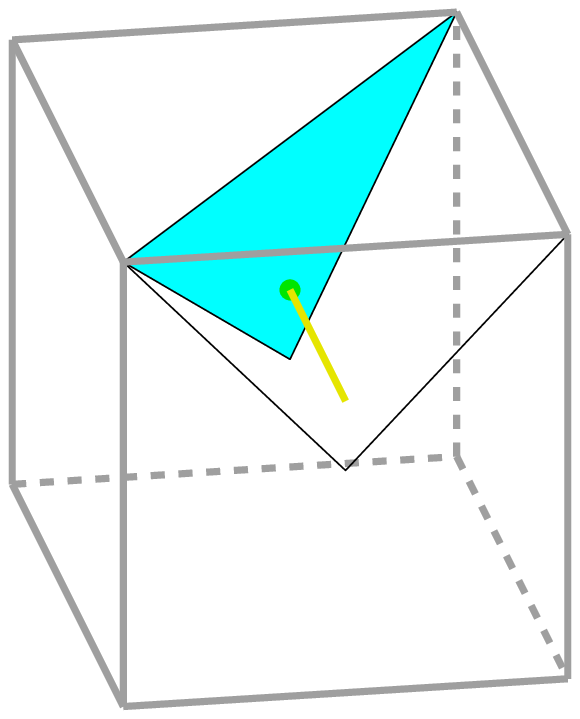}}
First take a look at triplets connected to the constant $2^{1-1/r}$. This value
does not depend on $p$ and $q$. So, we choose some point $(r,q,r)$ on the pictured segment
in the left image. We obtain, for all $p\geq r$,
\begin{align*}
\|[X,Y]\|_p & \leq \|[X,Y]\|_r \\
 & \leq 2^{1-1/r}\|X\|_q\|Y\|_r
\end{align*}
for arbitrary matrices $X$ and $Y$.
Moreover, for points $(p,p,r)$ as in the right image we get
\[\|[X,Y]\|_p\leq 2^{1-1/r}\|X\|_p\|Y\|_r\leq 2^{1-1/r}\|X\|_q\|Y\|_r\]
for any $q\leq p$.

\medskip
\begin{minipage}[b]{\pwdtt}
\includegraphics[scale=\psclt]{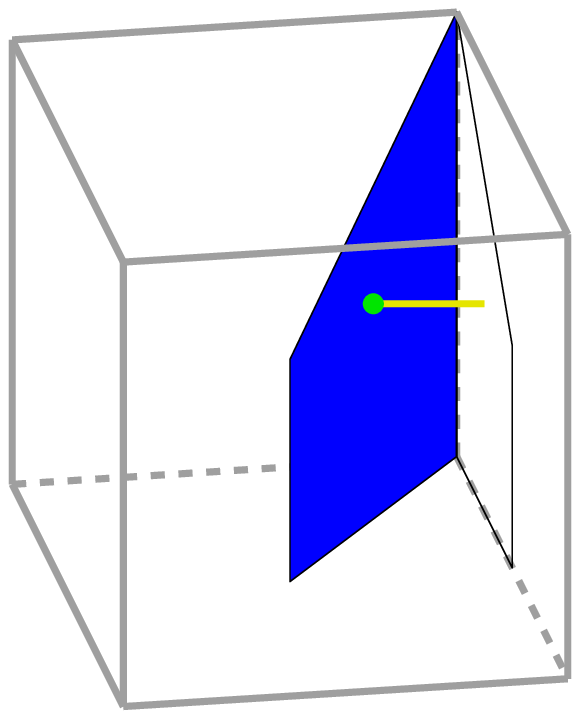} \includegraphics[scale=\psclt]{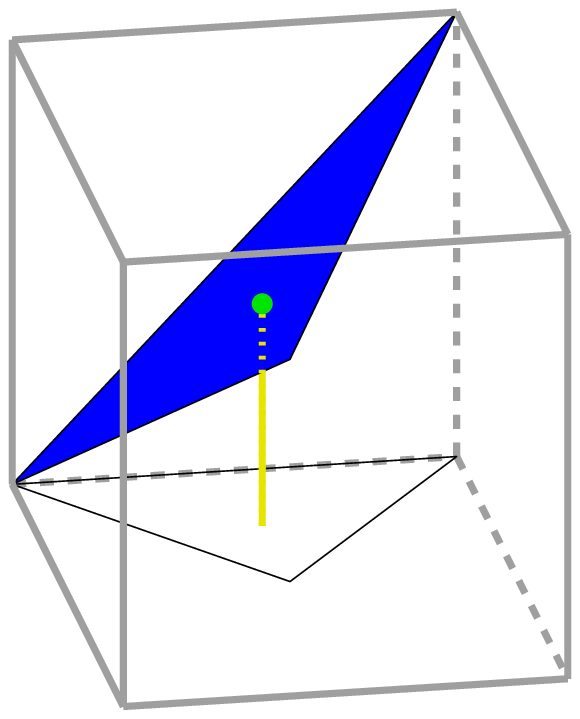}

\includegraphics[scale=\psclt]{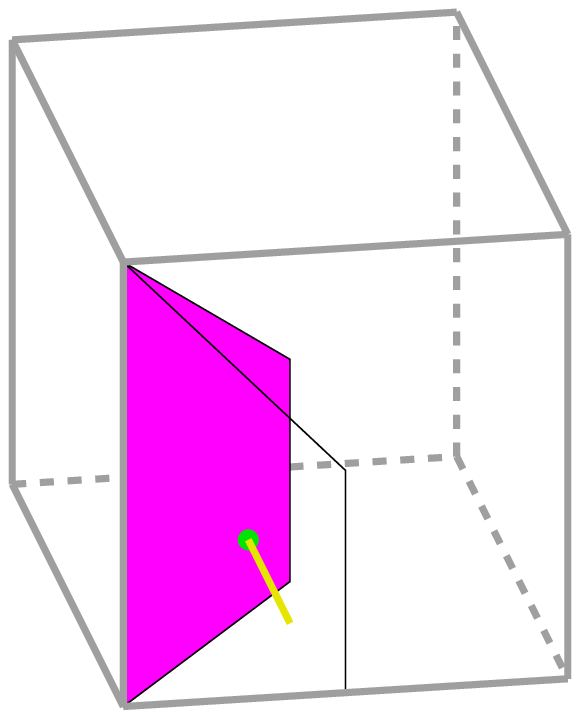} \includegraphics[scale=\psclt]{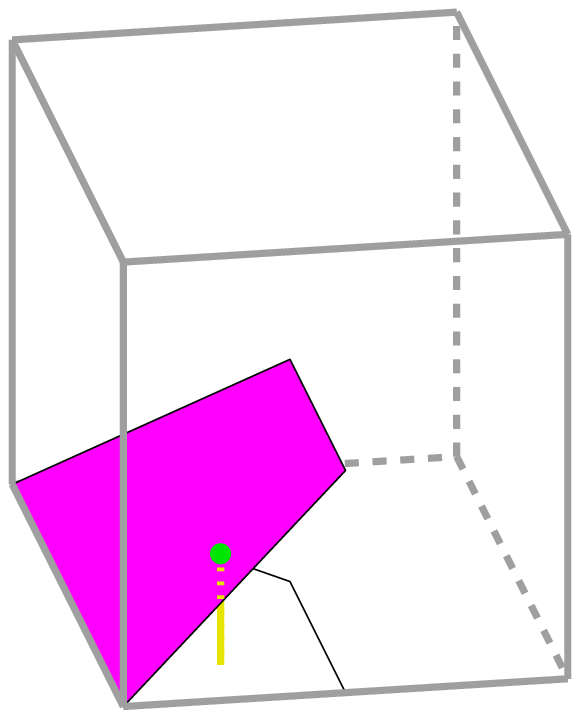}
\end{minipage}\begin{minipage}[b]{\twdtt}
For triplets belonging to $2^{1-1/q}$ we may argue in an analogous way for any $p\geq q$ (left)
and also for any $r\leq p$ (right).

As we only obtain upper bounds, we also need an example achieving equality.
One such is given by (\ref{eqspec2}).
Note that matrices of rank one are
essential for achieving equality in monotonicity relations. We will treat this
in more detail later.

\medskip

Similarly, for the segment where the constant is $2^{1/p}$, which is independent of $q$ and $r$, one can extend
the bound to $q\le p$ (left) and $r\leq p$ (right). Taking into account
(\ref{eqspec1}) we see that the value is sharp in these areas.
\end{minipage}

\medskip
\parpic[l]{\includegraphics[scale=\psclt]{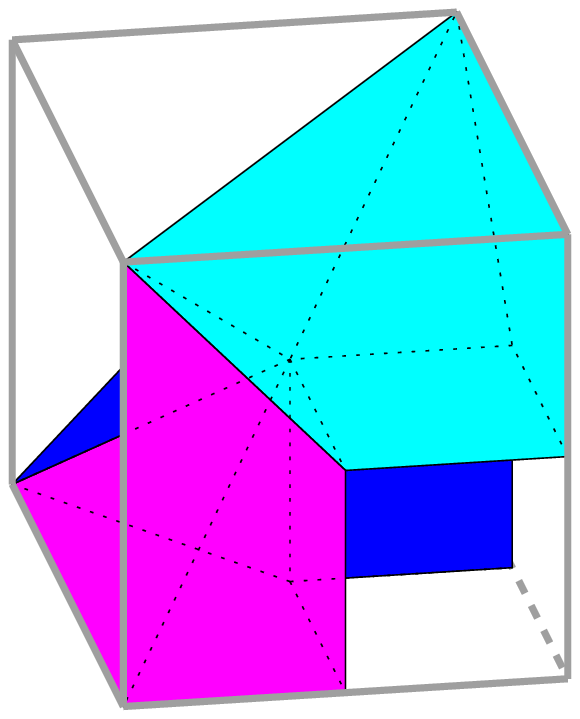}}

Summing up the results obtained so far, we get that the constant of Theorem \ref{thmCppr}
is valid also for a huge part of the general setting, namely for all $(p,q,r)$ with
$q\leq p$ or $r\leq p$, but not with all of $p>2,q<2$ and $r<2$.
Here we have one more indication that the areas (like the processes) are a lot easier
to visualize than to capture in formulas.

We point out the reflection symmetry of the areas and their values.
The light blue area is the image of the dark blue area under reflection about the plane $q=r$,
and the pink area is symmetric about that plane. One can even check that the value of $C_{p,q,r}$
equals the value in its mirror point.
This symmetry originates from the symmetry of $C$ under interchanging both $X$ with $Y$,
and $r$ with $q$,
as will be discussed in more detail at the end of this section (Proposition \ref{propCpqrsym}).

\subsection{Some more sophisticated techniques}

\medskip
\parpic[l]{\includegraphics[scale=\psclt]{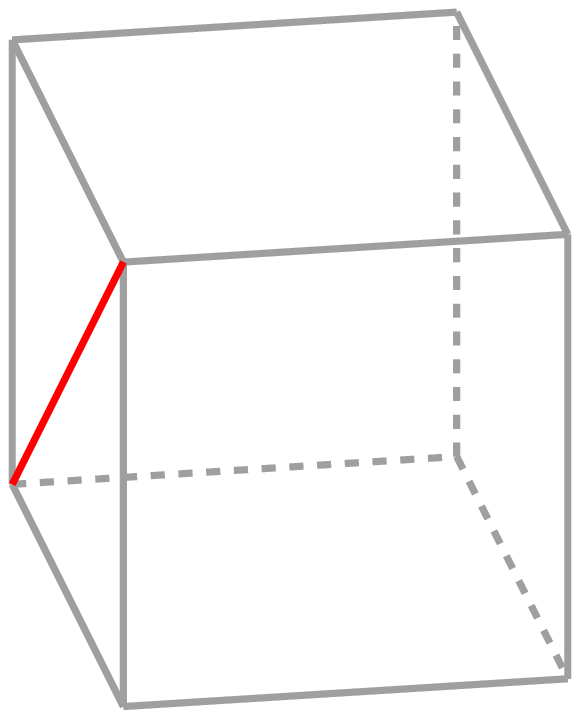}}\picskip{9}
For the next steps we need the values  for points $(1,q,q')$. These
can be obtained by a H\"{o}lder-type inequality which is true for
Schatten norms
\[\|XY\|_1\leq \|X\|_q \|Y\|_{q'},\]
whence, combined with the triangle inequality, one has
$\|XY-YX\|_1\leq 2\|X\|_q \|Y\|_{q'}$
giving $C_{1,q,q'}\leq 2$.
Example (\ref{eqspec1}) shows that equality can be achieved.
\vspace*{15pt}

\medskip
\parpic[l]{\includegraphics[scale=\psclt]{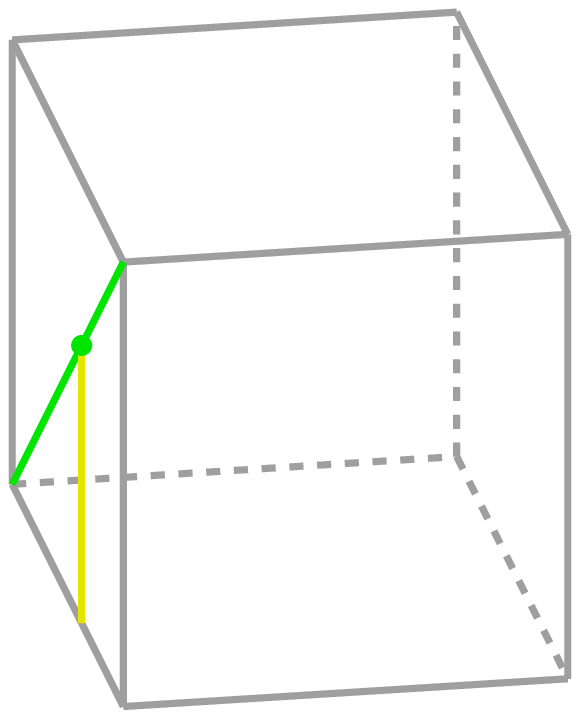}}

Now take any point $(1,q,q')$ from the line we just observed and apply
the monotonicity tool once more. We get
$C_{1,q,r}\leq 2$ for all $r\leq q'$.

Example (\ref{eqspec1}) again achieves equality here,
and the whole triangle admits the value 2.

The points in the triangle then serve as base points for the next interpolation step.\medskip\vspace*{40pt}

\bigskip
\parpic[l]{\includegraphics[scale=\psclt]{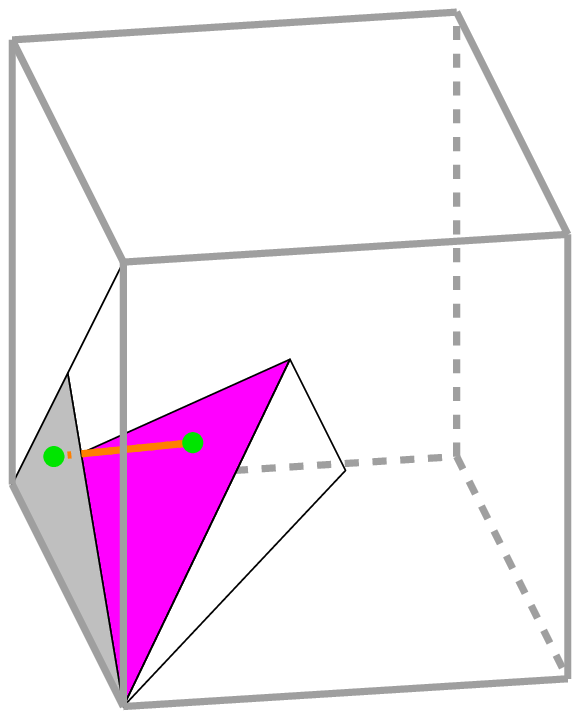} \includegraphics[scale=\psclt]{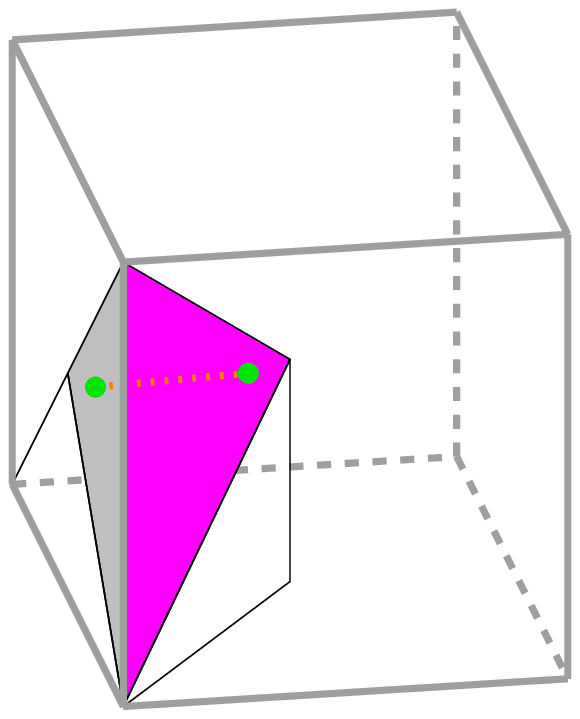}}\picskip{10}
Choose $q$ and $r$ arbitrarily in the grey triangle. We are going to interpolate only $p$
between 1 and $r$ (if $q\geq r$, left picture) or 1 and $q$ (if $r\geq q$, right picture).
For example, for the first case one has for (\ref{eqRTconvcomb})
\[\frac{1}{p}=\frac{1-\theta}{1}+\frac{\theta}{r}\]
which yields
\[C_{p,q,r}\leq 2^\frac{1/p-1/r}{1-1/r}\left(2^{1/r}\right)^\frac{1-1/p}{1-1/r} = 2^{1/p}.\]
Note that it doesn't matter whether we interpolate $K_X$ or $K_Y$, as both
$q$ and $r$ are fixed.

\parpic[l]{\includegraphics[scale=\psclt]{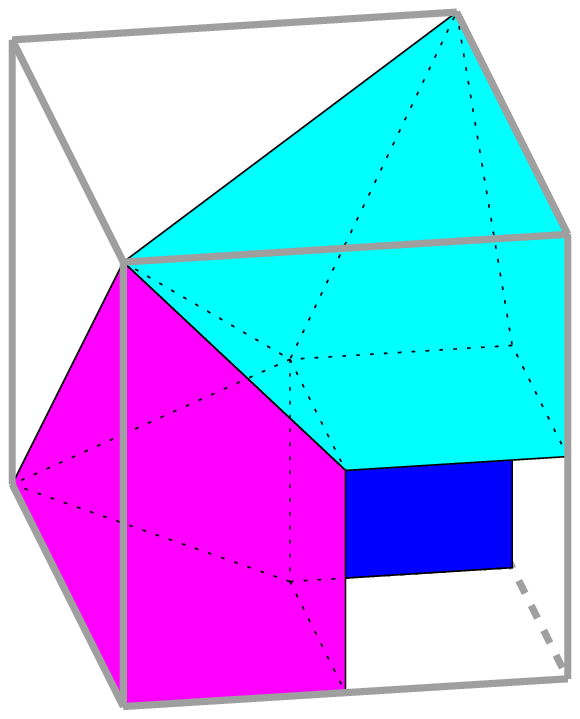}}\picskip{11}
For the second case we may proceed in an analogous way, or alternatively rely on the $q$-$r$-symmetry
already mentioned at the end of Section \ref{sec3.1}.\\
Also note that the area connected to ${\color[rgb]{1,0,1}2^{1/p}}$ is now of the same shape as the areas
of ${\color[rgb]{0,0,1}2^{1-1/q}}$ and ${\color[rgb]{0,1,1}2^{1-1/r}}$.\\
By this, we obtain two more symmetry planes, which are investigated in detail in Proposition \ref{propCpqrsym},
and which ensure the symmetry of the values and not only of the area's shape.

\vspace*{15pt}

\medskip
\parpic[l]{\includegraphics[scale=\psclt]{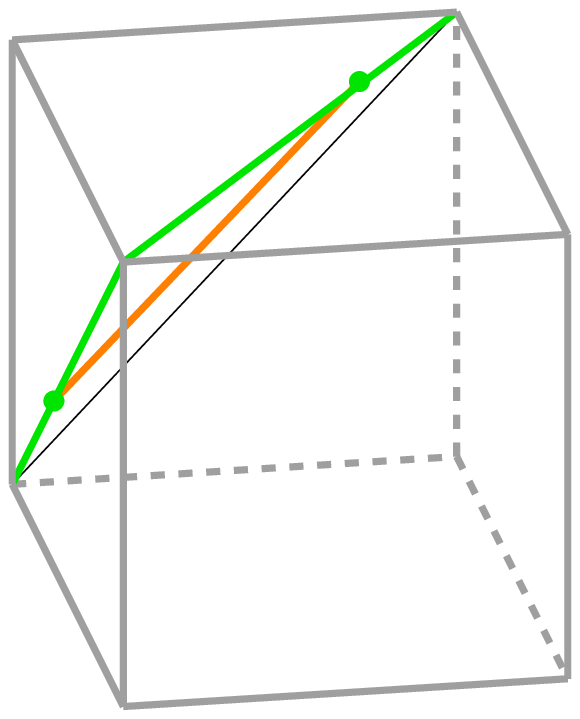}}\picskip{9}
Our next aim is to close the mould formed by the three areas for which the constant is known so far.
First we interpolate between the two points $(1,q,q')$ and $(q,q,\infty)$.
Both of them admit the constant 2, hence the points inbetween all share
this value.
The only remaining task is to determine which are the points inbetween.
Since $q$ is fixed, simple interpolation will work and requires
\[\frac{1}{p}=\frac{1-\theta}{1}+\frac{\theta}{q} \quad{\rm and}\quad \frac{1}{r}=\frac{1-\theta}{q'}+\frac{\theta}{\infty}.\]
Combining the latter we obtain the value 2 for all points satisfying
\[\frac{1}{p}=\frac{1}{q}+\frac{1}{r}.\]
We remark that $\ig\!^3$ maps the set of these triplets to a planar triangle.
After having done this calculation one gets an impression of the difficulties involved in the multilinear version,
when three equations come into play.

\medskip
\parpic[l]{\includegraphics[scale=\psclt]{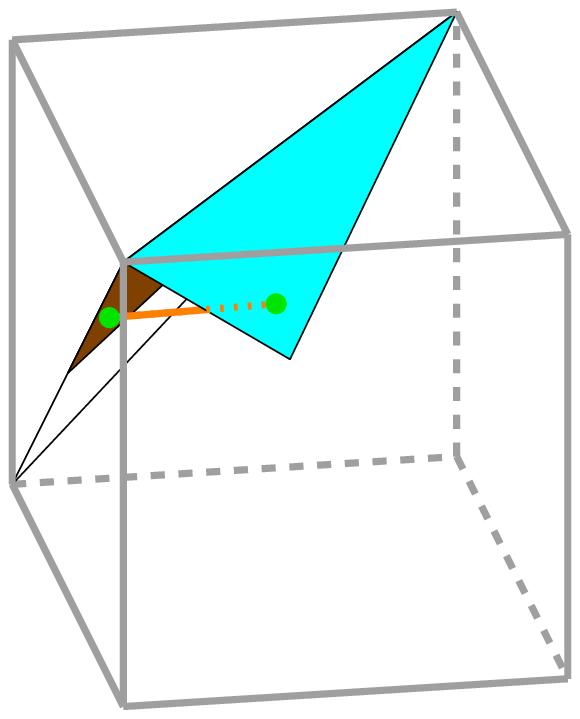} \includegraphics[scale=\psclt]{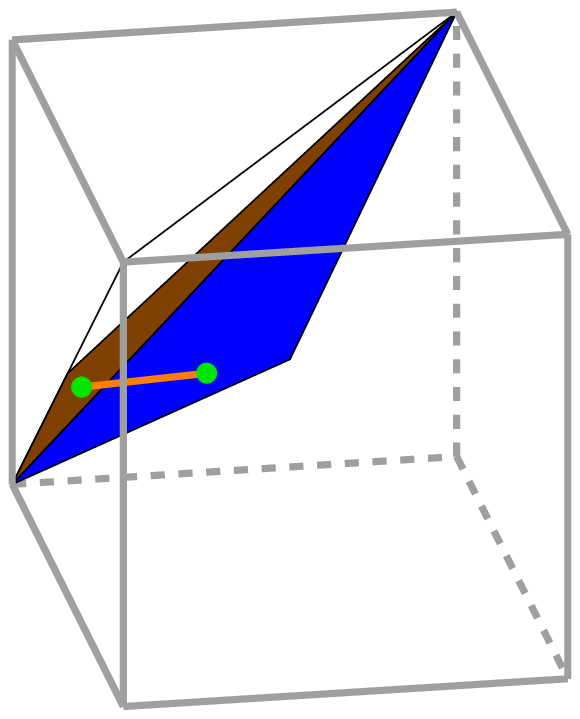}}\picskip{9}
Now we are in a position to close the gap between the plane that we have treated
and the known bodies, by means of interpolation. For this, again fix $q$ and $r$
as only $p$ will vary. Now choose $\tilde{p}$ such that $\frac{1}{\tilde{p}}=\frac{1}{q}+\frac{1}{r}$.
Hence, $(\tilde{p},q,r)$ lies in the brown triangle. The appropriate base point in the light blue
triangle is then given by $(q,q,r)$ (left image). For interpolants $(p,q,r)$ we
need to satisfy
\[\frac{1}{p}=\frac{1-\theta}{\tilde{p}}+\frac{\theta}{q}\]
which results in the bound
\[C_{p,q,r}\leq 2^{1+1/p-1/q-1/r}.\]
While this value seems to be rather exotic and maybe even perplexing
it is sharp nonetheless, as demonstrated by the example
\begin{equation}\label{eqspec3}
X=\left(\begin{array}{cc}
0 & 1 \\
1 & 0
\end{array}\right),\quad
Y=\left(\begin{array}{cc}
1 & 0 \\
0 & -1
\end{array}\right),\quad
XY-YX=\left(\begin{array}{cc}
0 & -2 \\
2 & 0
\end{array}\right).
\end{equation}

The second case is again done in a similar way or obtained by the $q$-$r$-symmetry.

\subsection{Trouble...}\label{sec3.3}

Knowledge of the constants for the plane $\frac{1}{p}=\frac{1}{q}+\frac{1}{r}$ successfully helped
to obtain the values $2^{1+1/p-1/q-1/r}$ in a triangular pyramid. So it is natural
to try the same for the pyramid opposite to it.
In order to perform the interpolation we need the
value $C_{1,\infty,\infty}$ of the pyramid's top. Unfortunately, this value is no longer
independent of the matrix size $d$. Thanks to the well-known inequalities
\begin{equation}\label{eqTriv1ii}
\|XY-YX\|_1\leq \|XY\|_1+\|YX\|_1\leq 2\|X\|_1\|Y\|_\infty\leq 2d\|X\|_\infty\|Y\|_\infty,
\end{equation}
we find a simple upper bound given by $2d$.

Using techniques similar
to those used in the last section, one gets the upper bound
\[C_{p,q,r}\leq 2d^{1/p-1/q-1/r}.\]
This follows in three steps: interpolating the line $p=1,r=\infty$, the plane $p=1$

\medskip
\,\hfill\includegraphics[scale=\psclt]{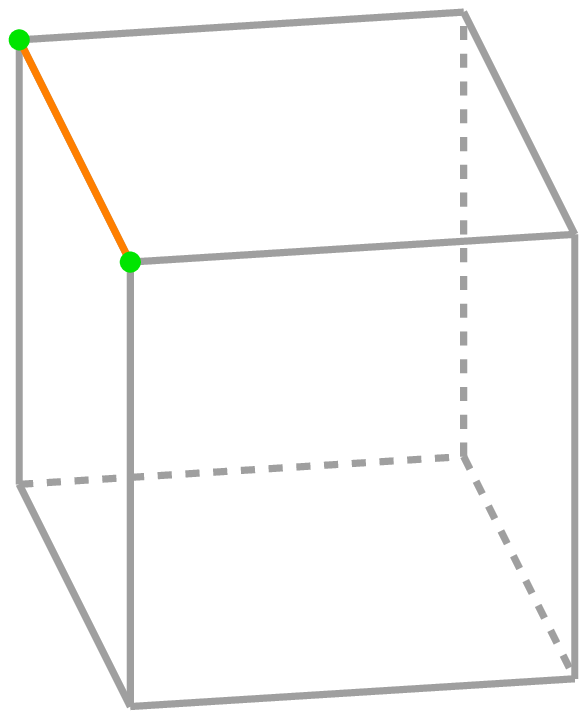} \hfill\begin{minipage}[b]{10pt}
$\Rightarrow$ \vspace*{40pt}
\end{minipage}\hfill
\includegraphics[scale=\psclt]{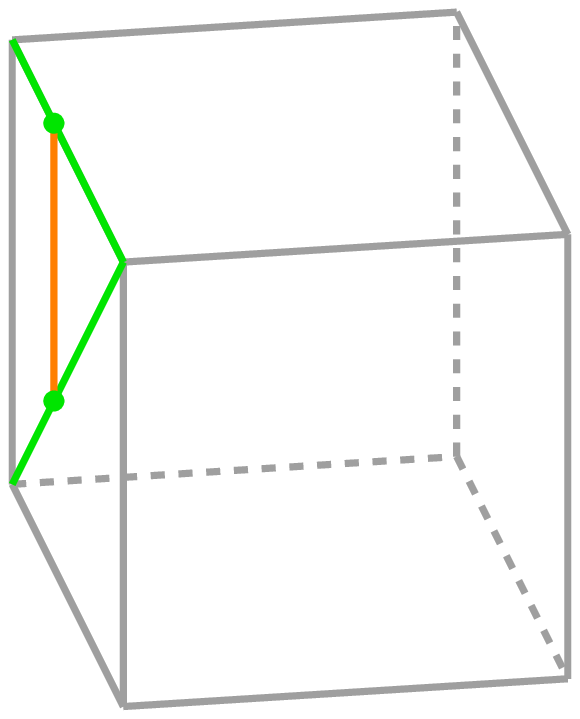} \hfill\begin{minipage}[b]{10pt}
$\Rightarrow$ \vspace*{40pt}
\end{minipage}\hfill
\includegraphics[scale=\psclt]{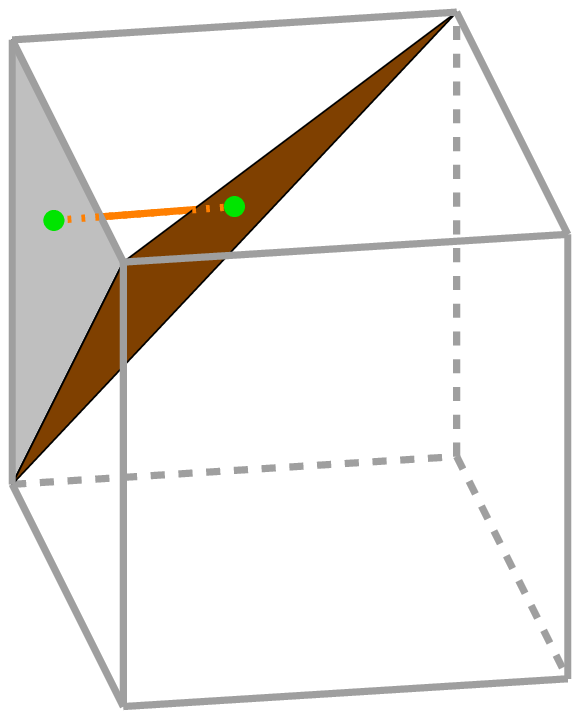}\hfill\,

\medskip
and finally the pyramid $\frac{1}{p}\geq\frac{1}{q}+\frac{1}{r}$.

At this point the interpolation method runs out of steam. Whereas for even $d$ the value
is shown to be sharp by the example
\[\tilde{X}=X\oplus X\oplus ...,\quad \tilde{Y}=Y\oplus Y\oplus...\]
with $2\times2$ matrices $X$ and $Y$ as in (\ref{eqspec3}),
we are unable to find an example when $d$ is odd. The reason may be
that in this case the estimate (\ref{eqTriv1ii}) is already not sharp.

\medskip
\parpic[l]{\includegraphics[scale=\psclt]{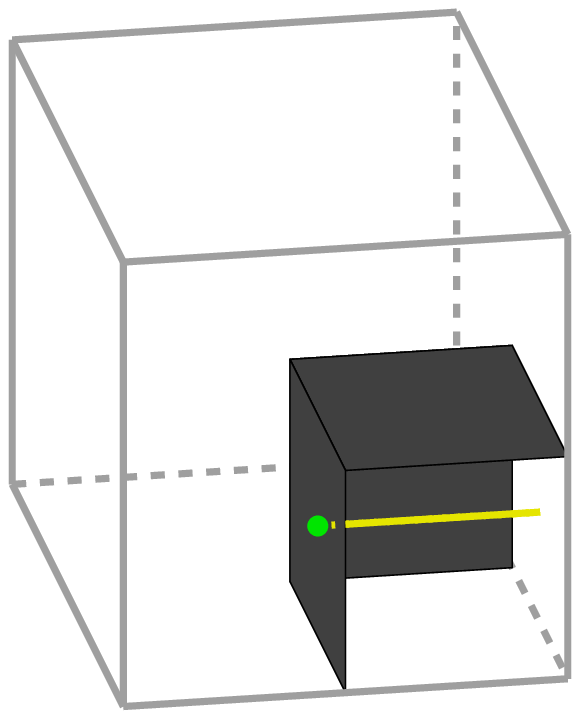}}\picskip{10}
The last area not yet investigated is the cube given by $p>2$ and $q,r<2$.
We do know the value of $C_{p,q,r}$ for three of its facets, namely $\sqrt{2}$.
The obvious method to apply is monotonicity. For instance, as indicated in the picture
we may write
\[\|XY-YX\|_p\leq \|XY-YX\|_2\leq \sqrt{2}\|X\|_q\|Y\|_r\]
for any $p>2$. By this, the upper bound $\sqrt{2}$ is extended to the
whole cube. Of course, one can use the monotonicity argument also with
reducing $q$ or $r$ based on the other facets instead.

\vspace*{15pt}

Sadly, this value is not sharp, and we can show this as follows. First we
observe that the value $\sqrt{2}$ is obtained solely by the knowledge of
$C_{2,2,2}=\sqrt{2}$, as the values on the facets themselves followed
from the value at the point $(2,2,2)$ using monotonicity. Now, we can use the fact that for $p_1>p_2$ equality
in $\|A\|_{p_1}\leq \|A\|_{p_2}$ holds if and only if ${\rm rank}\ A=1$. Hence, applying this for all indices,
we see that $X$, $Y$ and $XY-YX$ must all be matrices of rank one satisfying the equality
$\|XY-YX\|_2=\sqrt{2}\|X\|_2\|Y\|_2$. From Proposition 4.5 of \cite{BW2} we know
that without loss of generality two rank one matrices $X$ and $Y$ satisfy this
equality only if there are vectors $a,b$ such that $\|a\|_2=\|b\|_2=1, X=ab^*, Y=ba^*$
and $a^*b=0$. However, under those conditions, $XY-YX=aa^*-bb^*$ has rank two,
yielding a contradiction.

\subsection{The result, so far}

The previous steps obtained in this section (that is, the positive ones) add up to the following theorem.

\begin{theorem}\label{thmCpqr4}
For $(p,q,r)$ with $\frac{1}{p}\leq \frac{1}{q}+\frac{1}{r}$, excluding the octant
$p>2, q<2$ and $r<2$, one has
\[C_{p,q,r} = \max\{2^{1/p}, 2^{1-1/q}, 2^{1-1/r}, 2^{1+1/p-1/q-1/r}\}.\]

The four segments of $C_{p,q,r}$ corresponding to each of the four arguments of the maximum function
are given as follows:

\medskip
\,\hfill\includegraphics[scale=\psclt]{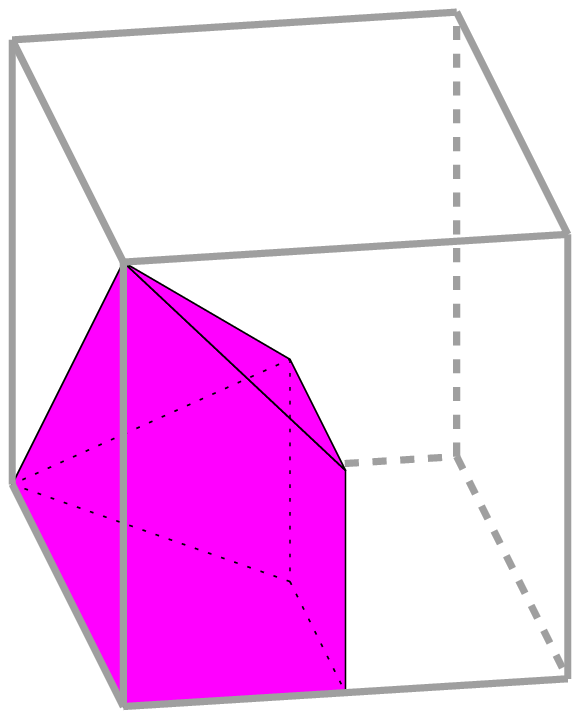}\quad \begin{minipage}[b]{70pt}
$2^{1/p}$ when\\

$q\leq p'$, $r\leq p'$,

$r\leq q'$ and

$p\leq 2$;
\end{minipage}\hfill \begin{minipage}[b]{70pt}
$2^{1-1/r}$ when\\

$r\geq p'$, $q\leq r$,

$q\leq p$ and

$r\geq 2$;
\vspace*{35pt}
\end{minipage} \quad\includegraphics[scale=\psclt]{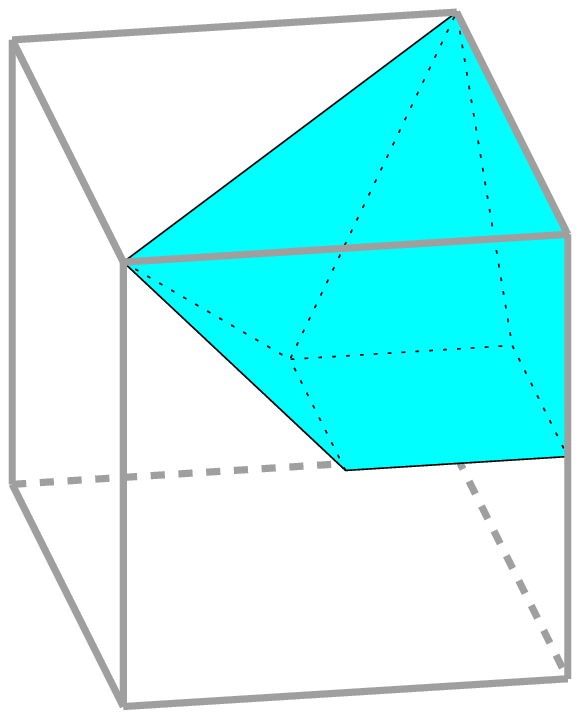}\hfill\,

\medskip
\,\hfill\includegraphics[scale=\psclt]{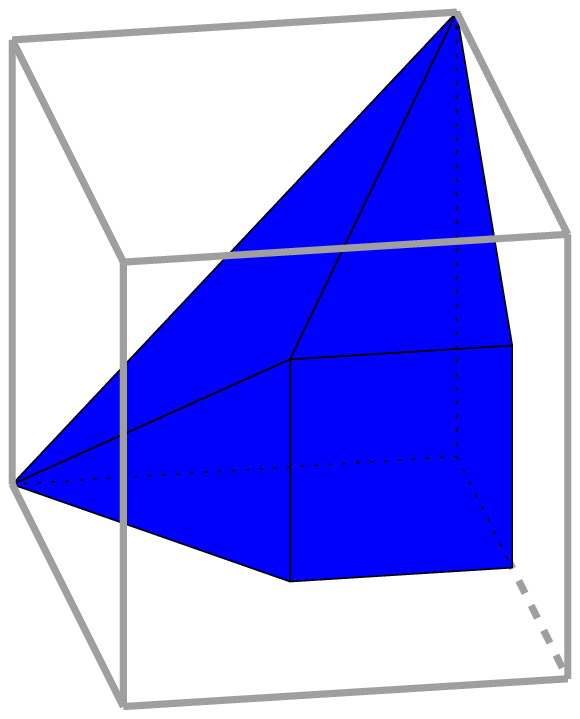}\quad \begin{minipage}[b]{70pt}
$2^{1-1/q}$ when\\

$q\geq p'$, $r\leq q$,

$r\leq p$ and

$q\geq 2$;
\end{minipage}\hfill \begin{minipage}[b]{70pt}
$2^{1+1/p-1/q-1/r}$

when\\

$\frac{1}{p}\leq\frac{1}{q}+\frac{1}{r}$,

$q\geq p$, $r\geq p$

and $r\geq p'$.
\vspace*{20pt}
\end{minipage} \quad\includegraphics[scale=\psclt]{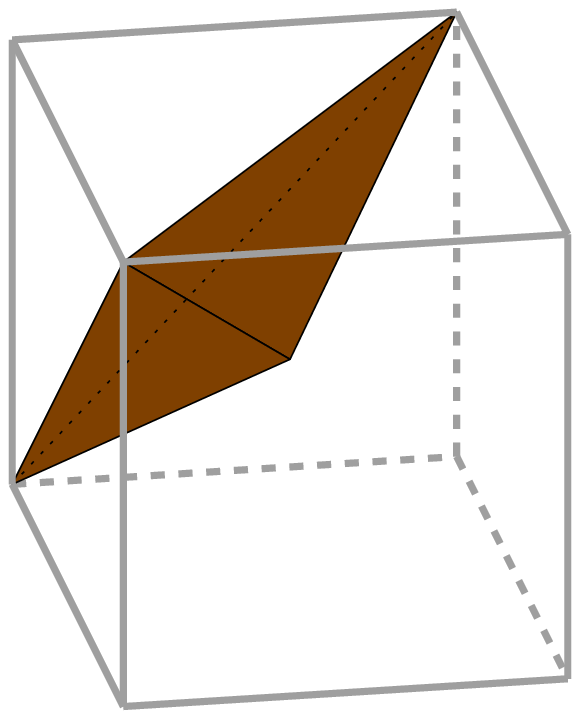}\hfill\,

\medskip
\begin{minipage}[b]{0.6\textwidth}
For $d\times d$ matrices of even size and $(p,q,r)$ with $\frac{1}{p}\geq \frac{1}{q}+\frac{1}{r}$
one has
\[C_{p,q,r} = 2d^{1/p-1/q-1/r}.\]

If $d$ is odd the latter is only an upper bound.
\end{minipage} \hfill\includegraphics[scale=\psclt]{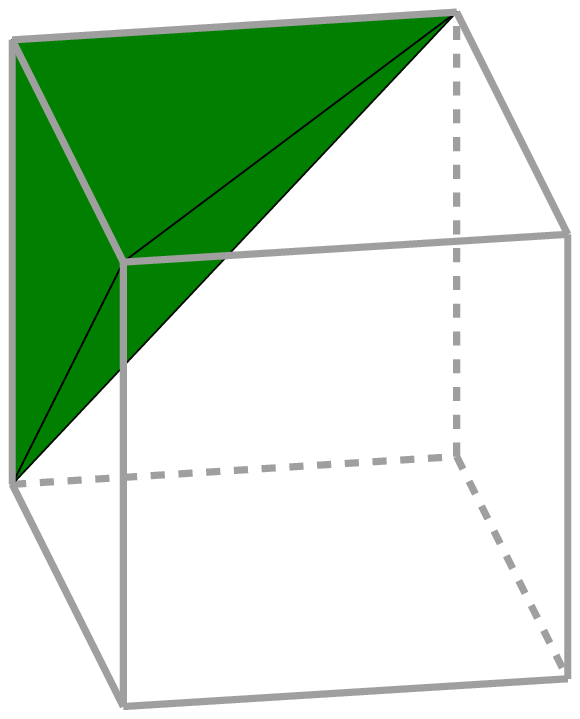}
\end{theorem}

Note that the constant for parameters in the region $\frac{1}{p}\leq \frac{1}{q}+\frac{1}{r}$
(i.e. the first four cases of Theorem \ref{thmCpqr4}) are independant of dimension.
Hence, the statement is also true in the infinite-dimensional setting of Schatten norms.

\medskip
The following result summarises all the symmetries we have encountered and
also encapsulates the duality arguments mentioned at the end of Section \ref{sec2.3}.

\begin{proposition}\label{propCpqrsym}
For any $(p,q,r)\in[1,\infty]^3$ one has
\[C_{p,q,r}=C_{p,r,q},\quad C_{p,q,r}=C_{r',q,p'},\quad C_{p,q,r}=C_{q',p',r}.\]
\end{proposition}

These three equalities represent the reflection symmetries of $C_{p,q,r}$ about the
planes $q=r$, $r=p'$, and $q=p'$, respectively:

\,\hfill\includegraphics[scale=\psclt]{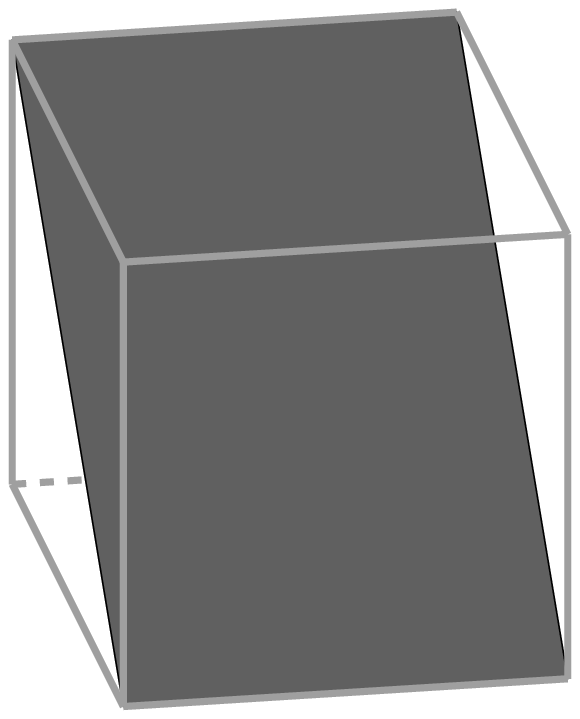}\hfill\includegraphics[scale=\psclt]{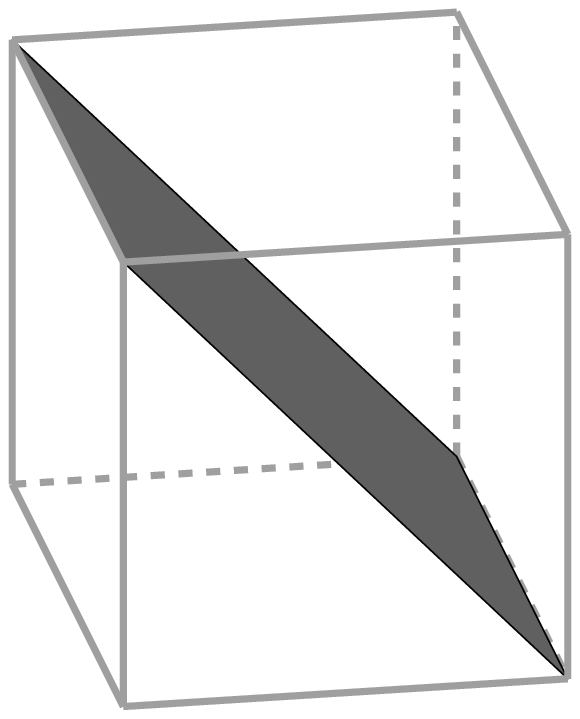}\hfill
\includegraphics[scale=\psclt]{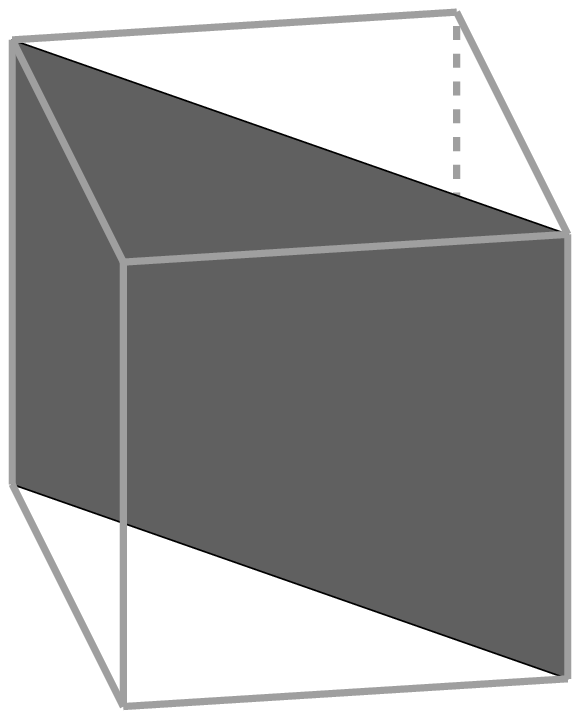}\hfill\,

The third picture generalises the duality statement from Section \ref{sec2.3}.

\medskip

\noindent\textit{Proof.}
The first equality is a mere consequence of $\|[X,Y]\|_p=\|[Y,X]\|_p$ and the
resulting possibility of changing the roles of $X$ and $Y$.

Now for the second equality, observe that for any fixed $X$ with $\|X\|_q=1$

\begin{eqnarray*}
&& \sup_Y \frac{\|K_X(Y)\|_p}{\|Y\|_r} = \sup_{\|Y\|_r=1} \sup_{\|W\|_{p'}=1} |\langle K_X(Y),W\rangle|\\
&& = \sup_{\|W\|_{p'}=1} \sup_{\|Y\|_r=1} |\langle Y,K_X^*(W)\rangle| = \sup_W \frac{\|K_X^*(W)\|_{r'}}{\|W\|_{p'}}
\end{eqnarray*}

and $K_X^* = -K_X$ imply the assertion.
Here, $\langle A,B\rangle=\tr B^*A$ denotes the inner product associated with
the Schatten classes.
The third equality is analogous or can be proved by combining
the first two equalities.
\qed

\medskip
The representations of norms as given in the last proof are called
variational characterisations and they will be of extraordinary use
in the following section, too.

\section{Extremal points\label{sec4}}
In the previous section we have squeezed the last drop out of the interpolation, monotonicity and duality methods, but two
areas in parameter space, a tetrahedron and a cube, still resist treatment. In the present section
we finally tackle these recalcitrant areas by finding the value of the constant in two specific points.
To do so, some new ideas are needed.

\subsection{The skeleton\label{sec4.1}}

\begin{figure}[h]
\centering\includegraphics[scale=0.75]{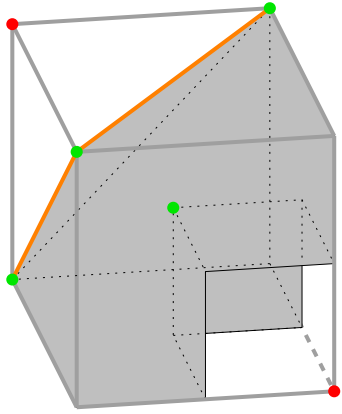}

\caption{Visualization of $(p,q,r)$ for which $C_{p,q,r}$ could be determined
by means of interpolation and monotonicity from the values of a couple of points (green).
The next logical targets are represented by red points.}\label{figKnown}
\end{figure}

In Figure \ref{figKnown} we depict all constellations $(p,q,r)$ covered so far in Sections
\ref{sec2} and \ref{sec3} by marking them in grey. All the values of the constant in these triplets
were the result of the knowledge of its value in only four points (or three, using symmetry), namely
$(2,2,2), (\infty,\infty,\infty)$ and $(1,1,\infty)$ or $(1,\infty,1)$ (marked green).
We also relied on the values of the points on the orange lines.
However, a closer look reveals that these may be obtained from
interpolation between two of the four green base points, too.

In Section \ref{sec3.3} we had a quick glance at the remaining two areas (white). In one situation,
monotonicity failed, while in the other the value at the interpolation base point
was likely not well-estimated (at least for odd-sized matrices).

In any case, the natural approach for carrying this further is to find the exact value of the constants $C_{1,\infty,\infty}$
($d$ odd) and $C_{\infty,1,1}$. These are the triplets marked red in the figure.

\subsection{The value of the constant at the corners\label{sec4.2}}
In this subsection we provide the value of the constant in the two corners just mentioned.
We prove the following theorem:
\begin{theorem}\label{thm1ii_i11}
For $d\times d$ matrices one has
\begin{enumerate}
\item[a)] $C_{1,\infty,\infty}=\left\{\begin{array}{ll}
    d |1+e^{i\pi/d}| = d\sqrt{2+2\cos(\pi/d)} & {\rm if\ } d {\rm\ is\ odd,}\\
    2d & {\rm if\ } d {\rm\ is\ even;}
\end{array}\right.$

\item[b)] $C_{\infty,1,1}=\sqrt{27}/4$.
\end{enumerate}
\end{theorem}

\noindent\textit{Proof of a).}

We only need to prove the formula for odd $d$, as the value
for even $d$ was already shown in Section \ref{sec3.3} in a much easier
way. However, as it requires no extra efforts, we nonetheless prove that particular result again
in the same fashion as for the odd case.

A variational characterisation for $C_{1,\infty,\infty}$ is given by
$$
C_{1,\infty,\infty} = \max_{X,Y} \{||XY-YX||_1: ||X||_\infty\le 1, ||Y||_\infty\le 1\}.
$$
Let us first fix $Y$.
The function to be maximised is convex in $X$, and the feasible set of $X$ is convex as well,
with extremal points given by the set of unitary matrices.
Thus, we can write:
$$
C_{1,\infty,\infty} = \max_{X,Y} \{||XY-YX||_1: X\mbox{ unitary }, ||Y||_\infty\le 1\}.
$$
A similar argument allows to conclude that $Y$ can also be restricted to the set of unitary matrices:
$$
C_{1,\infty,\infty} = \max_{X,Y \mbox{\small{ unitary}}} ||XY-YX||_1.
$$
In addition, the trace norm has a variational characterisation as well:
$$
||A||_1 = \max_{Z\mbox{\small{ unitary}}} |\trace ZA|.
$$
Thus we get a maximisation over three unitary matrices:
$$
C_{1,\infty,\infty} = \max_{X,Y,Z \mbox{\small{ unitary}}} |\trace Z(XY-YX)|.
$$

\bigskip

Every unitary matrix is unitarily equivalent to a diagonal matrix with all diagonal elements of modulus 1.
Applying this to $Y$, we get
\[Y=U\diag\left(e^{i\theta_1},e^{i\theta_2},\ldots,e^{i\theta_d}\right) U^*.\]
The matrix $U$ can be absorbed into $X$ and $Z$, so that w.l.o.g.\ we can restrict
$Y$ to be of this diagonal form. Indeed,
\begin{eqnarray*}
 && \trace (ZXY-ZYX) = \trace(ZXULU^* - ZULU^*X) \\
 && = \trace(ZUU^*XULU^* - ZULU^*XUU^*)
= \trace(Z'X'L-Z'LX'),
\end{eqnarray*}
where $Z'=U^*ZU$ and $X'=U^*XU$.

\bigskip

Then $[X,Y]$ can be rewritten as a Hadamard product: $XY-YX = A\circ X$, with $A$ a matrix with entries
$A_{jk}=e^{i\theta_k}-e^{i\theta_j}$.
The function to be maximised becomes
\beas
|\trace Z(XY-YX)| &=& |\trace Z(A\circ X)| = |\sum_{jk} Z_{kj}A_{jk}X_{jk}| \\
&\le& \sum_{jk} |Z_{kj}| \,\, |A_{jk}| \,\, |X_{jk}|.
\eeas
The Cauchy-Schwartz inequality leads to a further upper bound:
\beas
|\trace Z(XY-YX)|
&\le& \sum_{jk} |Z^T_{jk}| \,\, |A_{jk}| \,\, |X_{jk}| \\
&\le& \left(\sum_{jk} |Z^T_{jk}|^2 \,\, |A_{jk}|\right)^{1/2}
\,\,\left(\sum_{jk} |A_{jk}| \,\, |X_{jk}|^2\right)^{1/2}.
\eeas
Applying the maximisation over all unitary $X$ and $Z$ to both sides then yields
$$
\max_{X,Z}|\trace Z(XY-YX)| \le
\max_X \sum_{jk} |A_{jk}| \,\, |X_{jk}|^2,
$$
because both factors of the right-hand side could be maximised separately, and both maxima are equal.
Now note that the matrix with elements $|X_{jk}|^2$ is a doubly stochastic matrix (because $X$ is unitary).
Furthermore, the function to be maximised is linear in $|X_{jk}|^2$.
Hence, the maximum is achieved in extremal points of the set of doubly stochastic matrices.
By Birkhoff's theorem \cite{bhatia}, these are permutation matrices.
Thus we have a further reduction:
$$
\max_{X,Z}|\trace Z(XY-YX)| \le
\max_\pi \sum_{j} |A_{j\pi(j)}|,
$$
where the maximum is over all permutations $\pi$.
Observe that this inequality is actually an equality, as the left-hand side attains the right-hand side
for $Z^T$ and $X$ both equal to the permutation matrix representing $\pi$.

\bigskip

We are now left with calculating the maximum over all angles $\theta_j$ and all permutations $\pi$ of
$\sum_{j} |A_{j\pi(j)}| = \sum_j |e^{i\theta_{\pi(j)}}-e^{i\theta_j}|$.
This problem has a nice geometric interpretation.
The complex numbers $e^{i\theta_j}$ are points on the unit circle. The permutation $\pi$ maps every point
to another point, in a one-to-one fashion. If we draw edges from $e^{i\theta_j}$ to $e^{i\theta_{\pi(j)}}$
we obtain one or more polygons (in general, non-convex and self-intersecting), corresponding to the cycles of the permutation.
The problem is to distribute the points on the circle and choose the polygons so that the total length
of the edges (the total circumference of the polygon(s)) is maximised.

\bigskip

\begin{figure}[p]
\begin{center}
\hspace*{-25pt}\includegraphics[scale=0.5]{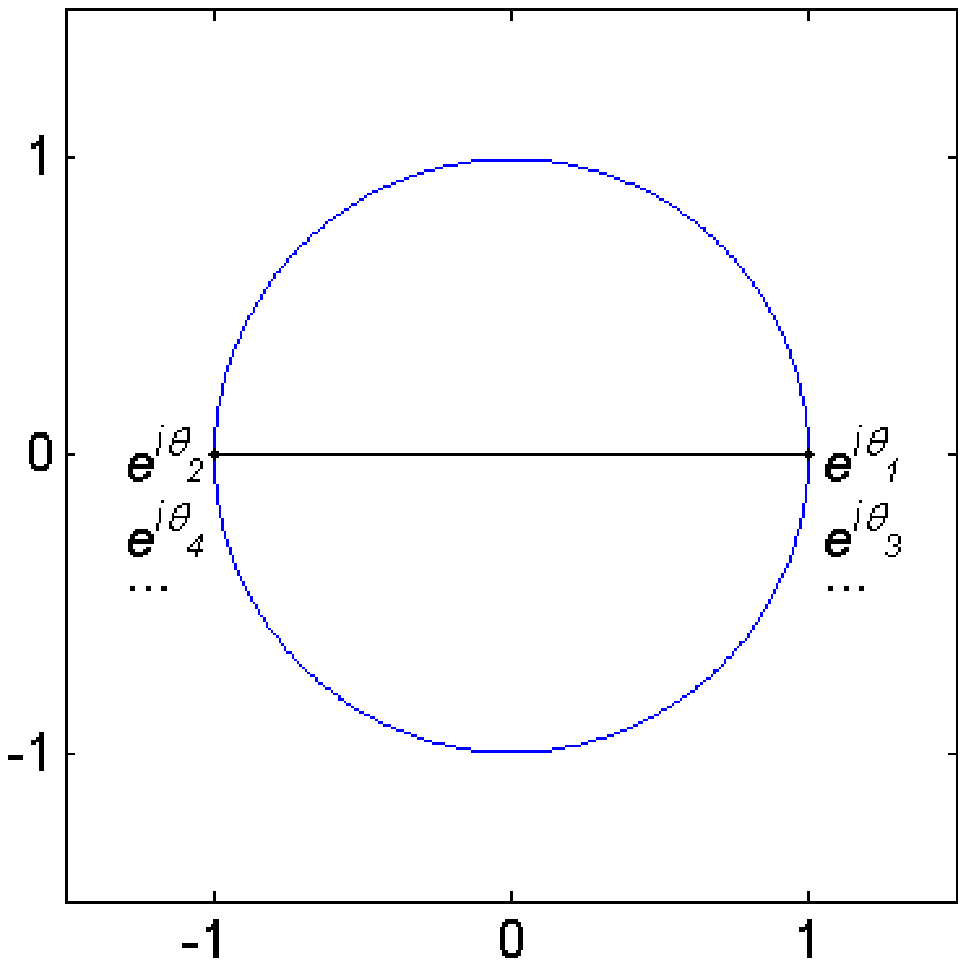}\hspace*{-25pt}
\includegraphics[scale=0.5]{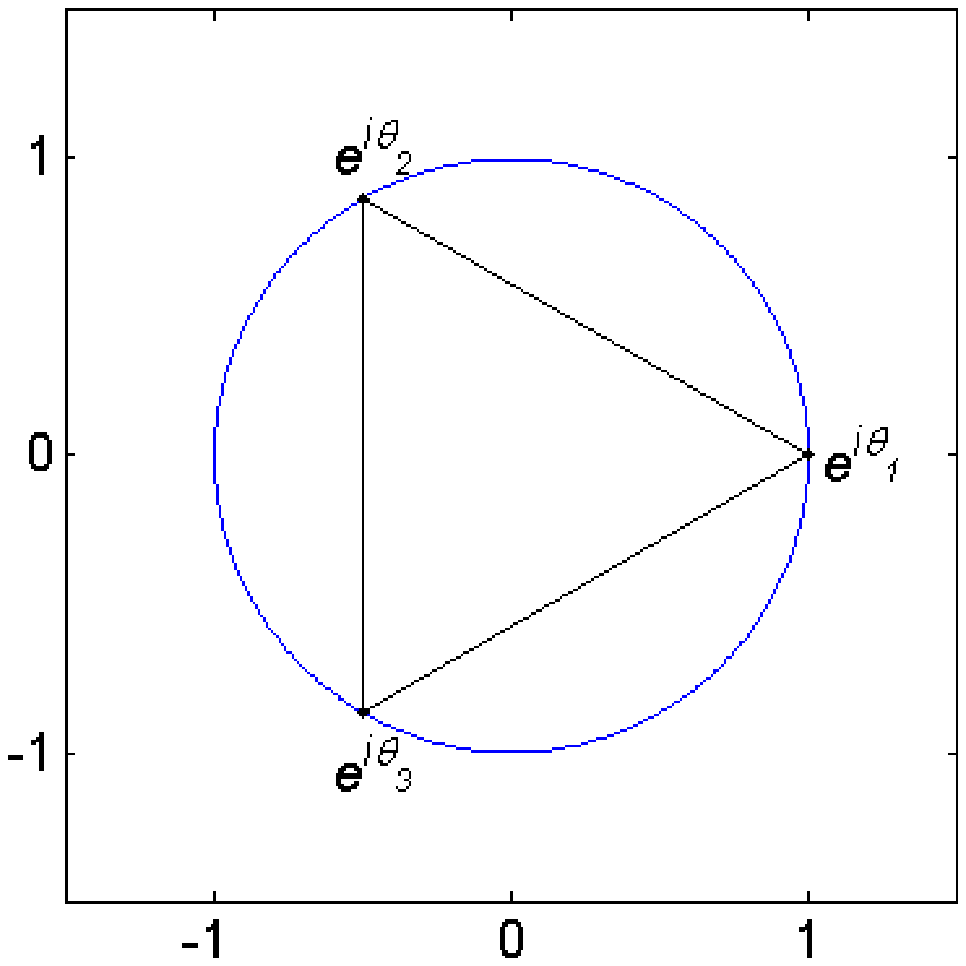}\hspace*{-25pt}
\caption{Distribution of $d$ points on the unit circle forming polygon(s) with
maximal total circumference: trivial solution for even $d$ (left) and obvious
configuration for $d=3$ (right).}\label{figStars23}
\end{center}
\end{figure}

\begin{figure}[p]
\begin{center}
\hspace*{-25pt}\includegraphics[scale=0.5]{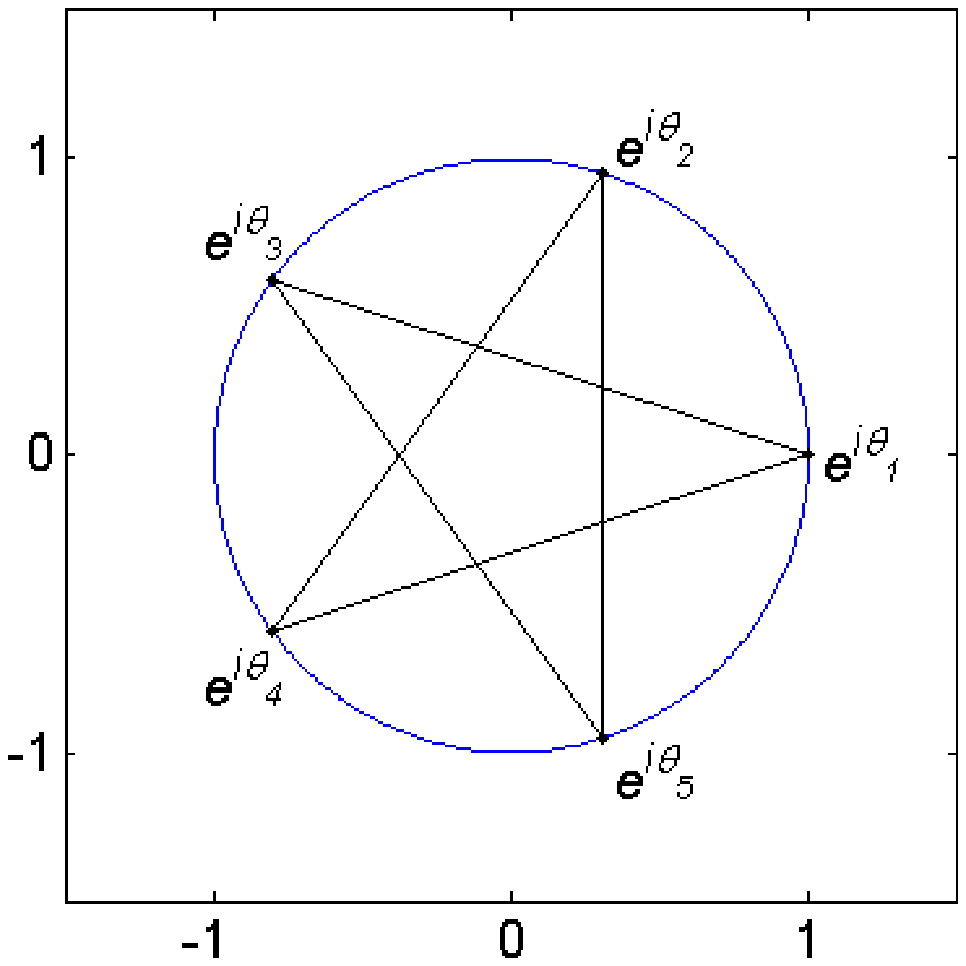}\hspace*{-25pt}
\includegraphics[scale=0.5]{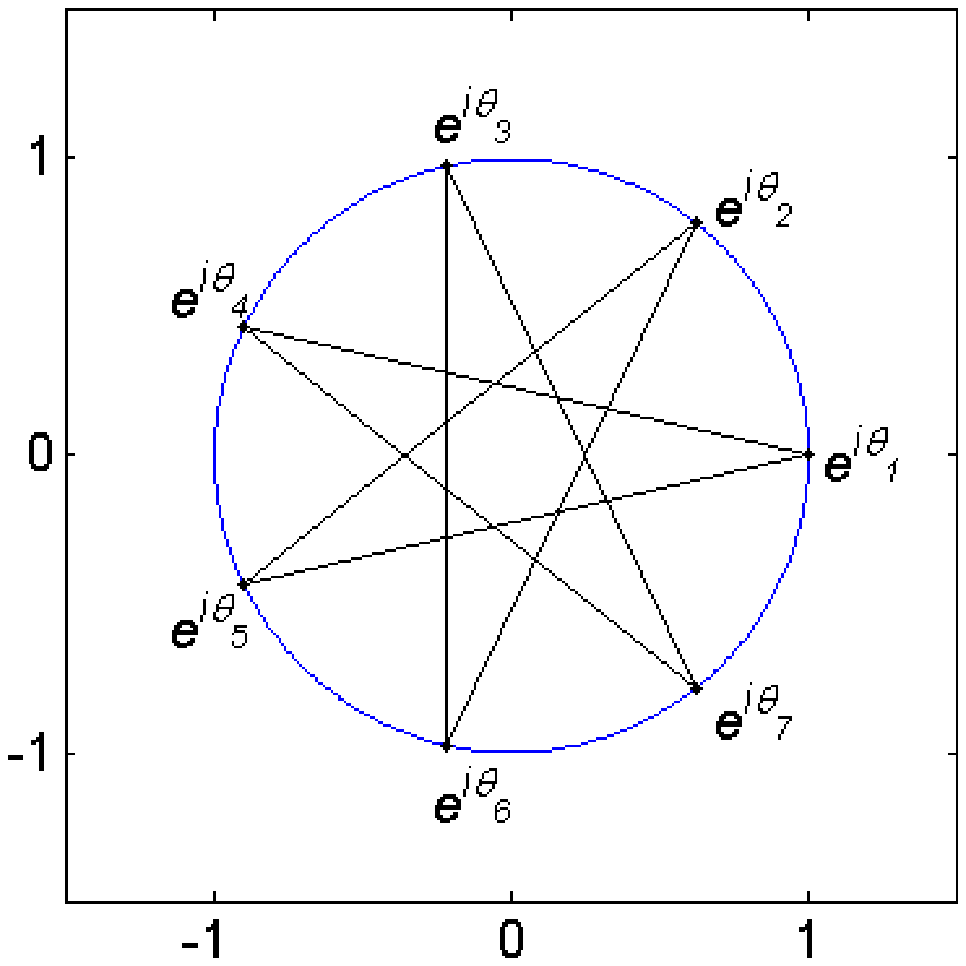}\hspace*{-25pt}
\caption{$\{d;((d-1)/2)\}$ star polygons with $d=5,7$ vertices.}\label{figStars57}
\end{center}
\end{figure}

For even $d$, the maximum is easy to find: $d/2$ points are equal to 1 while the others are $-1$, and the
permutation consists of $d/2$ 2-cycles. The maximum length is therefore $2d$.
See Figure \ref{figStars23} for an illustration. The odd case is not that simple.
Whereas $d=3$ can still easily be seen, larger sizes are more difficult.
It turns out that the maximal length is obtained when $\pi$ is a cyclic permutation (so that we have
only one polygon), the points are the $d$-th roots of unity, and they are connected in the shape of a star polygon
with Schl{\"a}fli-symbol $\{d;((d-1)/2)\}$ (see Figure \ref{figStars57}).
The upshot is that there are $d$ edges, and every edge has the same length $|1-e^{(d-1)\pi i/d}| = |1+e^{i\pi/d}|$.

We will prove this in two steps. First we calculate the maximal circumference $L(n)$ of a single $n$-polygon
(corresponding to $\pi$ being cyclic). Second, we show that $L(n)$ is superadditive: $L(\sum_i n_i)\ge\sum_i L(n_i)$.
That is, the total circumference does not increase by using a permutation $\pi$ consisting of several,
shorter cycles.

\bigskip

We first maximise $L(n)=\sum_{j=1}^n |e^{i\theta_{\pi(j)}}-e^{i\theta_j}|$ over all angles $\theta_j$
for $\pi$ a cyclic permutation.
As we can relabel the angles, it does not matter which cyclic permutation we take.
It therefore suffices to maximise
$\sum_{j=1}^n |e^{i\theta_{j-1}}-e^{i\theta_{j}}|$, with $\theta_0:=\theta_n$,
which is equal to $\sum_{j=1}^n |1-e^{i(\theta_{j}-\theta_{j-1})}|$.
Define $x_j=\theta_{j}-\theta_{j-1} \pmod{2\pi}$, so that $0\le x_j< 2\pi$.
We can now replace the maximisation over the angles by a maximisation over their differences $x_j$,
with the condition that $\sum_{j=1}^n x_j$ should be an integer multiple of $2\pi$ (because of the cyclicity
of $\pi$).
Noting also that $|1-e^{ix}| = \sqrt{2-2\cos x}$,
which in turn is equal to $2 \sin(x/2)$ over the interval $0\le x\le 2\pi$,
we then have the constrained maximisation
$$
L(n)=\max_{0\le x_1,\ldots,x_n\le 2\pi \atop
k\in \N, 0\le k\le n}
\left\{2\sum_{j=1}^n \sin(x_j/2): \sum_j x_j=2k\pi\right\}.
$$

From the concavity of the sine function over the interval $[0,\pi]$, we get
$$
\frac{1}{n} \sum_{j=1}^n \sin(x_j/2)
\le \sin\left(\frac{1}{2n} \sum_{j=1}^n x_j\right)
= \sin(k\pi/n),
$$
with equality if all $x_j$ are equal.
Therefore, the maximisation over the $x_j$ is readily done, and we get
$$
L(n) = \max_{k\in \N, 0\le k\le n} 2n\sin(k\pi/n).
$$
The remaining maximisation over $k$ is also easy:
for even $n$, $L(n)=2n$ (with $k=n/2$),
while for odd $n$ we get the smaller value
$$
L(n)=2n\sin((n-1)\pi/2n) = n\sqrt{2+2\cos(\pi/n)}.
$$

\bigskip

It remains to prove superadditivity of $L(n)$, i.e.
\[L(n)\ge L(n-k)+L(k).\]
For even $n$, this is simple:
either $k$ and $n-k$ are both even, in which case $L(n-k)+L(k)=2(n-k)+2k = 2n=L(n)$,
or they are both odd, in which case $L(n-k)+L(k)<2(n-k)+2k = 2n=L(n)$.

For odd $n$, consider $k$ odd and $n-k$ even, so that $L(n-k)=2(n-k)$.
We note that $L(n)-2n$ for odd $n$ is an increasing function of $n$.
Hence for odd $k$ and $n$, $L(n)-2n\ge L(k)-2k$, which directly implies
$L(n)\ge L(n-k)+L(k)$.
This ends the proof of a).

\bigskip

\noindent\textit{Proof of b).}
A variational characterisation of $C_{\infty,1,1}$ is
$$
C_{\infty,1,1} = \max_{X,Y}\{||\,XY-YX\,||_\infty: ||X||_1,||Y||_1\le1\}.
$$
Now note that $||\,XY-YX\,||_\infty$ is convex in $X$, and that the set of $X$ such that $||X||_1\le1$
is a convex set with extremal points the rank 1 matrices $X=uv^*$, where $u$ and $v$ are normalised vectors.
Thus,
\[\max_X \{||\,XY-YX\,||_\infty: ||X||_1\le1\}\]
is achieved for $X$ of the form $X=uv^*$.

Similarly, the latter is a convex function in $Y$
(the pointwise maximum of two convex functions is again convex) and, therefore, is also maximal for $Y$ of the form
$Y=ab^*$, where $a$ and $b$ are normalised vectors.
Hence,
$$
C_{\infty,1,1} = \max_{u,v,a,b} ||\,uv^*ab^*-ab^*uv^*\,||_\infty.
$$
The norm itself also has a variational expression:
\[||A||_\infty = \max_{p,q}|p^* Aq|,\]
where
$p$ and $q$ are also normalised vectors.
We thus end up with a maximisation over 6 normalised vectors:
\beas
C_{\infty,1,1} &=& \max_{u,v,a,b,p,q} |p^*(uv^*ab^*-ab^*uv^*)q| \\
&=& \max_{u,v,a,b,p,q} | (p,u) (v,a) (b,q) - (p,a)(b,u)(v,q)|.
\eeas
It is in principle possible to perform this maximisation over each of the 6 vectors in turn,
but the calculations immediately become very long-winded.
A much better approach is to focus attention to the inner products directly.

W.l.o.g.\ we can restrict the values of all inner products to be real, which can be done
simply by considering real vectors only.
It is easily seen that $| (p,u) (v,a) (b,q) - (p,a)(b,u)(v,q)|$ cannot be made bigger
by allowing complex valued inner products.
Thus, let $(p,u)=\cos\alpha$, $(v,a)=\cos\beta$ and $(b,q)=\cos\gamma$,
and $(p,a)=\cos\delta$, $(b,u)=\cos\eta$ and $(v,q)=\cos\theta$.
The point to observe now is that of these angles exactly 5 can be chosen independently,
while the remaining one is then subject to an inequality, as illustrated here:
$$
v \stackrel{\beta}{\longleftrightarrow} a
\stackrel{\delta}{\longleftrightarrow} p
\stackrel{\alpha}{\longleftrightarrow} u
\stackrel{\eta}{\longleftrightarrow} b
\stackrel{\gamma}{\longleftrightarrow} q.
$$
In this example, $\theta$, the angle between $v$ and $q$, is restricted to be less than the sum of all other
angles (which is not a restriction if that sum is larger than $\pi$).
Thus we get
$$
C_{\infty,1,1} = \max_{\alpha,\beta,\gamma,\delta,\eta,\theta\ge0}
\{ | \cos\alpha\cos\beta\cos\gamma - \cos\delta\cos\eta\cos\theta|:
0\le \theta \le \alpha+\beta+\gamma+\delta+\eta \}.
$$

\begin{lemma}\label{lem1}
$$
\max_{\alpha,\beta: \alpha+\beta=x} \cos\alpha\cos\beta = \cos^2(x/2).
$$
and
$$
\min_{\alpha,\beta: \alpha+\beta=x} \cos\alpha\cos\beta = -\sin^2(x/2).
$$

\end{lemma}
\textit{Proof of Lemma \ref{lem1}.}
\beas
\max_{\alpha,\beta: \alpha+\beta=x} \cos\alpha\cos\beta
&=& \max_\alpha \cos\alpha\cos(x-\alpha) \\
&=& (\cos x)/2 +\max_\alpha \cos(x-2\alpha)/2 \\
&=& (\cos x+1)/2 = \cos^2(x/2).
\eeas
The minimum is given by $(\cos x-1)/2=-\sin^2(x/2)$.
\qed

\begin{lemma}\label{lem2}
For $-\pi\le x\le \pi$,
$$
\max_{\alpha,\beta,\gamma: \alpha+\beta+\gamma=x} \cos\alpha\cos\beta\cos\gamma = \cos(x/3)^3;
$$
For $0\le x\le 2\pi$,
$$
\min_{\alpha,\beta,\gamma: \alpha+\beta+\gamma=x} \cos\alpha\cos\beta\cos\gamma = -\cos((x-\pi)/3)^3.
$$
The maximal and minimal values outside these intervals are obtained by
periodical extension.
\end{lemma}
\textit{Proof of Lemma \ref{lem2}.}
By applying Lemma \ref{lem1}, we get
\beas
\max_{\alpha,\beta,\gamma: \alpha+\beta+\gamma=x} \cos\alpha\cos\beta\cos\gamma
&=& \max_{y,\gamma: y+\gamma=x}  \cos\gamma \max_{\alpha,\beta: \alpha+\beta=y} \cos\alpha\cos\beta \\
&=& \max_{y,\gamma: y+\gamma=x}  \cos\gamma \cos(y/2)^2 \\
&=& \max_y \cos(x-y) \cos(y/2)^2.
\eeas
The stationary points of $\cos(x-y) \cos(y/2)^2$ as function of $y$ are $y=\pi$ and $y=2(x+k\pi)/3$,
yielding the values
$0$ and $\cos((x+k\pi)/3)^2 \cos[(x-2k\pi)/3]$.
The maximum of these values is $\cos(x/3)^3$ for $-\pi\le x\le \pi$, while
the maximum outside this interval is obtained by periodical extension.
The minimum is calculated in a similar way.
\qed

\medskip
With this lemma we are thus led to replace the maximisation over the 6 angles by
a single maximisation: we maximise the first term over angles $\alpha,\beta,\gamma$ subject to
$\alpha+\beta+\gamma=x$, and minimise the second term over angles $\delta,\eta,\theta$
subject to $\theta-\delta-\eta=x$ (as the sign of $\delta$ and $\eta$ is irrelevant
in $\cos\delta\cos\eta\cos\theta$ we can use Lemma \ref{lem2} here too). This leads to
$$
C_{\infty,1,1} = \max_{0\le x\le\pi} \cos^3(x/3) + \cos^3((x-\pi)/3).
$$
The maximum is achieved for $x=\pi/2$ and equal to
$\sqrt{27}/4$. This ends the proof of b).
\qed

\subsection{Interpolation revisited\label{sec4.3}}

The major hope behind Theorem \ref{thm1ii_i11} is of course that we might be able to use interpolation
to close the two gaps for the unknown triplets $(p,q,r)$.

\parpic[l]{\includegraphics[scale=\psclt]{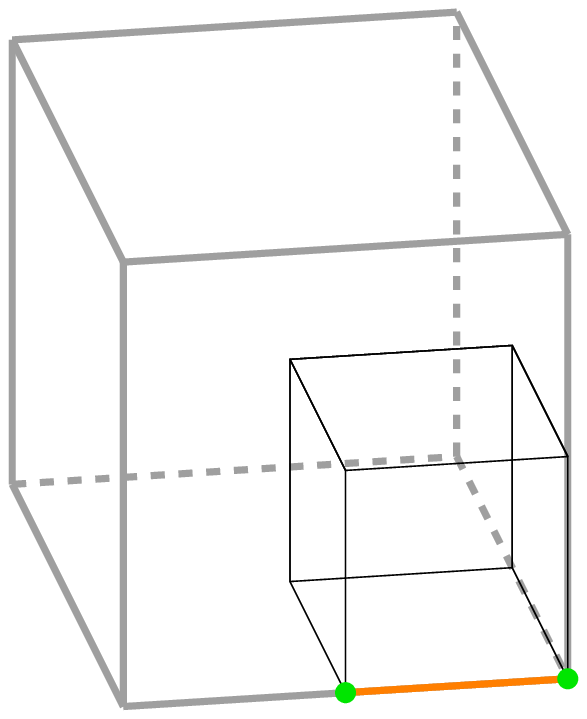}}\picskip{8}
For $p\geq 2, q,r\leq 2$ (the cube in the lower right of the illustrations) interpolation turns out
to be the wrong method, at least with the present data. We demonstrate this for the line $q=r=1$ $(p\geq 2)$.
Even the exact value of $C_{\infty,1,1}$ does not force interpolation bounds
to be sharp in this case.
We obtain
\begin{equation}\label{eqp11int}
C_{p,1,1}\leq 2^{1/p}\left(\frac{\sqrt{27}}{4}\right)^{1-2/p}.
\end{equation}

\medskip

However,
as in the proof of Theorem \ref{thm1ii_i11}~b) one can show that for these points
the $X$ and $Y$ achieving the maximum are matrices of rank one.
Hence, $XY-YX$ has at most two non-zero singular values.
Combining the knowledge of
\begin{center}
\begin{tabular}{lll}
$C_{2,1,1}$ &\,:\,& $\sqrt{\sigma_1^2+\sigma_2^2}\leq\sqrt{2}$\quad\quad and\\
$C_{\infty,1,1}$ &\,:\,& $\sigma_2\leq \sigma_1\leq\sqrt{27}/4$
\end{tabular}
\end{center}
already yields a better upper estimate than (\ref{eqp11int}), in fact
\begin{equation}\label{eqp112}
C_{p,1,1}\leq \left(\left(\sqrt{27}/4\right)^p + \sqrt{2-\left(\sqrt{27}/4\right)^2}^p \right)^{1/p}.
\end{equation}

\medskip

Recalling the example (\ref{eqspec1}) given in Section \ref{sec2.1} we may ensure
\begin{equation}\label{eqp11triv}
C_{p,1,1}\geq 2^{1/p}.
\end{equation}
Alas, this is a worse lower bound for $p\rightarrow\infty$ as it tends to $1<\sqrt{27}/4\approx 1.229$.

\medskip

The trickier example of two (normed) rank one matrices from \cite[page 1880]{BW2} gives us the curve
\begin{equation}\label{eqR1ExSv}
(\sigma_1,\sigma_2) = \sqrt{2}\frac{\sqrt{8\cos\phi\,\sin\phi}}{1+2\,\cos\phi\,\sin\phi}\,(\cos\phi,\sin\phi)
\quad {\rm with}\quad \phi \in [0,\pi/4]
\end{equation}
for possible singular values of $XY-YX$.

\begin{figure}[p]
\centering\includegraphics[scale=0.75]{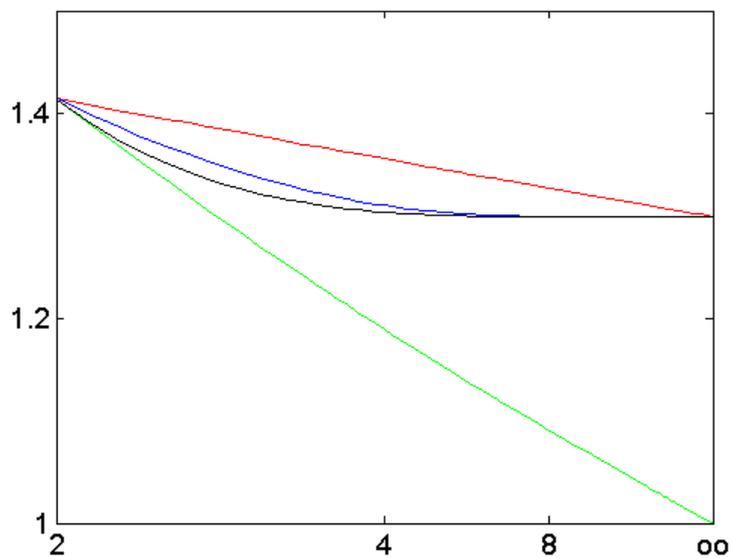}

\caption{Estimates for $C_{p,1,1}$: upper bounds (\ref{eqp11int}) (red) and (\ref{eqp112}) (blue) and lower bounds
(\ref{eqp11triv}) (green) and (\ref{eqR1ExSv}) (black).}\label{figp11Bounds}
\end{figure}

By choosing a point on the curve (\ref{eqR1ExSv}) with $p$-norm as large as possible
we obtain a very good lower bound to $C_{p,1,1}$, which is numerically approximated
in Figure~\ref{figp11Bounds}. Moreover, we conjecture that the resulting value is
equal to the constant $C_{p,1,1}$ for $p>2$. The estimates given by the upper and lower
bounds (also pictured in the figure) are already very tight.

\medskip

Note that $C_{\infty,q,1}$ and $C_{\infty,1,r}$ can be determined by duality from $C_{p,1,1}$.
Recall the symmetries of Proposition \ref{propCpqrsym} for that purpose.

\bigskip

\parpic[l]{\includegraphics[scale=\psclt]{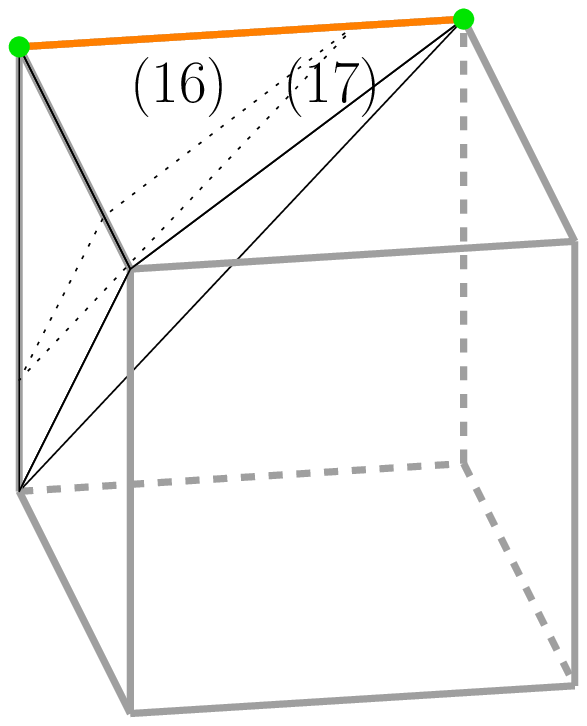}}\picskip{7}
Similarly, for odd-sized $d\times d$ matrices in the upper left pyramid there seems to be one more
plane of cusps determined by the example used for $C_{1,\infty,\infty}$ in Theorem~\ref{thm1ii_i11}~a)
\[X=\left(\begin{array}{cc}
 \cO & I_{\lceil d/2\rceil} \\
 I_{\lfloor d/2\rfloor} & \cO
 \end{array}\right), Y=\left(\begin{array}{ccc}
 e^{i\theta_1} & & \\
 & \ddots & \\
 & & e^{i\theta_d}
 \end{array}\right)
 \]
yielding the value
\begin{equation}\label{eqpii1}
C_{p,\infty,\infty}\geq d^{1/p}\sqrt{2+2\cos(\pi/d)}
\end{equation}
on the one hand, as well as the value
\begin{equation}\label{eqpii2}
C_{p,\infty,\infty}\geq 2 (d-1)^{1/p}
\end{equation}
on the other hand, given by padding an example matrix of even size $d-1$ with a zero line and column
\[X=\left(\begin{array}{cccccc}
 0 & 1 & & & & \\
 1 & 0 & & & & \\
  & & \ddots & & & \\
  & & & 0 & 1 & \\
  & & & 1 & 0 & \\
  & & & & & 0
 \end{array}\right), Y=\left(\begin{array}{cccccc}
 1 & & & & & \\
  & -1 & & & & \\
  & & \ddots & & & \\
  & & & 1 & & \\
  & & & & -1 & \\
  & & & & & 0
 \end{array}\right).\]

\begin{figure}[p]
\centering\includegraphics[scale=0.4]{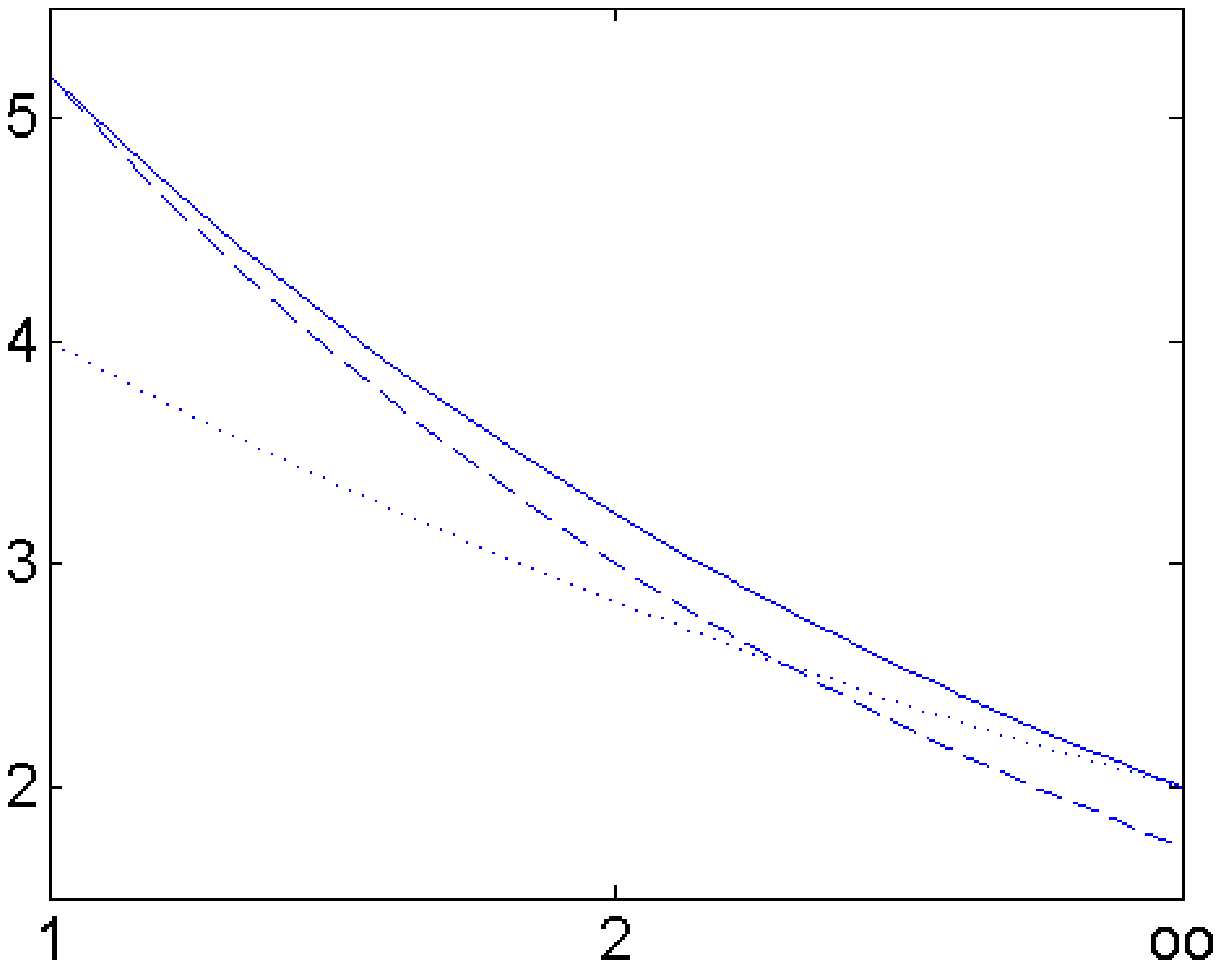} \includegraphics[scale=0.4]{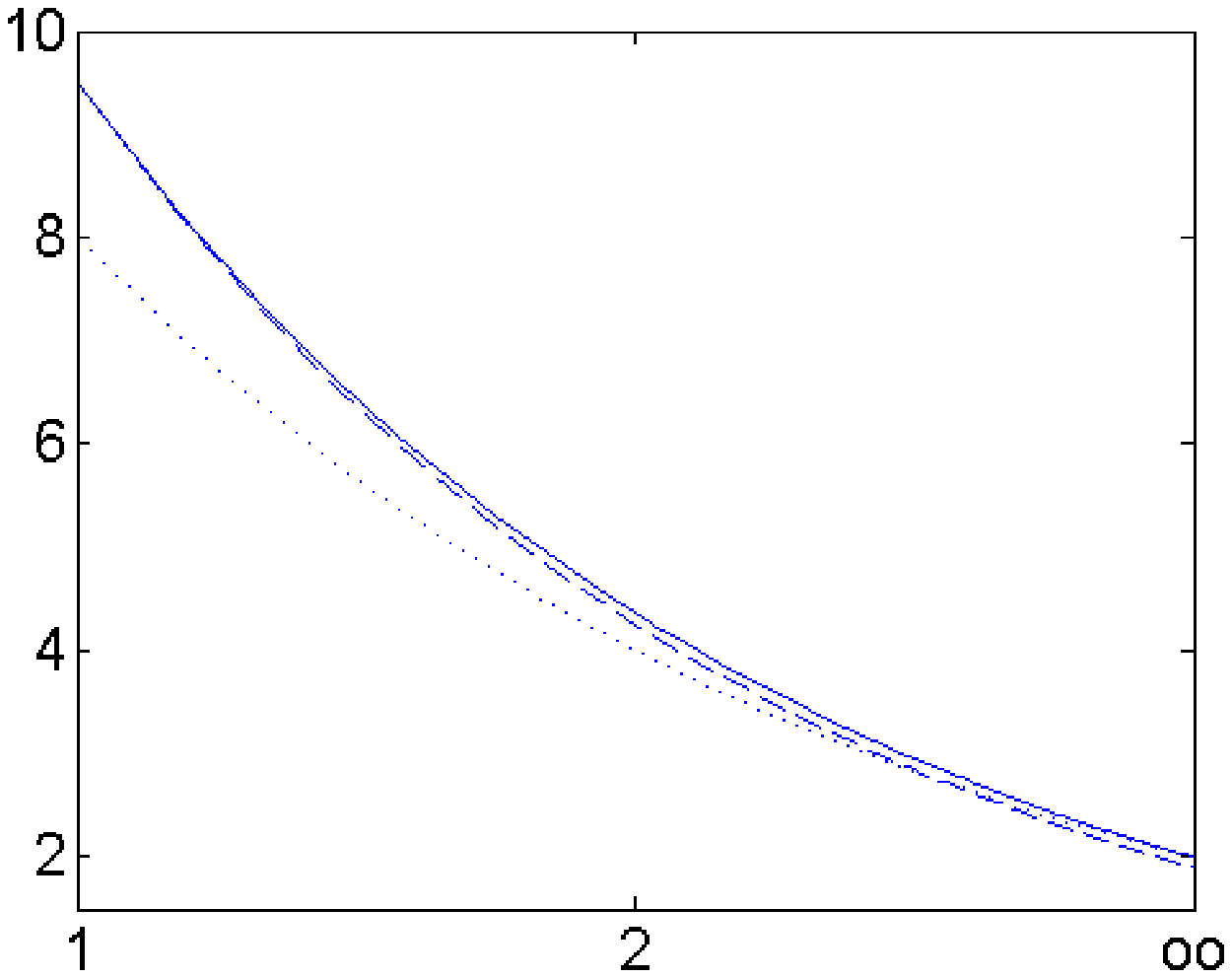}

\caption{Estimates for $C_{p,\infty,\infty}$: upper bound (\ref{eqpiiint}) (solid)
and lower bounds (\ref{eqpii1}) (dashed) and (\ref{eqpii2}) (dotted)
for $d=3$ (left) and $d=5$ (right).}\label{figpiiBounds}
\end{figure}

In the picture we indicated the areas of the pyramid where the examples yielding (\ref{eqpii1}) or
(\ref{eqpii2}) represent the largest known lower bounds.
Assuming that these two values are equal to $C_{p,\infty,\infty}$ and that
the two examples achieve the pyramid's values, it is left as a simple exercise
in interpolation to show that the interface boundary surface is really mapped to
a plane by $\ig\!^3$.
The second example is indeed a direct sum of the $2\times 2$ matrices (\ref{eqspec3})
from Section \ref{sec3.3},
padded with an additional row and column of zeroes to get an odd dimension;
recall that for even dimension these matrices did achieve equality.

The interpolation bound
\begin{equation}\label{eqpiiint}
C_{p,\infty,\infty}\leq \left(d \sqrt{2+2\cos(\pi/d)}\right)^{1/p}\cdot 2^{1-1/p}
\end{equation}
is likely not sharp.

From Figure \ref{figpiiBounds} we are tempted to conjecture that the difference between the bounds (\ref{eqpiiint})
and $\max\{(\ref{eqpii1}),(\ref{eqpii2})\}$ vanishes as $d\rightarrow\infty$.
Moreover, the index $p_0$ for which (\ref{eqpii1}) and (\ref{eqpii2}) coincide
seems to tend to infinity when the size $d$ is increased.

\section{Maximality\label{sec5}}

In \cite{BW2}, after giving the first general proof for $C_{2,2,2}=\sqrt{2}$, the notion
of maximality was introduced. This notion was subsequently extended to $p$-maximality in \cite{W}.
Consistent with these definitions we want to investigate the maximality problem in the general
context and call a pair $(X,Y)$ of $d\times d$ matrices $(p,q,r)$-maximal if both $X$ and $Y$
are non-zero and satisfy (\ref{eqNIpqr}) with equality, i.e.\
\[\|XY-YX\|_p = C_{p,q,r} \|X\|_q\|Y\|_r.\]
In contrast to \cite{W}, we are only looking here at the Schatten norms.

A characterization of $(2,2,2)$-maximality, which is called \textit{maximality} in \cite{BW2}
and \textit{Schatten 2-maximality} in \cite{W}, was recently given in \cite{CVW}.
This result will serve as a basis for further investigations to derive criteria for maximality
in the $(p,q,r)$ case, in combination with the tools we have used in Sections \ref{sec2} and \ref{sec3}
to obtain the exact values of the bound $C_{p,q,r}$.

\medskip

First of all, we will see that the method of monotonicity imposes strong restrictions.

\begin{lemma}\hfill\,\label{lemMaxMon}
\begin{enumerate}
\item[a)] If $C_{\tilde{p},q,r}=C_{p,q,r}$ was obtained by monotonicity via increasing $p<\tilde{p}$ and $(X,Y)$ is
$(\tilde{p},q,r)$-maximal then $\rg(XY-YX)=1$ and $(X,Y)$ is $(p,q,r)$-maximal.

\item[b)] If $C_{p,\tilde{q},r}=C_{p,q,r}$ was obtained by monotonicity via decreasing $q>\tilde{q}$ and $(X,Y)$ is
$(p,\tilde{q},r)$-maximal then $\rg X=1$ and $(X,Y)$ is $(p,q,r)$-maximal.
An analogous statement is true for $r$ and $Y$.
\end{enumerate}
\end{lemma}

\noindent\textit{Proof.}
The monotonicity argument in a) works as follows:
\[\frac{\|XY-YX\|_{\tilde{p}}}{\|X\|_q\|Y\|_r} \leq \frac{\|XY-YX\|_p}{\|X\|_q\|Y\|_r} \leq C_{p,q,r}.\]

Hence, if $(X,Y)$ is $(\tilde{p},q,r)$-maximal, the left-hand side equals $C_{\tilde{p},q,r}$
and this implies the $(p,q,r)$-maximality of the pair since all of the inequalities
in the chain become equalities. Moreover, we get
\[\|XY-YX\|_{\tilde{p}} = \|XY-YX\|_p\]
for $\tilde{p}>p$ which is only possible if the corresponding matrix has rank one.

The proof of b) is similar.
\qed

\medskip

In \cite{W} we argued that some properties are preserved by interpolation.
Furthermore, we used the fact that especially a rank one structure is left untouched.
Of course, this argument only works if the obtained interpolation bounds are sharp.

\begin{lemma}Let $1<p,q,r<\infty$.\label{lemMaxInt}

\begin{enumerate}
\item[a)] If $C_{p,q,r}$ is obtained by interpolation connected to a base, for which all matrices $X$
of a maximal pair $(X,Y)$ admit rank one, then also the matrices $X$ of a $(p,q,r)$-maximal pair $(X,Y)$
must have rank one. Similar statements hold for $Y$ and $XY-YX$.

\item[b)] If $C_{p,q,r}$ is obtained by interpolation between any point and $(2,2,2)$
(directly or via several steps)
and $(X,Y)$ is $(p,q,r)$-maximal then $(X,Y)$ is $(2,2,2)$-maximal.

\item[c)] If $C_{p,q,r}$ is obtained by interpolation connected to a base, for which all matrices $X$
of a maximal pair $(X,Y)$ are unitarily similar to matrices of the type
$\left(\begin{array}{cc}
0 & x_{12} \\
x_{21} & 0
\end{array}\right) \oplus \cO$
with $|x_{12}|=|x_{21}|$, then also all matrices $X$ of a $(p,q,r)$-maximal pair $(X,Y)$
have this property. Similar statements hold for $Y$ and $XY-YX$.
\end{enumerate}
\end{lemma}

\noindent\textit{Proof.}
The key point in the proofs is that if $(X,Y)$ is a maximal pair with respect to an interpolated triplet,
then an appropriately modified pair $(\tilde{X},\tilde{Y})$ is maximal with respect to the
base point triplet.

An analysis in \cite[proof of Proposition 8]{W} showed that the matrix $\tilde{X}$ is actually a scaled version
of $X$ in the sense that every entry (i.e.\ a complex number) keeps its complex argument,
but has its absolute value raised to a specific power
(one of us calls this operation a \textit{polar power}; see \cite{kapolar}).
More precisely, if
$x_{jk}=re^{i\varphi}$ then $\tilde{x}_{jk}=r^Pe^{i\varphi}$.

Clearly, an entry with the value 0 is not altered in any way by this procedure.
Moreover, the claim of a) was already proven true and applied in \cite{W} based on these ideas.

\medskip

Now, for any interpolation connected to the base point $(2,2,2)$ or any other point that has been obtained
by such a process, the scaled pair needs to be maximal in the original sense.
This statement is true if the interpolation process is the usual Riesz-Thorin theorem (complex version)
or the tensor argument extension, since the tensor structure is unharmed
by the scaling procedure.

By Theorems 3.1 and 3.2 of \cite{CVW} all these pairs are given by:
\[U\tilde{X}U^*=\tilde{X}_0\oplus\cO,\quad U\tilde{Y}U^*=\tilde{Y}_0\oplus\cO\]
with $\tilde{X}_0,\tilde{Y}_0\in\C^{2\times 2}$ and
\[0=\tr \tilde{X}_0=\tr \tilde{Y}_0=\tr \tilde{Y}_0^*\tilde{X}_0.\]

The only information that we had obtained in \cite{BW2} about
matrices of a maximal pair was that they should have rank at most two,
which was not enough to obtain meaningful restrictions for interpolants in \cite{W}. But
with the simultaneous unitary similarity to essentially $2\times 2$ matrices,
it is now easy to see that with $\tr \tilde{X}_0=0$ also $\tr X_0=0$ is given.
The last conclusion is only possible as the trace is now the sum of only two
entries, or equivalently we have the relation
\[\tr \tilde{X}_0=0 \Leftrightarrow \tilde{x}_{11}=-\tilde{x}_{22}\]
which is kept by scaling the modulus back to $X$.

For transferring the orthogonality of $\tilde{X}_0$ and $\tilde{Y}_0$ to $X_0$
and $Y_0$ we furthermore
need the well-known statement that a trace zero matrix is unitarily similar to a matrix whose diagonal elements
are all zero (\cite{HJI}, p.\ 77). Hence, without loss of generality we may assume
$\tilde{X}_0=\left(\begin{array}{cc}
0 & \tilde{x}_{12} \\
\tilde{x}_{21} & 0
\end{array}\right)$, implying for the scalar product
\[\tr \tilde{Y}_0^*\tilde{X}_0 =0=\tilde{x}_{12}\overline{\tilde{y}_{12}}+\tilde{x}_{21}\overline{\tilde{y}_{21}}.\]
Of course, the latter is also preserved by scaling, since there are again only two summands
in which exactly one component of $X$ and one of $Y$ appear as factors.
Since $(X,Y)$ obeys the same relations as $(\tilde{X},\tilde{Y})$ specified above and
the theorems in \cite{CVW} yield necessary and sufficient conditions,
we obtain the $(2,2,2)$-maximality of the pair.

\medskip
The claim of c) can be shown in a similar but even simpler fashion.
\qed

\medskip

{\bf Remark.} The part c) in Lemma \ref{lemMaxInt} is a generalization of a) since
every rank one matrix with trace zero is unitarily similar to a matrix of the type
$\left(\begin{array}{cc}
0 & x_{12} \\
0 & 0
\end{array}\right) \oplus \cO$. The zero trace will automatically be given in combination
with b). Note that for the described type of matrices $X$ one has $\sigma(X)=(c,c,0,...)$
for some $c>0$. Hence, $X$ is unitarily similar to a multiple of a unitary $2\times 2$
matrix that is padded with zeros.

As all matrices of maximal pairs connected to b) admit rank not greater than two,
we are able to apply the strong estimate
\[\|A\|_p\leq\|A\|_q\leq 2^{1/q-1/p}\|A\|_p\quad\forall p\geq q\]
for Schatten norms. In general, the constant 2 in the second inequality would have been
the rank or even the size $d$.
Such estimates were crucial in \cite{W} to determine $(1,1,1)$-maximal pairs
and will also be of use in the following.

\medskip
Both Lemmas imply that maximal pairs can only be found in a very limited range.
We will see that only the boundary of the parameter space may need a separate treatment,
but will mostly fit with the results of the interior. Lemma \ref{lemMaxInt} b) implies that moreover
$(2,2,2)$-maximality can be expected to be richer than others (excluding possibly cases
like $(\infty,\infty,\infty)$ at the boundary).
Before proceeding with the consequences of these two results we need to introduce
a new drawing convention we'll adhere to.

\medskip

Up to now we had no problems to picture sets of points $(p,q,r)$, as all of them were closed sets,
i.e.\ points, lines with end-points or complete bodies containing all of its bounding facets.
However, for visualizing areas connected to shared properties of maximality we will encounter
open sets. In order to visualize them in a comprehensible way we only draw lines and points
instead of colored facets, and in the following way:

\medskip

\begin{minipage}[t]{\pwdt}
\includegraphics[scale=\pscl]{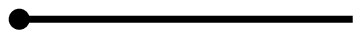}\vspace*{8pt}

\includegraphics[scale=\pscl]{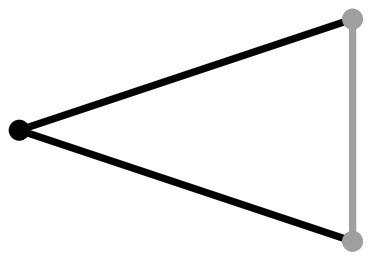}\vspace*{8pt}

\includegraphics[scale=\pscl]{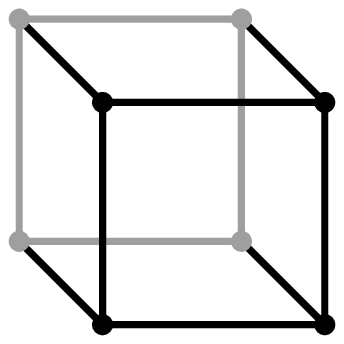}\vspace*{16pt}

\includegraphics[scale=\pscl]{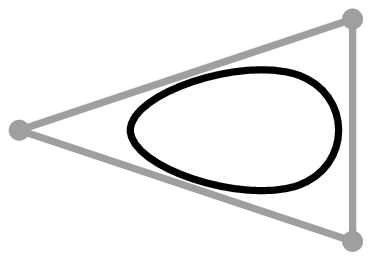}\vspace*{-5pt}

\includegraphics[scale=\pscl]{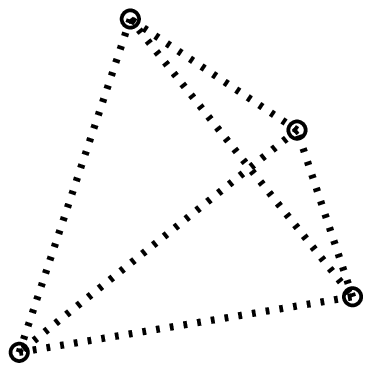}
\end{minipage}\begin{minipage}[t]{\twdt}
This marks all points on the line except the right end.

\medskip
This picture marks all points of the triangle, excluding the grey edge at the right.
The grey line itself contains its end-points.

\bigskip
In three-dimensional space this marks the complete body enclosed by the facets
of the same colour and their neighbours. If a line or facet is colored differently
it is excluded. For instance, the image on the left marks, in black, the whole cube except
the back facet and its boundaries.

\bigskip
The oval marks the interior of a facet, i.e.\ excluding the grey boundary.

\bigskip
Dotted lines mark the interior
of the three-dimensional body, i.e.\ excluding its surface (facets, edges and vertices).
\end{minipage}

\begin{theorem}\label{thmCpqrEq}
For $(p,q,r)$ in the respective areas, one has:
A pair $(X,Y)$ with $Z=XY-YX$ is $(p,q,r)$-maximal if and only if
there exist a unitary $U\in\C^{d\times d}$ and $X_0,Y_0,Z_0\in\C^{2\times 2}$
with $\tr X_0=\tr Y_0=0=\tr(Y_0^*X_0)$ such that
\[UXU^*=X_0\oplus\cO,\quad UYU^*=Y_0\oplus\cO,\quad UZU^*=Z_0\oplus\cO\]

and moreover:

\smallskip

\includegraphics[scale=0.5]{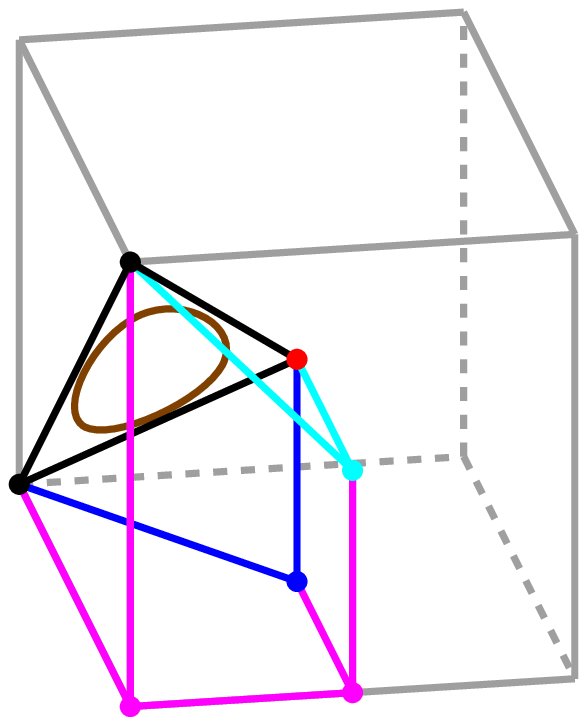}
\begin{minipage}[b]{30pt}
{\color[rgb]{1,0,0}1)}

{\color[rgb]{0,1,1}2a)}

{\color[rgb]{0,0,1}2b)}

{\color[rgb]{1,0,1}3a)}

{\color[rgb]{0.5,0.25,0}4c)}
\end{minipage}
\includegraphics[scale=0.5]{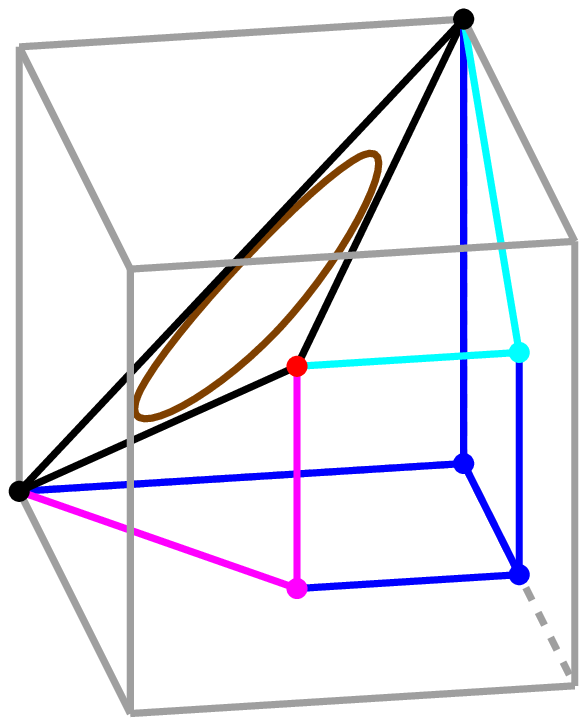}
\begin{minipage}[b]{30pt}
{\color[rgb]{1,0,0}1)}

{\color[rgb]{1,0,1}2b)}

{\color[rgb]{0,1,1}2c)}

{\color[rgb]{0,0,1}3b)}

{\color[rgb]{0.5,0.25,0}4a)}
\end{minipage}
\includegraphics[scale=0.5]{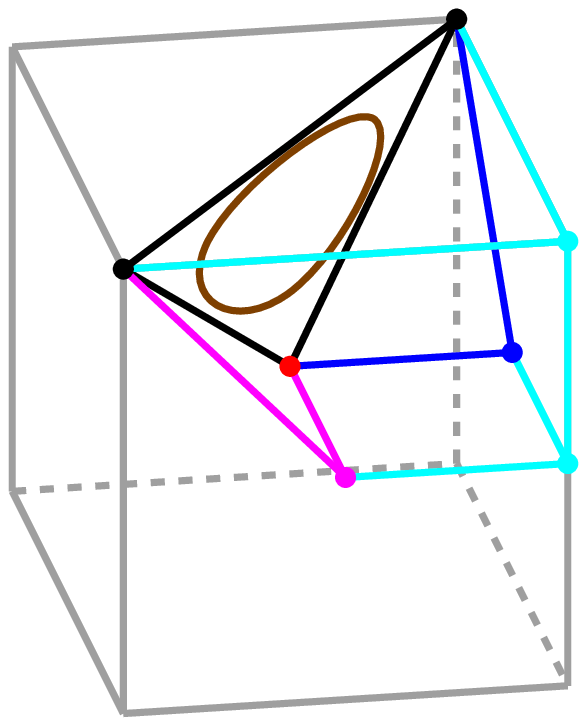}
\begin{minipage}[b]{30pt}
{\color[rgb]{1,0,0}1)}

{\color[rgb]{1,0,1}2a)}

{\color[rgb]{0,0,1}2c)}

{\color[rgb]{0,1,1}3c}

{\color[rgb]{0.5,0.25,0}4b)}
\end{minipage}

\begin{enumerate}
\item[1)] $X_0,Y_0$ are arbitrary otherwise,
\end{enumerate}

\begin{minipage}{0.5\textwidth}
\begin{enumerate}
\item[2a)] $\rg X_0=1$,

\item[2b)] $\rg Y_0=1$,

\item[2c)] $\rg Z_0=1$,
\end{enumerate}
\end{minipage}\begin{minipage}{0.5\textwidth}
\begin{enumerate}
\item[3a)] $\rg X_0=\rg Y_0=1$,

\item[3b)] $\rg Y_0=\rg Z_0=1$,

\item[3c)] $\rg X_0=\rg Z_0=1$,
\end{enumerate}
\end{minipage}

\begin{minipage}[b]{0.67\textwidth}
\begin{enumerate}
\item[4a)] $X_0$ is a non-zero multiple of a unitary matrix,

\item[4b)] $Y_0$ is a non-zero multiple of a unitary matrix,

\item[4c)] $Z_0$ is a non-zero multiple of a unitary matrix,\\

\item[5)] $X_0,Y_0$ and $Z_0$ are all a non-trivial multiples of unitary matrices.
\end{enumerate}
\end{minipage}\hspace*{15pt}\includegraphics[scale=0.5]{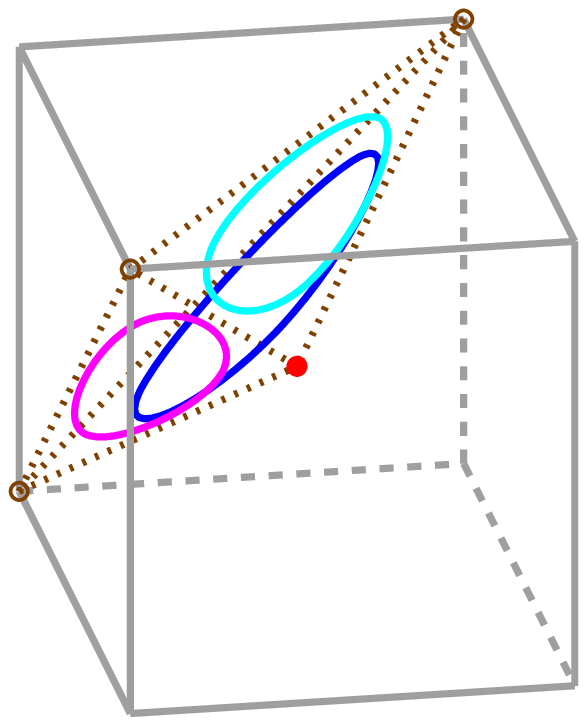}
\begin{minipage}[b]{30pt}
{\color[rgb]{1,0,0}1)}

{\color[rgb]{0,1,1}4b)}

{\color[rgb]{0,0,1}4a)}

{\color[rgb]{1,0,1}4c)}

{\color[rgb]{0.5,0.25,0}5)}
\end{minipage}
\end{theorem}

A look at Theorem \ref{thmCpqr4} should make clear why we won't describe
the regions of parameter space by (in)equalities at this point.

\medskip

\noindent\textit{Proof.}

\begin{minipage}[b]{\pwdt}
\includegraphics[scale=\psclt]{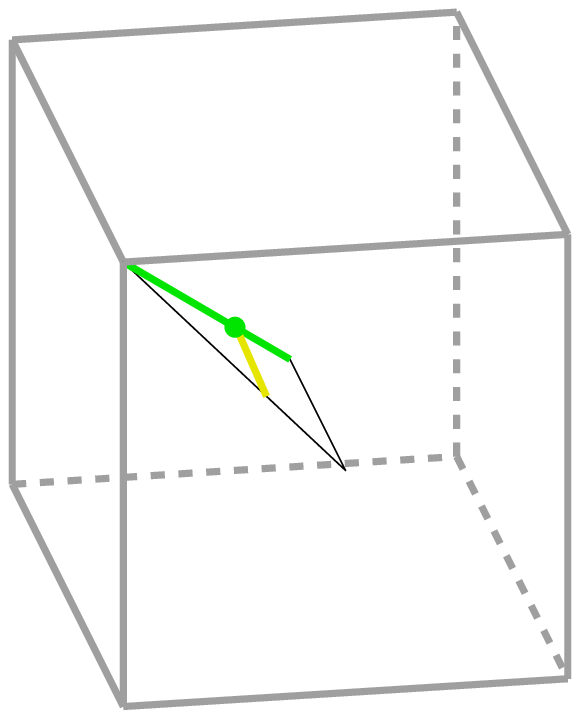}

\includegraphics[scale=\psclt]{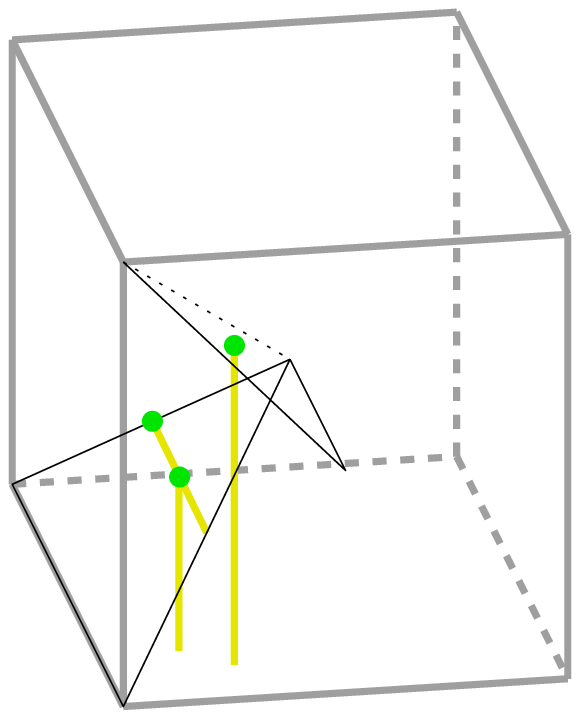}
\end{minipage}\begin{minipage}[b]{\twdt}
The origin of the investigations is {\it 1)} and has been proven in \cite{CVW}
as Theorems 3.1 and 3.2.

\medskip
The parts {\it 2)} are results of Lemma \ref{lemMaxMon}. For this recall
the construction of the values in these areas by monotonicity
in Section \ref{sec3.1}. For instance, in the case {\it 2a)} $q$ was decreased.
As a consequence $X$ must be a rank one matrix.

\medskip

The parts {\it 3)} are similarly easy. But, beginning with a triangle from {\it 2)},
a second monotonicity along another direction is applicable.

\medskip

As pictured for {\it 3a)} we decrease $r$, yielding $\rg Y=1$ additionally.
Observe that, except for $p=1$, Lemma \ref{lemMaxInt} and Lemma \ref{lemMaxMon}
(for transferring to $q=1$ and $r=1$) grant the similarity to $2\times 2$ matrices.

\medskip
For the last segment of this area remember that we closed the gap by interpolation with one base point
in the triangle at $p=1$.
\end{minipage}

\smallskip
\begin{minipage}[b]{\pwdt}
\includegraphics[scale=\psclt]{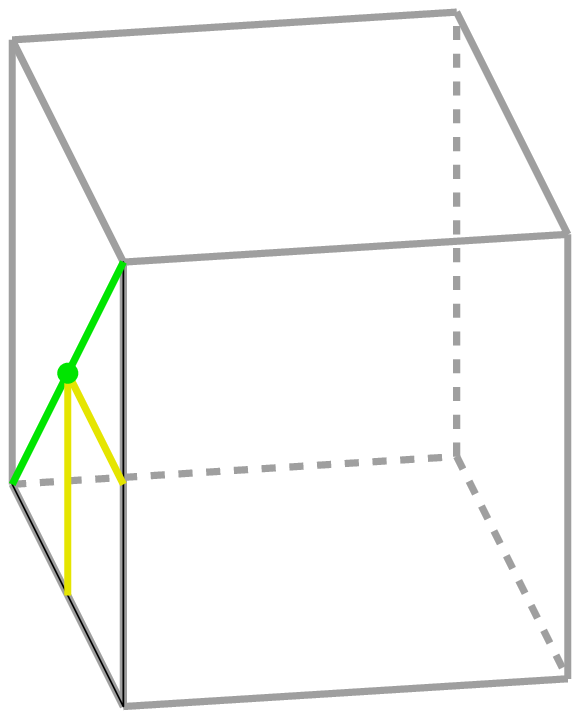}
\end{minipage}\begin{minipage}[b]{\twdt}
These base points can be handled by monotonicity with two
different directions, granting $\rg X=\rg Y=1$. The similarity relation is now simply verified directly.
Then by Lemma \ref{lemMaxInt} all interpolants inherit this property.

\medskip

The other two parts {\it 3b)} and {\it 3c)} may be proven in a similar fashion or can be shown
by a symmetry argument.
Take a look at the proof of Proposition \ref{propCpqrsym} and observe that properties of
maximal pairs are indeed swapped between $X$, $Y$ and $XY-YX$ as stated.
\end{minipage}

\begin{minipage}[b]{\pwdt}
\includegraphics[scale=\psclt]{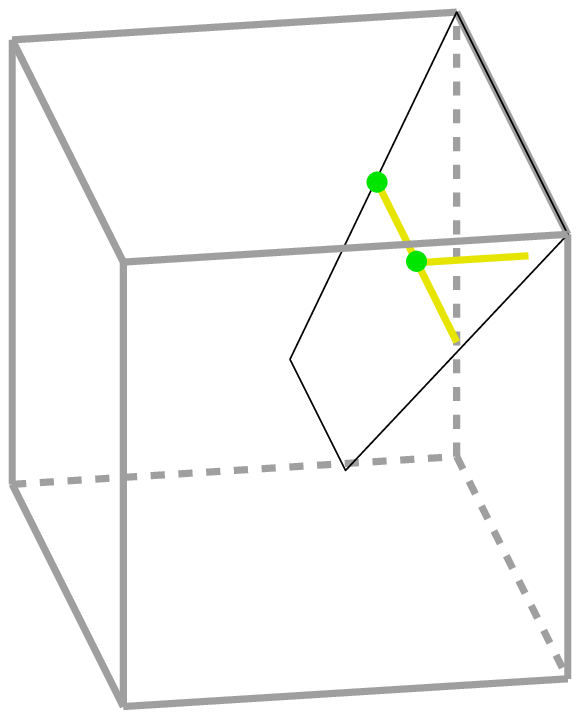}

\includegraphics[scale=\psclt]{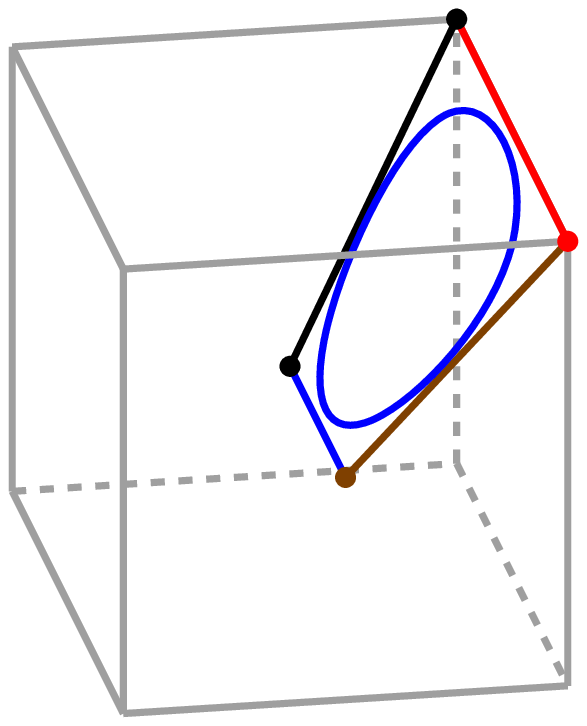}
\end{minipage}\begin{minipage}[b]{\twdt}
A small hint for direct use with {\it 3c)}:

\medskip

In one half of the area it is easy to see that $\rg X=1$ and $\rg Z=1$
with help of monotonicity.

\medskip

In the blue area we know for sure the $(2,2,2)$-maximality. Hence,
\[\sigma(Y)=(y_1,y_2,0,...),\quad\sigma(Z)=(z_1,z_2,0,...)\]
and because of $\rg X=1$ (due to monotonicity), one may assume for maximal pairs
\[2^{1-1/r}=\frac{(z_1^r+z_2^r)^{1/r}}{(y_1^r+y_2^r)^{1/r}} \leq
\frac{\sqrt{z_1^2+z_2^2}}{2^{1/r-1/2}\sqrt{y_1^2+y_2^2}}=2^{1-1/r}\]
yielding $\|Z\|_r=\|Z\|_2$ or equivalently $\rg Z=1$. At
this point the rank-specific estimates from the last remark\linebreak\vspace*{-10pt}
\end{minipage}

came into play.

The brown line inherits the properties of the blue area thanks to monotonicity, by Lemma \ref{lemMaxMon}.
For the red line, a simple calculation gives
\[\|Z\|_1\leq 2\|X\|_1\|Y\|_\infty=2\|X\|_q\|Y\|_\infty=\|Z\|_\infty,\]
which also results in $\rg Z=1$. In this case, $\rg X=1$ was again the result
of Lemma \ref{lemMaxMon}.

Now, having for the facet two rank one matrices in any maximal pair,
this is also true for the second half of the {\it 3c)} area by Lemma \ref{lemMaxInt}.

\medskip

For {\it 4c)}, i.e.\ the triplets $(p,q,q')$, check due to the $(q,q,q')$-maximality
and subsequently the $(2,2,2)$-maximality
in one half of the triangle
that $Z$ has only two non-zero singular values. Hence, we can write
\[2^{1/p}=\frac{\|Z\|_p}{\|X\|_q\|Y\|_{q'}} \leq \frac{2^{1/p-1/q}\|Z\|_q}{\|X\|_q\|Y\|_{q'}}\leq 2^{1/p}\]
yielding
\[\|Z\|_p=2^{1/p-1/q}\|Z\|_q\]
which results in $\sigma(Z)=(c,c,0,...)$ for some $c>0$ as claimed.
In the second half of the triangle we have $(q',q,q')$-maximality and the conclusions are
analogous.

\medskip

\parpic[l]{\includegraphics[scale=\psclt]{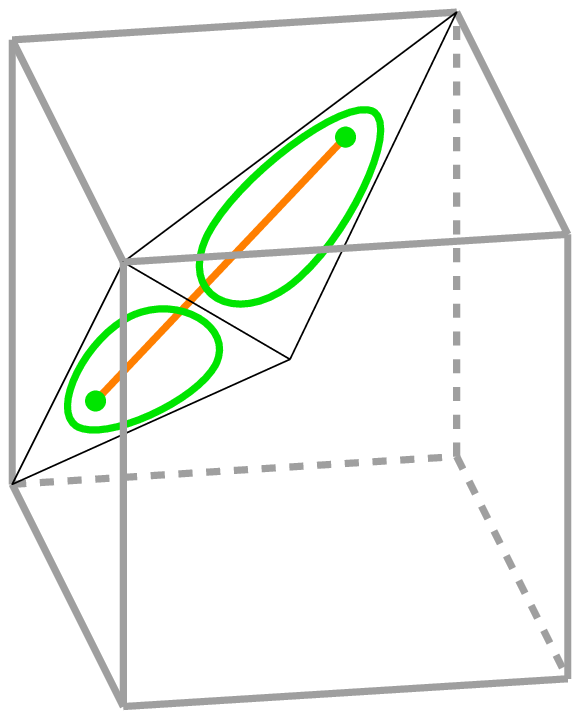}}\picskip{12}
The area for {\it 5)} can be handled by interpolation. However, here we will do this
in a different way than in Section \ref{sec3}, by interpolating two indices simultaneously
($p$ and $r$ or $p$ and $q$) and keeping the other index fixed.
Here the ordinary interpolation (for $K_X$ or $K_Y$) already suffices,
in line with what had already been observed in Section \ref{sec2.3}. \\
The picture illustrates one such
process for fixed $X$ and $q$. By this process, the unitarity properties of $Y_0$ and $Z_0$ are inherited
by the interpolants, too.

\vspace*{30pt}

Another interpolation process (fixing $Y$ and $r$) and the first one
complement each other in a fitting manner, determining also the properties of $X_0$ and $Z_0$.
In the end, since every point of the interior is an interpolant with respect to two interpolations,
all three matrices are of the asserted type.
\qed

\medskip
We remark that a property of {\it 3)} automatically implies the respective property
of {\it 4)} for the complementing variable, e.g. {\it 3a)} $\Rightarrow$ {\it 4c)},
but the converse is not true. Notice that one bounding facet is not covered by the
methods of the proof.


\section{Conclusions\label{sec6}}

In accordance with the title we have chosen for this paper, we want to point out
several occasions at which the `river of convexity' crossed our way.
First of all, we perused a specialized convexity theorem
in the form of Riesz-Thorin interpolation in Sections \ref{sec2} and \ref{sec3}.
We have seen that in many cases this theorem, in its usual form, does an excellent job.
Even for the bilinear operator called commutator it becomes applicable by fixing variables.
One major issue could be efficiently solved by applying this theorem to some unusual structures.
By interpolating along lines that (taken together) build up planes we establish new
bases for subsequent interpolation steps. In summary, these axis-oriented processes
are able to cover even more complicated regions of parameter space and may give strong estimates.

Furthermore, we demonstrated that it is possible to illustrate a bunch of (in)equalities and descriptive processes
in an intuitive way. We hope the reader enjoyed using this graphical tool
rather than having to comb through a vast array of formulas.

In Section \ref{sec4} we encountered convexity multiple times. We have seen convex
functions (and their concave counter-part), convex sets and properties related to both
of them with regard to extremal points. Here, we also have drawn connections to a visually appealing
geometrical problem.

\section*{Acknowledgments}
We are greatly indebted to Alexei Karlovich for bringing reference \cite{BS} to our attention  after
we presented the basic ideas of our approach during the WATIE 2009 conference.
When first coming in contact with interpolation
theory we quickly appreciated its power and applicability to our investigations. The careful
reader may have noticed at which point our basic approach first ran into problems, namely at the 'diagonal part'.
We were aware that for processing this specific segment a 'bilinear interpolation theorem' would be
extraordinarily helpful. By reading the title of Riesz' original paper it now becomes clear why our
search for such a particular result was doomed for failure from the very beginning. Our work-around for
compensating the lack of a reference has at least one consequence: besides the primary
interpretation as a theorem for multilinear operators, we may also look at the result as a theorem
for linear operators acting on sets that are not spaces.
We wanted to tell this story in order to
ensure bigger greed for the future.

\section*{Recommended Reading}

We divided the bibliographic section into two parts,
the actual references we required in the proofs and the listing of the papers
related to the topic that are interesting for obtaining further information
but not yet necessary to understand the present paper. This also gives
a historical overview on the developments.

\medskip
Origins of the commutator problem:
\begin{itemize}
\item A. B\"{o}ttcher, D. Wenzel,
How big can the commutator of two matrices be and how big is it typically?
Linear Algebra Appl. 403 (2005) 216--228.

{\it The paper that started the topic. We see an explanation for the observation
that it is hardly possible to find matrices with a big commutator.
Further $C_{2,2,2}=\sqrt{2}$ is verified when restricted to special classes
of matrices.}

\item L. L\'{a}szl\'{o}, Proof of B\"{o}ttcher and Wenzel's conjecture on commutator norms for 3-by-3 matrices,
Linear Algebra Appl. 422 (2007) 659--663.

{\it The first proof of the Schatten 2 problem for general matrices of size greater than 2.
Sadly the ideas seem not to be portable to larger matrices.}

\item S.-W. Vong, X.-Q. Jin,
Proof of B\"{o}ttcher and Wenzel's conjecture,
Oper. Matrices 2 (2008) 435--442.

{\it The first proof comprising all real $d\times d$ matrices
and a demonstration that the result can be shown elementary.
The paper is moreover an excellent example that `elementary' should not be confused with `short'
or `trivial'.}

\item Z.-Q. Lu, Proof of the normal scalar curvature conjecture,
Available from: arXiv:0711.3510v1 [math.DG] 22 November 2007.

{\it Another proof for real $d\times d$ matrices. This one is interesting
from an operator theoretic point of view and gave some ideas that influenced
the characterization of maximality found in \cite{CVW}.}

\end{itemize}

The development afterwards continued with \cite{BW2}, \cite{Aud}, \cite{CVW} and \cite{W}.

\medskip

The beginnings of interpolation theory:
\begin{itemize}
\item M. Riesz, Sur les maxima des formes bilin\'{e}aires et sur les fonctionnelles lin\'{e}aires,
Acta Math. \textbf{49}, 465--497 (1926). In french.

{\it In this paper a (from today's point of view) quite special convexity theorem is deduced.
It turned out to be only the first result in a line of generalisations.}

\item G.O.~Thorin, Convexity theorems generalizing those of M.~Riesz and
Hadamard with some applications,
Comm.\ Sem.\ Math.\ Univ.\ Lund [Medd.\ Lunds Univ.\ Mat.\ Sem.] \textbf{9}, 1--58 (1948).

{\it This paper gave rise to the whole mathematical sector of interpolation theory by
extending Riesz' result to the complex numbers, with an idea that rightly may be called ingenious.
In J.E.~Littlewood's words: it is one of the
most impudent ideas in mathematics.}

\item Pham The Lai, L'analogue dans ${\cal C}^p$ des th\'eor\`emes de convexit\'e
de M.~Riesz et G.O.~Thorin, Studia Math.\ \textbf{46}, 111--124 (1973). In french.

{\it We included a reference to this work as it demonstrates how Thorin's proof is adapted for
Schatten norms of (finite or infinite) matrices, which would be enough for our needs.}
\end{itemize}

\end{document}